\documentclass{amsart}
\usepackage{graphicx}
\usepackage[pdfauthor={Richard J. Mathar},pdfkeywords={Tatami mats, Domino tilings, combinatorics},pdfsubject={Mathematics, Combinatorics},colorlinks=true,citecolor=blue]{hyperref}

\usepackage{amsmath}

\newtheorem{thm}{Theorem}
\newtheorem{defn}{Definition}
\newtheorem{rem}{Remark}
\newtheorem{conj}{Conjecture}
\newtheorem{exa}{Example}

\begin{document}

\title[Paving Floors with Rectangular Tiles]{Paving Rectangular Regions with Rectangular Tiles: Tatami and non-Tatami Tilings.}

\author{Richard J. Mathar} 
\email{mathar@mpia.de}
\urladdr{http://www.mpia.de/~mathar}
\address{Hoeschstr. 7, 52372 Kreuzau, Germany}


\date{\today}

\begin{abstract}
The number of complete tilings of $m\times n$ floors for tiles
of shape $1\times 2$, $1\times 3$, $1\times 4$ and $2\times 3$
is computed numerically for floors up to width $m=9$ and variable
floor lengths $n$. Counts are obtained for two classes, for fixed
tile stack orientation on one hand and for counts up to rotations
and reflections on the other hand. Counts are refined by the
number of points on the floor where 4 tiles meet, i.e., by
the degree of violation of the requirement for Tatami tilings.
\end{abstract}

\maketitle 

\section{Definitions} 

\subsection{Rectangular floors and tiles}
We consider floors of rectangular dimension $m\times n$ laid out by
unit squares. We count full covers with tiles of width $t_m$ and length $t_n$.
The tiles may be of mixed orientation---of which there are two,
aligning $t_n$ with $n$ or with $m$.
Only one tile shape is considered at a time; hybrid covers by
dominos ($1\times 2$ tiles) and unit-tiles, for example, will not
be discussed.

An obvious requirement of each complete tiling of the floor
is that the ratio of the floor area by the tile area is an integer:
\begin{defn} (Matching Condition)
\begin{equation}
nm \equiv 0\quad \mod (t_nt_m).
\end{equation}
\end{defn}
Another appropriate constraint is that $t_n$ and $t_m$ are coprime, since
otherwise a trivial congruent scaling of the floor
geometry and dimensions by the largest common factor
would generate essentially the same results.

We write $T(n,m)$ for the counts of different tilings
of the $m \times n$ floor, but do not add labels that make the
dependence on the side lengths $t_n$ and $t_m$ of the tile explicit.

\subsection{Tatami condition}
Four tiles may meet at between $0$ and the maximum of
$\lfloor (n-1)/t_n\rfloor \lfloor(m-1)/t_m\rfloor$
and
$\lfloor (n-1)/t_m\rfloor \lfloor(m-1)/t_n\rfloor$
points at corners
of the unit base squares.
(This upper limit is poor in the sense that is an upper estimate assuming that
all tiles are aligned vertically or horizontally.)
At points of that type, the edges of the tiles form
a 4-way crossing.
\begin{defn} (Counts by number of 4 tile meets)
$T_t(n,m)$ denotes the number of tilings which contain $t$
points where 4 tiles meet.
\end{defn}
If we call tilings where no four tiles meet ``Tatami tilings,''
the number of Tatami tilings is $T_0(n,m)$.
The $T_t(n,m)$ with $t>0$ count tilings that violate the Tatami property.

This work shows tables of $T_t(n,m)$ where the floor width $m$ 
and the tile shape $t_n$ and $t_m$ are fixed, where the floor length $n$
increases along the table rows, and where $t$ increases along
the table columns (Tables \ref{tab.run2}--{\ref{tab.run9i_3_2}}).
The leftmost column in the tables shows $n$ whenever the matching condition
is fulfilled. The second column in the
tables are the row sums,
\begin{equation}
T(n,m) = \sum_{t\ge 0} T_t(n,m).
\label{eq.rsum}
\end{equation}

Ordinary generating functions downward columns are
noted as follows:
\begin{defn} (Generating function for tilings where $t$ points exists where 4
tiles meet)
\begin{equation}
T_t(z,m) \equiv \sum_{n\ge 0} T_t(n,m)z^n .
\end{equation}
\end{defn}

\begin{defn} (Generating function for unrestricted tilings)
\begin{equation}
T(z,m) = \sum_{n\ge 0} T(n,m)z^n. 
\end{equation}
\end{defn}
Mutual insertion of the previous three equations shows that the
generating function of the row sums is the sum over the generating
functions of the columns,
\begin{equation}
T(z,m) 
= \sum_{t\ge 0} T_t(z,m).
\label{eq.gfsumrule}
\end{equation}
Floor tilings may be equivalent (congruent) in the
sense that reflections of the entire stack of tiles along
the horizontal and/or vertical axis through the middle of the floor
displays another tiling. This group of symmetry operations
of the rectangle will be augmented
by the 90-degree rotations if the floor is a square, i.e., if $m=n$.
If only one representative of the set of 4 or 8 congruent tilings
is counted, an overbar is added to the capital $T$ to 
denote the number of incongruent tilings:
\begin{defn} (Counts by number of 4 tile meets)
$\bar T_t(n,m)$ denotes the number of incongruent tilings which contain $t$
points where 4 tiles meet.
\end{defn}
Row sums and generating functions for counts
of incongruent tilings are defined as for the full counts:
\begin{defn} (Number of incongruent tilings of the $m\times n$ floor)
\begin{equation}
\bar T(n,m) = \sum_{t\ge 0} \bar T_t(n,m).
\end{equation}
\end{defn}

\begin{defn} (Generating function of the incongruent tilings which
contain $t$ points where 4 tiles meet)
\begin{equation}
\bar T_t(z,m) = \sum_{n\ge 0} \bar T_t(n,m)z^n.
\end{equation}
\end{defn}

\begin{defn} (Generating function for the number of incongruent tilings)
\begin{equation}
\bar T(z,m) = \sum_{n\ge 0} \bar T(n,m)z^n .
\end{equation}
\end{defn}

The reader should keep in mind
that these notations do not show the dependence on the tile's dimensions.
We count coverage by $1\times 2$ tiles in Section \ref{sec.12},
by $1\times 3$ tiles in Section \ref{sec.13},
by $1\times 4$ tiles in Section \ref{sec.14},
and by
by $2\times 3$ tiles in Section \ref{sec.23}; each section
defines a different set of $T$, $\bar T$ and associated generating functions.

\subsection{Classification by Slide Line Count}
Tilings of a $m \times n$ floor may be
stacks of tilings of the $m_1 \times n$ floor and of the $m_2 \times n$
floor where $m=m_1+m_2$ (or stacks of three or more such tilings).
The resultant stack then has one or more ``slide lines'' that run
parallel to the long edge of the floor and do not cut through any of the 
tiles.

Similar to the degree of violation
of the Tatami property, the number of slides lines (between 0 and $m-1$, inclusive)
allows a (rough) classification of all tilings of the $m\times n$ floor
with tiles of shape $t_m\times t_n$.

\begin{defn}
$\hat T_s(n,m)$ denotes the number of tilings of the $m\times n$
floor with $s$ slide lines.
\end{defn}

\begin{defn} \label{def.slid}
$\hat {\bar T}_s(n,m)$ denotes the number of incongruent tilings
of the $m\times n$
floor with $s$ slide lines.
\end{defn}

Equivalent to (\ref{eq.rsum}), tables of $\hat T_s(n,m)$
have row sums $T(n,m)$ and generating functions along the columns:
\begin{eqnarray}
T(n,m)&=&\sum_{s=0}^{m-1} \hat T_s(n,m);\\
T(z,m)&=&\sum_{s=0}^{m-1} \hat T_s(z,m).
\label{eq.sumrul2}
\end{eqnarray}
Section \ref{sec.sli} is dedicated to tables of tilings refined according
to their number of side lines.

\section{Domino Tiling} \label{sec.12}
Tilings with $1\times 2$ tiles (also known as dominos) are
represented by Tables \ref{tab.run2}--\ref{tab.run9i}.
\subsection{Results (full count)}

The values for Tatami tilings (domino tilings)---represented by
the columns $T_0(n,m)$ in Tables \ref{tab.run2}--\ref{tab.run9}
and $\bar T_0(n,m)$ in \ref{tab.run3i}--\ref{tab.run9i} are well understood
by a previous enumeration by Ruskey and Woodcock
\cite{RuskeyEJC16}.
The row sums $T(n,m)$ of Tables \ref{tab.run3}--\ref{tab.run9}
have already been reported by Klarner and Pollack
\cite{KlarnerDM32,HockDAM8,StanleyDAM12,StrehlAAM27}.

We turn to a discussion of individual tables.

Row sums of Table \ref{tab.run2} are the Fibonacci numbers,
sequence A000045 in the Encyclopedia of Integer Sequences \cite[A000045]{EIS}.
The generating function for column 0 is \cite[A068921]{EIS}
\begin{thm} (Table \ref{tab.run2})
\begin{equation}
T_0(z,2) = \frac{1+z^2}{1-z-z^3}.
\end{equation}
\end{thm}
\begin{conj} (Table \ref{tab.run2})
\begin{equation}
T_1(z,2) = z^4\frac{1}{(1-z-z^3)^2}.
\end{equation}
\end{conj}
\begin{conj} (Table \ref{tab.run2})
\begin{equation}
T_2(z,2) = z^6\frac{1-z}{(1-z-z^3)^3}.
\end{equation}
\end{conj}

\small
\begin{table}
\caption{Number $T(n,2)$ and $T_t(n,2)$ of domino tilings of $2\times n$ boards.
}
\begin{tabular}{rr|rrrrrrrrrrrrrrrrrrr}
$n$ & & 0 & 1 & 2 & 3 & 4 & 5 & 6 & 7 &8 & 9 & 10\\
\hline
1&1&1&0&0&0&0&0&0&0&0&0&0\\
2&2&2&0&0&0&0&0&0&0&0&0&0\\
3&3&3&0&0&0&0&0&0&0&0&0&0\\
4&5&4&1&0&0&0&0&0&0&0&0&0\\
5&8&6&2&0&0&0&0&0&0&0&0&0\\
6&13&9&3&1&0&0&0&0&0&0&0&0\\
7&21&13&6&2&0&0&0&0&0&0&0&0\\
8&34&19&11&3&1&0&0&0&0&0&0&0\\
9&55&28&18&7&2&0&0&0&0&0&0&0\\
10&89&41&30&14&3&1&0&0&0&0&0&0\\
11&144&60&50&24&8&2&0&0&0&0&0&0\\
12&233&88&81&43&17&3&1&0&0&0&0&0\\
13&377&129&130&77&30&9&2&0&0&0&0&0\\
14&610&189&208&132&57&20&3&1&0&0&0&0\\
15&987&277&330&224&108&36&10&2&0&0&0&0\\
16&1597&406&520&379&193&72&23&3&1&0&0&0\\
17&2584&595&816&633&342&143&42&11&2&0&0&0\\
18&4181&872&1275&1047&605&264&88&26&3&1&0&0\\
19&6765&1278&1984&1722&1052&485&182&48&12&2&0&0\\
20&10946&1873&3077&2814&1808&891&345&105&29&3&1&0\\
21&17711&2745&4758&4570&3088&1602&654&225&54&13&2&0\\
22&28657&4023&7337&7385&5232&2843&1242&436&123&32&3&1\\
23&46368&5896&11286&11880&8796&5014&2298&850&272&60&14&2\\
24&75025&8641&17322&19029&14699&8760&4193&1663&537&142&35&3\\
25&121393&12664&26532&30363&24426&15167&7606&3155&1074&323&66&15\\
26&196418&18560&40563&48279&40371&26084&13650&5900&2159&648&162&38\\
27&317811&27201&61908&76518&66404&44571&24250&10976&4188&1327&378&72\\
\end{tabular}
\label{tab.run2}
\end{table}
\normalsize

Column $T(n,m)$ of Table \ref{tab.run3} is
\cite[A001835]{EIS}:
\begin{thm} (Table \ref{tab.run3})
\begin{equation}
T(z,3)= \frac{1-z^2}{1-4z^2+z^4}.
\label{eq.run3}
\end{equation}
\end{thm}
Column $T_0(n,m)$ is \cite[A068922]{EIS}:
\begin{thm} (Table \ref{tab.run3})
\begin{equation}
T_0(z,3)=z^2-1+\frac{2}{1-z^2-z^4}.
\end{equation}
\end{thm}
For the next two columns we have:
\begin{conj} (Table \ref{tab.run3})
\begin{equation}
T_1(z,3) = 2z^6\frac{3+3z^2-3z^4-2z^6}{(1-z^2-z^4)^2}.
\end{equation}
\end{conj}
\begin{conj} (Table \ref{tab.run3})
\begin{equation}
T_2(z,3) =
z^4+2z^6\frac{6+11z^2-9z^4-12z^6-z^8+z^{10}}{ (1-z^2-z^4)^3 }.
\end{equation}
\end{conj}

\small
\begin{table}
\caption{Number $T(n,3)$ and $T_t(n,3)$ of domino tilings of $3\times n$ boards.
}
\begin{tabular}{rr|rrrrrrrrrrrrrrrrrrr}
$n$ & & 0 & 1 & 2 & 3 & 4 & 5 & 6 & 7 &8 \\
\hline
2&3&3&0&0&0&0&0&0&0&0\\
4&11&4&6&1&0&0&0&0&0&0\\
6&41&6&18&12&4&1&0&0&0&0\\
8&153&10&36&58&32&12&4&1&0&0\\
10&571&16&74&156&174&92&40&14&4&1\\
12&2131&26&142&384&578&518&284&128&50&16\\
14&7953&42&268&860&1646&1998&1578&886&422&170\\
16&29681&68&494&1838&4202&6408&6672&4912&2816&1374\\
18&110771&110&898&3780&10024&18238&23500&22004&15564&9010\\
20&413403&178&1612&7566&22732&47852&73190&83316&72300&49900\\
22&1542841&288&2866&14816&49638&118242&208586&279346&289268&237802\\
24&5757961&466&5054&28512&105190&279056&556128&854582&1030190&992446\\
26&21489003&754&8852&54080&217586&634978&1407596&2435918&3348484&3712826\\
28&80198051&1220&15414&101338&441146&1402550&3417114&6564072&10121734&12706318\\
30&299303201&1974&26706&187932&879436&3022324&8016016&16898784&28848726&40432874\\
32&1117014753&3194&46068&345410&1728056&6377980&18272816&41888806&78334170&121156904\\
\label{tab.run3}
\end{tabular}
\end{table}
\normalsize

Column $T(n,m)$ of Table \ref{tab.run4} is \cite[A005178]{EIS}:
\begin{thm} (Table \ref{tab.run4})
\begin{equation}
T(z,4) =
\frac{1-z^2}{ 1-z-5z^2-z^3+z^4 }.
\end{equation}
\end{thm}
Column $T_0(n,m)$ is \cite[A068923]{EIS}:
\begin{thm} (Table \ref{tab.run4})
\begin{equation}
T_0(z,4) =
-1+3z^2+2z^3+\frac{2+z+z^2+z^4}{1-z^3-z^5}.
\end{equation}
\end{thm}

\small
\begin{table}
\caption{Number $T(n,4)$ and $T_t(n,4)$ of domino tilings of $4\times n$ boards.
}
\begin{tabular}{rr|rrrrrrrrrrrrrrr}
$n$ & & 0 & 1 & 2 & 3 & 4 & 5 & 6 & 7 &8 \\
\hline
1&1&1&0&0&0&0&0&0&0&0\\
2&5&4&1&0&0&0&0&0&0&0\\
3&11&4&6&1&0&0&0&0&0&0\\
4&36&2&14&18&2&0&0&0&0&0\\
5&95&3&20&43&26&3&0&0&0&0\\
6&281&3&24&74&113&56&10&1&0&0\\
7&781&3&30&123&248&272&88&15&2&0\\
8&2245&5&34&173&477&717&596&192&41&9\\
9&6336&5&38&252&792&1581&1962&1260&358&72\\
10&18061&6&53&318&1254&3036&4848&4906&2614&780\\
11&51205&8&64&404&1864&5361&10426&13826&11720&5479\\
12&145601&8&79&553&2600&8817&20258&32969&37265&26921\\
13&413351&11&98&717&3600&13656&36536&70299&97856&95707\\
14&1174500&13&116&937&4960&20364&61894&137367&227287&275616\\
15&3335651&14&150&1207&6760&29639&99486&251723&479984&694119\\
16&9475901&19&181&1532&9132&42467&154532&435514&944178&1576953\\
17&26915305&21&220&1989&12146&60088&233572&721548&1747606&3315355\\
18&76455961&25&277&2525&16126&83677&346240&1155602&3079704&6537383\\
19&217172736&32&330&3203&21298&115525&504476&1800597&5218642&12216799\\
20&616891945&35&413&4075&27895&158158&724090&2746127&8557790&21856138\\
21&1752296281&44&506&5126&36468&214484&1028182&4110548&13666350&37689958\\
22&4977472781&53&608&6491&47369&289373&1443894&6059733&21341722&63027592\\
23&14138673395&60&762&8157&61352&387604&2009478&8818495&32701156&102707437\\
\end{tabular}
\label{tab.run4}
\end{table}
\normalsize

Column $T(n,m)$ of Table \ref{tab.run5} is
\cite[A003775]{EIS}; column $T_0(n,m)$ is \cite[A068924]{EIS}.

\small
\begin{table}
\caption{Number $T(n,5)$ and $T_t(n,5)$ of domino tilings of $5\times n$ boards.
}
\begin{tabular}{rr|rrrrrrrrrrrrrrrrrrr}
$n$ & & 0 & 1 & 2 & 3 & 4 & 5 & 6 & 7 &8 \\
\hline
2&8&6&2&0&0&0&0&0&0&0\\
4&95&3&20&43&26&3&0&0&0&0\\
6&1183&2&32&147&332&343&220&92&14&1\\
8&14824&2&38&271&1046&2695&3730&3370&2206&1061\\
10&185921&4&38&422&2302&8144&21058&35753&41254&35275\\
12&2332097&4&68&532&4074&19405&65662&169536&319122&439701\\
14&29253160&6&82&864&6206&37520&163410&538759&1390882&2792397\\
16&366944287&8&114&1224&10120&64464&339114&1382583&4458614&11576972\\
\end{tabular}
\label{tab.run5}
\end{table}
\normalsize

Column $T(n,m)$ of Table \ref{tab.run6} is
\cite[A028468]{EIS} and column $T_0(n,m)$ is \cite[A068925]{EIS}, so
the two generating functions are known:
\begin{thm} (Table \ref{tab.run6})
\begin{equation}
T(z,6) =
\frac{1}{13}
[\frac{5+12z+3z^2}{1+5z+6z^2+z^3}
+\frac{6-10z+3z^2}{1-6z+5z^2-z^3}
+\frac{1}{1+z}
+\frac{1}{1-z}
].
\end{equation}
\end{thm}
\begin{thm} (Table \ref{tab.run6})
\begin{equation}
T_0(z,6) =
-1+9z^2+6z^3+2z^4
+\frac{1}{13}[
\frac{12-9z}{1-z+z^2}
+
\frac{14+24z+19z^2+7z^3+4z^4}{1+z-z^3-z^4-z^5}
].
\end{equation}
\end{thm}

\small
\begin{table}
\caption{Number $T(n,6)$ and $T_t(n,6)$ of domino tilings of $6\times n$ boards.
}
\begin{tabular}{rr|rrrrrrrrrrrrrrr}
$n$ & & 0 & 1 & 2 & 3 & 4 & 5 & 6 & 7 &8 &9 & 10 \\
\hline
1&1&1&0&0&0&0&0&0&0&0&0&0\\
2&13&9&3&1&0&0&0&0&0&0&0&0\\
3&41&6&18&12&4&1&0&0&0&0&0&0\\
4&281&3&24&74&113&56&10&1&0&0&0&0\\
5&1183&2&32&147&332&343&220&92&14&1&0&0\\
6&6728&2&24&202&714&1326&1842&1508&782&288&38&2\\
7&31529&2&38&255&1214&3242&6126&7909&6784&3888&1576&438\\
8&167089&1&37&290&1665&6315&15372&27943&36693&35232&24777&12630\\
9&817991&1&42&350&2288&10342&32878&73852&126190&167827&167658&125352\\
10&4213133&2&39&429&2990&15183&57599&162389&350681&588547&783469&820190\\
11&21001799&3&24&470&3756&21304&93950&313656&817194&1673001&2739234&3644387\\
12&106912793&4&39&567&4624&28961&141741&543976&1654825&4022461&7941693&12855987\\
13&536948224&3&46&525&5408&37887&203462&880052&3045892&8594250&19871300&37953976\\
\end{tabular}
\label{tab.run6}
\end{table}
\normalsize

Column $T(n,m)$ of Table \ref{tab.run7} is
\cite[A028469]{EIS}.

\small
\begin{table}
\caption{Number $T(n,7)$ and $T_t(n,7)$ of domino tilings of $7\times n$ boards.
}
\begin{tabular}{rr|rrrrrrrrrrrrrrr}
$n$ & & 0 & 1 & 2 & 3 & 4 & 5 & 6 & 7 &8 &9 & 10 \\
\hline
2&21&13&6&2&0&0&0&0&0&0&0&0\\
4&781&3&30&123&248&272&88&15&2&0&0&0\\
6&31529&2&38&255&1214&3242&6126&7909&6784&3888&1576&438\\
8&1292697&2&36&372&2606&11163&37162&90677&170694&244705&265598&220083\\
10&53175517&0&38&470&3802&24137&113316&404254&1148114&2630509&4931018&7521981\\
12&2188978117&2&28&528&5470&40512&239344&1130007&4329026&13626282&35621942&78462193\\
\end{tabular}
\label{tab.run7}
\end{table}
\normalsize

Column $T(n,m)$ of Table \ref{tab.run8} is
\cite[A028470]{EIS}.

\small
\begin{table}
\caption{Number $T(n,8)$ and $T_t(n,8)$ of domino tilings of $8\times n$ boards.
}
\begin{tabular}{rr|rrrrrrrrrrrrrrr}
$n$ & & 0 & 1 & 2 & 3 & 4 & 5 & 6 & 7 &8 &9 & 10 \\
\hline
1&1&1&0&0&0&0&0&0&0&0&0&0\\
2&34&19&11&3&1&0&0&0&0&0&0&0\\
3&153&10&36&58&32&12&4&1&0&0&0&0\\
4&2245&5&34&173&477&717&596&192&41&9&1&0\\
5&14824&2&38&271&1046&2695&3730&3370&2206&1061&324&68\\
6&167089&1&37&290&1665&6315&15372&27943&36693&35232&24777&12630\\
7&1292697&2&36&372&2606&11163&37162&90677&170694&244705&265598&220083\\
8&12988816&2&34&310&3110&16712&68504&224484&569884&1127972&1798390&2307764\\
9&108435745&2&42&445&4058&24514&119062&453678&1403192&3503161&7121656&11916877\\
10&1031151241&1&29&460&4313&30449&176077&785735&2864328&8638123&21830440&46070190\\
\end{tabular}
\label{tab.run8}
\end{table}
\normalsize

Column $T(n,m)$ of Table \ref{tab.run9} is
\cite[A028471]{EIS}.

\small
\begin{table}
\caption{Number $T(n,9)$ and $T_t(n,9)$ of domino tilings of $9\times n$ boards.
}
\begin{tabular}{rr|rrrrrrrrrrrrrrr}
$n$ & & 0 & 1 & 2 & 3 & 4 & 5 & 6 & 7 &8 &9 & 10 \\
\hline
2&55&28&18&7&2&0&0&0&0&0&0&0\\
4&6336&5&38&252&792&1581&1962&1260&358&72&14&2\\
6&817991&1&42&350&2288&10342&32878&73852&126190&167827&167658&125352\\
8&108435745&2&42&445&4058&24514&119062&453678&1403192&3503161&7121656&11916877\\
\end{tabular}
\label{tab.run9}
\end{table}
\normalsize

\clearpage
\subsection{Results (incongruent)}
Counts of domino tilings where only one representative of the
roto-reflected copies of each tiling is counted are shown
in Tables \ref{tab.run2i}--\ref{tab.run9i}.

Table \ref{tab.run2i} is characterized by:
\begin{thm} (Table \ref{tab.run2i} \cite[A060312]{EIS})
\begin{equation}
\bar T(z,2)
= -z^2+\frac{1}{2}[\frac{1}{1-z-z^2}
+\frac{1+z+z^2}{1-z^2-z^4}].
\end{equation}
\end{thm}
\begin{thm} (Table \ref{tab.run2i} \cite[A068927]{EIS})
\begin{equation}
\bar T_0(z,2) 
= -z^2+\frac{1}{2}[
\frac{1+z+z^2+z^5}{1-z^2-z^6}
+\frac{1+z^2}{1-z-z^3}
].
\end{equation}
\end{thm}
\begin{conj} (Table \ref{tab.run2i})
\begin{equation}
\bar T_1(z,2) = z^4\frac{(1+z+z^2)(1-z)^2 }{ (1-z^2-z^6)(1-z-z^3)^2 }
\end{equation}
\end{conj}

\small
\begin{table}
\caption{Number $\bar T(n,2)$ and $\bar T_t(n,2)$ of incongruent domino tilings of $2\times n$ boards.
}
\begin{tabular}{rr|rrrrrrrrrrrrrrr}

$n$ & & 0 & 1 & 2 & 3 & 4 & 5 & 6 & 7 &8 &9 &10 &11\\

\hline
1&1&1&0&0&0&0&0&0&0&0&0&0&0\\
2&1&1&0&0&0&0&0&0&0&0&0&0&0\\
3&2&2&0&0&0&0&0&0&0&0&0&0&0\\
4&4&3&1&0&0&0&0&0&0&0&0&0&0\\
5&5&4&1&0&0&0&0&0&0&0&0&0&0\\
6&9&6&2&1&0&0&0&0&0&0&0&0&0\\
7&12&8&3&1&0&0&0&0&0&0&0&0&0\\
8&21&12&6&2&1&0&0&0&0&0&0&0&0\\
9&30&16&9&4&1&0&0&0&0&0&0&0&0\\
10&51&24&16&8&2&1&0&0&0&0&0&0&0\\
11&76&33&25&13&4&1&0&0&0&0&0&0&0\\
12&127&49&42&24&9&2&1&0&0&0&0&0&0\\
13&195&69&65&40&15&5&1&0&0&0&0&0&0\\
14&322&102&106&70&30&11&2&1&0&0&0&0&0\\
15&504&145&165&115&54&19&5&1&0&0&0&0&0\\
16&826&214&263&196&99&39&12&2&1&0&0&0&0\\
17&1309&307&408&322&171&73&21&6&1&0&0&0&0\\
18&2135&452&642&535&306&137&46&14&2&1&0&0&0\\
19&3410&653&992&870&526&246&91&25&6&1&0&0&0\\
20&5545&960&1545&1426&910&454&176&56&15&2&1&0&0\\
21&8900&1393&2379&2300&1544&808&327&114&27&7&1&0&0\\
22&14445&2046&3678&3723&2626&1438&626&224&64&17&2&1&0\\
23&23256&2978&5643&5965&4398&2519&1149&429&136&31&7&1&0\\
24&37701&4371&8675&9564&7365&4409&2106&842&273&75&18&2&1\\
25&60813&6376&13266&15222&12213&7605&3803&1586&537&163&33&8&1\\
26&98514&9354&20302&24219&20210&13091&6842&2972&1086&331&84&20&2\\
27&159094&13665&30954&38324&33202&22324&12125&5503&2094&668&189&37&8\\
28&257608&20041&47198&60582&54412&37926&21414&10128&4017&1380&390&96&21\\
29&416325&29307&71770&95372&88678&63966&37455&18366&7621&2716&805&220&39\\
\end{tabular}
\label{tab.run2i}
\end{table}
\normalsize

Column $\bar T_0(n,3)$ of Table \ref{tab.run3i} is
found in \cite[A068928]{EIS}, which shows
\begin{thm} (Table \ref{tab.run3i})
\begin{equation}
\bar T_0(z,3) =
z^2+\frac{1}{2}
[\frac{1+z^2+z^4}{1-z^4-z^8}
+\frac{1}{1-z^2-z^4}
].
\end{equation}
\end{thm}
For the row sums we formulate
\begin{conj} (Table \ref{tab.run3i})
\begin{equation}
\bar T(z,3) =
\frac{1}{8}[\frac{1}{1+z}+\frac{1}{1-z}]
+\frac{1+2z^2-z^6}{2(1-4z^4+z^8)}
+\frac{1-z^2}{4(1-4z^2+z^4)}.
\label{eq.run3i}
\end{equation}
\end{conj}

\small
\begin{table}
\caption{Number $\bar T(n,3)$ and $\bar T_t(n,3)$ of incongruent domino tilings of $3\times n$ boards.
}
\begin{tabular}{rr|rrrrrrrrrrrrrrr}
$n$ & & 0 & 1 & 2 & 3 & 4 & 5 & 6 & 7 &8 &9 \\
\hline
2&2&2&0&0&0&0&0&0&0&0&0\\
4&5&2&2&1&0&0&0&0&0&0&0\\
6&14&2&5&5&1&1&0&0&0&0&0\\
8&46&4&10&17&9&4&1&1&0&0&0\\
10&156&5&19&43&45&27&10&5&1&1&0\\
12&561&9&37&102&148&137&73&35&13&5&1\\
14&2037&12&68&223&414&514&398&232&106&47&15\\
16&7525&21&126&473&1058&1626&1679&1248&709&353&148\\
18&27874&30&226&960&2512&4596&5885&5549&3899&2282&1131\\
20&103741&51&407&1919&5699&12027&18329&20908&18106&12535&7287\\
22&386386&76&719&3732&12421&29645&52174&69982&72350&59598&40361\\
24&1440946&127&1270&7182&26329&69916&139112&213893&257658&248374&196028\\
26&5374772&195&2217&13571&54419&158927&351963&609370&837226&928729&845193\\
28&20054945&322&3864&25437&110347&350976&854462&1641695&2530761&3177495&3292433\\
30&74835209&504&6683&47075&219901&755959&2004146&4225650&7212468&10109805&11772446\\
\end{tabular}
\label{tab.run3i}
\end{table}
\normalsize

Column $\bar T_0(n,4)$ of Table \ref{tab.run4i} is
found in \cite[A068929]{EIS}, which provides
\begin{thm} (Table \ref{tab.run4i})
\begin{equation}
\bar T_0(z,4) =
2z^2+z^3-z^4 +
\frac{1}{2}[
z\frac{1+z+2z^3+z^4+z^5+z^8}{1-z^6-z^{10}}
+\frac{2+z+z^2+z^4}{1-z^3-z^5}
]
.
\end{equation}
\end{thm}

\small
\begin{table}
\caption{Number $\bar T(n,4)$ and $\bar T_t(n,4)$ of incongruent domino tilings of $4\times n$ boards.
}
\begin{tabular}{rr|rrrrrrrrrrrrrrr}
$n$ & & 0 & 1 & 2 & 3 & 4 & 5 & 6 & 7 &8 &9 \\
\hline
1&1&1&0&0&0&0&0&0&0&0&0\\
2&4&3&1&0&0&0&0&0&0&0&0\\
3&5&2&2&1&0&0&0&0&0&0&0\\
4&9&1&3&4&1&0&0&0&0&0&0\\
5&33&2&6&15&8&2&0&0&0&0&0\\
6&98&2&9&25&34&22&5&1&0&0&0\\
7&230&2&10&38&67&80&26&6&1&0&0\\
8&658&3&12&57&134&203&166&62&16&4&1\\
9&1725&3&13&77&214&428&512&346&102&24&5\\
10&4876&4&19&100&347&819&1272&1300&702&230&61\\
11&13378&5&22&124&505&1421&2690&3577&3009&1451&428\\
12&37794&5&30&171&715&2347&5236&8493&9535&6947&3007\\
13&105761&6&35&221&983&3597&9383&17960&24818&24353&15439\\
14&299221&8&43&292&1369&5386&15924&35057&57610&69816&60038\\
15&844219&8&56&372&1857&7809&25514&64024&121268&175164&187455\\
16&2392040&11&68&480&2530&11214&39713&110741&238496&397460&508485\\
17&6773154&12&83&617&3359&15854&59923&183201&440817&834301&1236712\\
18&19211023&14&108&792&4489&22117&88979&293450&776790&1644337&2768166\\
19&54485124&17&125&1006&5930&30525&129574&456942&1315446&3070613&5787975\\
20&154636939&20&162&1284&7813&41885&186189&697095&2157279&5492248&11435056\\
21&438909205&24&196&1617&10217&56806&264434&1043191&3444055&9467905&21540497\\
\end{tabular}
\label{tab.run4i}
\end{table}
\normalsize

Column $\bar T_0(n,5)$ of Table \ref{tab.run5i} is
\cite[A068930]{EIS}.

\small
\begin{table}
\caption{Number $\bar T(n,5)$ and $\bar T_t(n,5)$ of incongruent domino tilings of $5\times n$ boards.
}
\begin{tabular}{rr|rrrrrrrrrrrrrrr}
$n$ & & 0 & 1 & 2 & 3 & 4 & 5 & 6 & 7 &8 &9 & 10 \\
\hline
2&5&4&1&0&0&0&0&0&0&0&0&0\\
4&33&2&6&15&8&2&0&0&0&0&0&0\\
6&329&1&8&43&88&96&58&30&4&1&0&0\\
8&3818&1&10&75&266&704&945&869&561&277&85&21\\
10&46878&1&10&113&579&2070&5288&9024&10354&8904&5815&3054\\
12&584386&2&17&141&1026&4903&16447&42564&79871&110221&114726&94591\\
14&7318152&2&21&221&1557&9441&40889&134928&347878&698731&1081974&1313568\\
16&91752831&3&29&318&2535&16202&84840&346053&1114901&2895523&6060041&10207144\\
\end{tabular}
\label{tab.run5i}
\end{table}
\normalsize

Column $\bar T_0(n,6)$ of Table \ref{tab.run6i} is
\cite[A068931]{EIS}.

\small
\begin{table}
\caption{Number $\bar T(n,6)$ and $\bar T_t(n,6)$ of incongruent domino tilings of $6\times n$ boards.
}
\begin{tabular}{rr|rrrrrrrrrrrrrrr}
$n$ & & 0 & 1 & 2 & 3 & 4 & 5 & 6 & 7 &8 &9 & 10 \\
\hline
1&1&1&0&0&0&0&0&0&0&0&0&0\\
2&9&6&2&1&0&0&0&0&0&0&0&0\\
3&14&2&5&5&1&1&0&0&0&0&0&0\\
4&98&2&9&25&34&22&5&1&0&0&0&0\\
5&329&1&8&43&88&96&58&30&4&1&0&0\\
6&930&1&4&29&96&181&247&211&105&48&7&1\\
7&8121&1&11&69&310&834&1558&2032&1733&1024&405&126\\
8&42873&1&13&82&438&1642&3926&7143&9339&9018&6335&3296\\
9&206420&1&13&97&587&2635&8279&18624&31728&42279&42204&31700\\
10&1060866&1&12&127&781&3903&14570&40942&88162&147973&196828&206307\\
11&5265647&2&7&130&970&5423&23645&78788&204821&419280&686092&913146\\
12&26782279&2&12&164&1207&7403&35771&136737&414961&1007932&1988590&3219095\\
\end{tabular}
\label{tab.run6i}
\end{table}
\normalsize

\small
\begin{table}
\caption{Number $\bar T(n,7)$ and $\bar T_t(n,7)$ of incongruent domino tilings of $7\times n$ boards.
}
\begin{tabular}{rr|rrrrrrrrrrrrrrr}
$n$ & & 0 & 1 & 2 & 3 & 4 & 5 & 6 & 7 &8 &9 & 10 \\
\hline
2&12&8&3&1&0&0&0&0&0&0&0&0\\
4&230&2&10&38&67&80&26&6&1&0&0&0\\
6&8121&1&11&69&310&834&1558&2032&1733&1024&405&126\\
8&324617&1&9&98&661&2841&9347&22828&42804&61449&66562&55266\\
10&13303375&0&10&123&958&6083&28402&101293&287287&658240&1233331&1881672\\
\end{tabular}
\end{table}
\normalsize

\small
\begin{table}
\caption{Number $\bar T(n,8)$ and $\bar T_t(n,8)$ of incongruent domino tilings of $8\times n$ boards.
}
\begin{tabular}{rr|rrrrrrrrrrrrrrr}
$n$ & & 0 & 1 & 2 & 3 & 4 & 5 & 6 & 7 &8 &9 & 10 \\
\hline
1&1&1&0&0&0&0&0&0&0&0&0&0\\
2&21&12&6&2&1&0&0&0&0&0&0&0\\
3&46&4&10&17&9&4&1&1&0&0&0&0\\
4&658&3&12&57&134&203&166&62&16&4&1&0\\
5&3818&1&10&75&266&704&945&869&561&277&85&21\\
6&42873&1&13&82&438&1642&3926&7143&9339&9018&6335&3296\\
7&324617&1&9&98&661&2841&9347&22828&42804&61449&66562&55266\\
8&1629189&1&6&43&400&2133&8647&28309&71540&141600&225402&289306\\
9&27129182&1&12&119&1027&6191&29862&113815&351290&877086&1781753&2981718\\
\end{tabular}
\end{table}
\normalsize

\small
\begin{table}
\caption{Number $\bar T(n,9)$ and $\bar T_t(n,9)$ of incongruent domino tilings of $9\times n$ boards.
}
\begin{tabular}{rr|rrrrrrrrrrrrrrr}
$n$ & & 0 & 1 & 2 & 3 & 4 & 5 & 6 & 7 &8 &9 & 10 \\
\hline
2&30&16&9&4&1&0&0&0&0&0&0&0\\
4&1725&3&13&77&214&428&512&346&102&24&5&1\\
6&206420&1&13&97&587&2635&8279&18624&31728&42279&42204&31700\\
8&27129182&1&12&119&1027&6191&29862&113815&351290&877086&1781753&2981718\\
\end{tabular}
\label{tab.run9i}
\end{table}
\normalsize

\clearpage
\section{Tiling with $1\times 3$ tiles}\label{sec.13}
Tilings with tiles of shape $1\times 3$ are
represented by Tables \ref{tab.run3_3_1}--\ref{tab.run9i_3_1}.
\subsection{Results (full count)}

The row sums $T(n,3)$ in Table \ref{tab.run3_3_1} 
are found in \cite[A000930]{EIS}.
\begin{thm} (Table \ref{tab.run3_3_1})
\begin{equation}
T(z,3)
= \frac{1}{1-z-z^3}.
\label{eq.run3_3_1}
\end{equation}
\end{thm}

Column $T_0(n,3)$ appears to
be \cite[A003269]{EIS}:
\begin{thm} (Table \ref{tab.run3_3_1})
\begin{equation}
T_0(z,3)=
\frac{1+z^3}{1-z-z^4}.
\end{equation}
\end{thm}
\begin{conj} (Table \ref{tab.run3_3_1})
\begin{equation}
T_t(z,3)=
\left\{
\begin{array}{ll}
z^{3(t/2+1)}\frac{(1-z)^{t/2-1}}{(1-z-z^4)^{t/2+1}},& \mathrm{even }\, t>0\\
0, & \mathrm{odd}\, t.
\end{array}
\right.
\end{equation}
\end{conj}
The previous three equations comply with the sum rule (\ref{eq.gfsumrule}).

\small
\begin{table}
\caption{Number $T(n,3)$ and $T_t(n,3)$ of tilings of $3\times n$ boards with 
$1\times 3$ tiles.
}
\begin{tabular}{rr|rrrrrrrrrrrrrrr}
$n$ & & 0 & 1 & 2 & 3 & 4 & 5 & 6 & 7 &8 &9 & 10 & 11 & 12 & 13 & 14 \\
\hline
1&1&1&0&0&0&0&0&0&0&0&0&0&0&0&0&0 \\
2&1&1&0&0&0&0&0&0&0&0&0&0&0&0&0&0 \\
3&2&2&0&0&0&0&0&0&0&0&0&0&0&0&0&0 \\
4&3&3&0&0&0&0&0&0&0&0&0&0&0&0&0&0 \\
5&4&4&0&0&0&0&0&0&0&0&0&0&0&0&0&0 \\
6&6&5&0&1&0&0&0&0&0&0&0&0&0&0&0&0 \\
7&9&7&0&2&0&0&0&0&0&0&0&0&0&0&0&0 \\
8&13&10&0&3&0&0&0&0&0&0&0&0&0&0&0&0 \\
9&19&14&0&4&0&1&0&0&0&0&0&0&0&0&0&0 \\
10&28&19&0&7&0&2&0&0&0&0&0&0&0&0&0&0 \\
11&41&26&0&12&0&3&0&0&0&0&0&0&0&0&0&0 \\
12&60&36&0&19&0&4&0&1&0&0&0&0&0&0&0&0 \\
13&88&50&0&28&0&8&0&2&0&0&0&0&0&0&0&0 \\
14&129&69&0&42&0&15&0&3&0&0&0&0&0&0&0&0 \\
15&189&95&0&64&0&25&0&4&0&1&0&0&0&0&0&0 \\
16&277&131&0&97&0&38&0&9&0&2&0&0&0&0&0&0 \\
17&406&181&0&144&0&60&0&18&0&3&0&0&0&0&0&0 \\
18&595&250&0&212&0&97&0&31&0&4&0&1&0&0&0&0 \\
19&872&345&0&312&0&155&0&48&0&10&0&2&0&0&0&0 \\
20&1278&476&0&459&0&240&0&79&0&21&0&3&0&0&0&0 \\
21&1873&657&0&672&0&368&0&134&0&37&0&4&0&1&0&0 \\
22&2745&907&0&979&0&565&0&223&0&58&0&11&0&2&0&0 \\
23&4023&1252&0&1422&0&867&0&356&0&99&0&24&0&3&0&0 \\
24&5896&1728&0&2062&0&1320&0&563&0&175&0&43&0&4&0&1 \\
25&8641&2385&0&2984&0&1995&0&894&0&301&0&68&0&12&0&2 \\
26&12664&3292&0&4308&0&3003&0&1419&0&492&0&120&0&27&0&3 \\
27&18560&4544&0&6206&0&4510&0&2228&0&798&0&220&0&49&0&4 \\
28&27201&6272&0&8925&0&6752&0&3466&0&1304&0&389&0&78&0&13 \\
29&39865&8657&0&12816&0&10071&0&5368&0&2130&0&648&0&142&0&30 \\
30&58425&11949&0&18376&0&14972&0&8294&0&3431&0&1074&0&269&0&55 \\
31&85626&16493&0&26310&0&22201&0&12764&0&5467&0&1800&0&487&0&88 \\
32&125491&22765&0&37620&0&32844&0&19549&0&8673&0&3015&0&824&0&165 \\
33&183916&31422&0&53728&0&48475&0&29818&0&13729&0&4964&0&1392&0&322 \\
34&269542&43371&0&76648&0&71381&0&45341&0&21630&0&8074&0&2387&0&595 \\
35&395033&59864&0&109230&0&104892&0&68748&0&33882&0&13080&0&4089&0&1020 \\
36&578949&82629&0&155507&0&153844&0&103928&0&52824&0&21151&0&6862&0&1753 \\
37&848491&114051&0&221184&0&225239&0&156652&0&82076&0&34016&0&11364&0&3070 \\
38&1243524&157422&0&314325&0&329202&0&235504&0&127113&0&54342&0&18757&0&5367 \\
39&1822473&217286&0&446320&0&480371&0&353204&0&196175&0&86364&0&30917&0&9160
\end{tabular}
\label{tab.run3_3_1}
\end{table}
\normalsize

The row sums $T(n,4)$ of Table \ref{tab.run4_3_1} are \cite[A049086]{EIS}:
\begin{thm} (Table \ref{tab.run4_3_1})
\begin{equation}
T(z,4)=
\frac{1-2z^3+z^6}{1-5z^3+3z^6-z^9}.
\label{eq.run4_3_1}
\end{equation}
\end{thm}
For its column $T_0$ we suspect
\begin{conj} (Table \ref{tab.run4_3_1})
\begin{equation}
T_0(z,4)=
-1+z^3+\frac{2}{1-z^3-2z^6-z^9},
\end{equation}
which means these are essentially the sequence \cite[A002478]{EIS}
multiplied by 2.
\end{conj}
The denominator of the generating function
of this conjecture claims the recurrence
$T_0(n,4)=T_0(n-3,4)+2T_0(n-6,4)+T_0(n-9,4)$. This is interpreted with
the aid of Table \ref{tab.run431} as follows: To tile the $4\times n$
floor, take
\begin{enumerate}
\item
a tiling of the $4\times (n-3)$ floor and attach a $4\times 3$
super-tile with the one (out of 2) orientations that avoid a 4-crossing
at the interface, or
\item
take a tiling of the $4\times (n-6)$ floor
and attach a $4\times 6$ super-tile which come in two shapes, one with
a block of 2 vertical tiles at the left and one with a block of 2 vertical
tiles at the right, or
\item
take a tiling of the $4\times (n-9)$ floor
and attach the $4\times 9$ super-tile with a slide line
with the one (out of 2) orientations---second line in Table \ref{tab.run431}---
that avoid a 4-crossing at the interface.
\end{enumerate}
The claim is essentially that the final $4\times 9$ super-tile of the $4\times n$
Tatami tiling with $1\times 3$ tiles is always one of the shapes in Table \ref{tab.run431}.

\begin{table}
\caption{The $T_0(9,4)=12$ Tatami tilings of row $n=9$ in Table \ref{tab.run4_3_1} with
$\bar T_0(9,4)=4$ incongruent forms of row $n=9$ in Table \ref{tab.run4i_3_1}.}
\begin{tabular}{llll}
\includegraphics[width=0.25\textwidth]{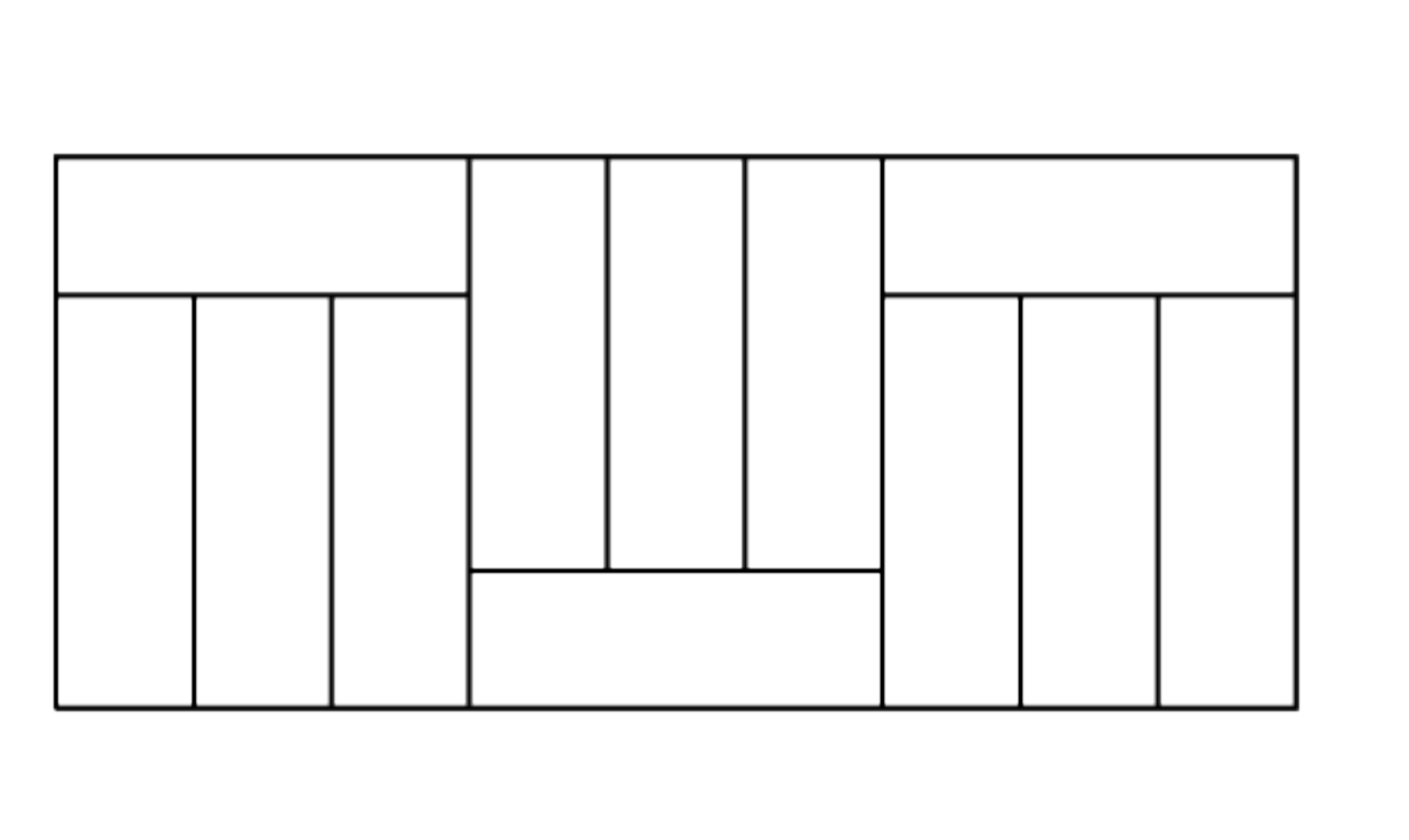}
\includegraphics[width=0.25\textwidth]{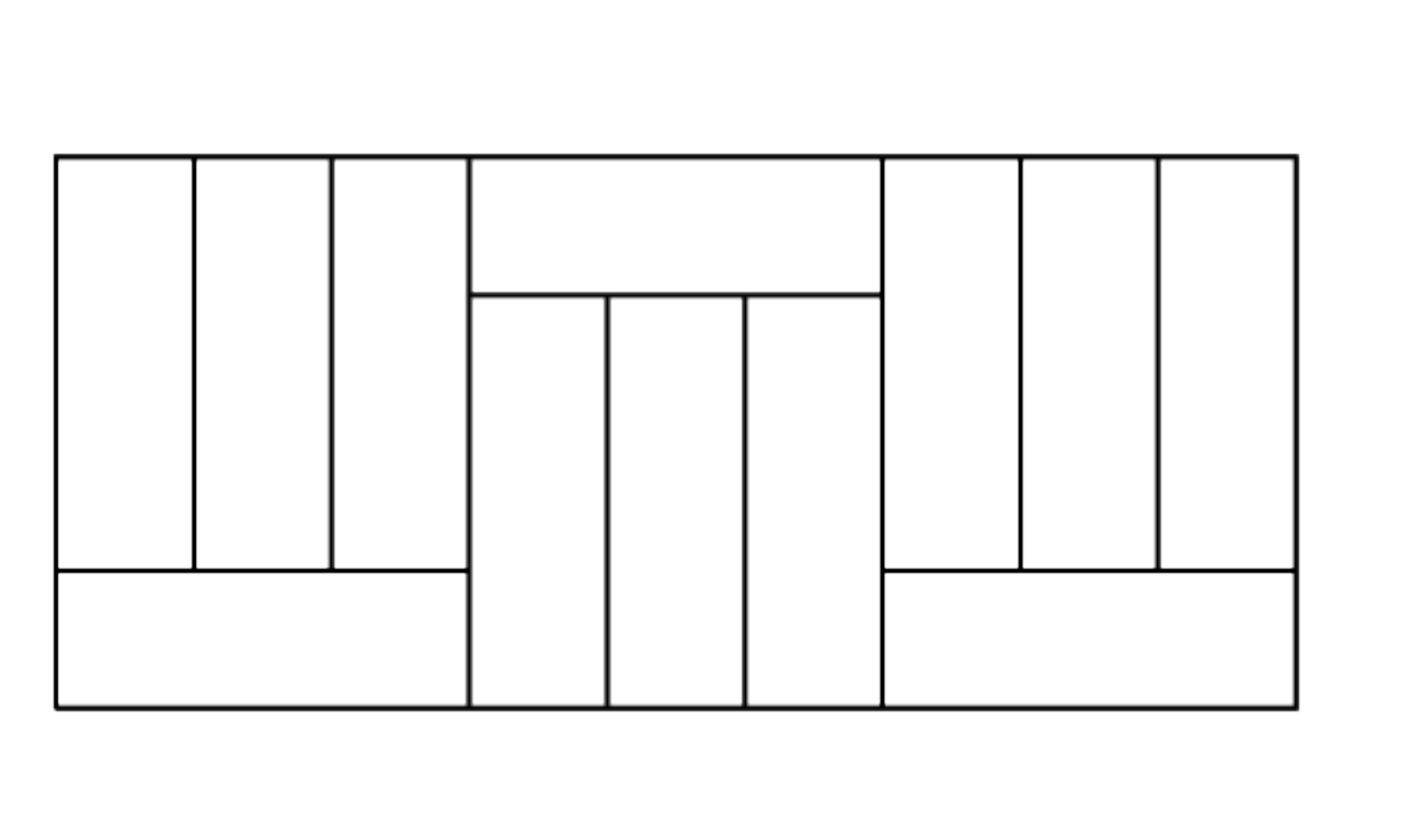}
\\
\hline
\includegraphics[width=0.25\textwidth]{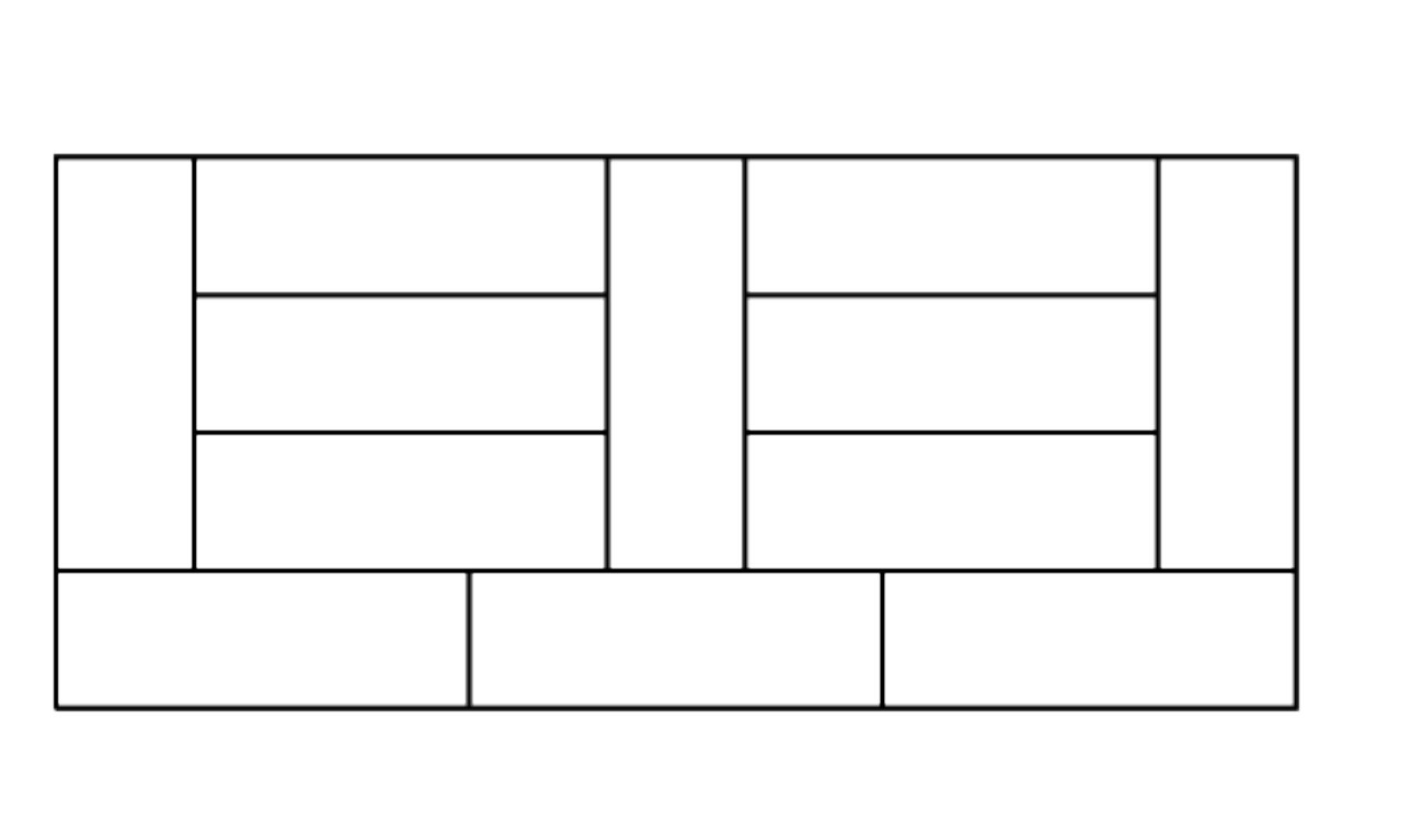}
\includegraphics[width=0.25\textwidth]{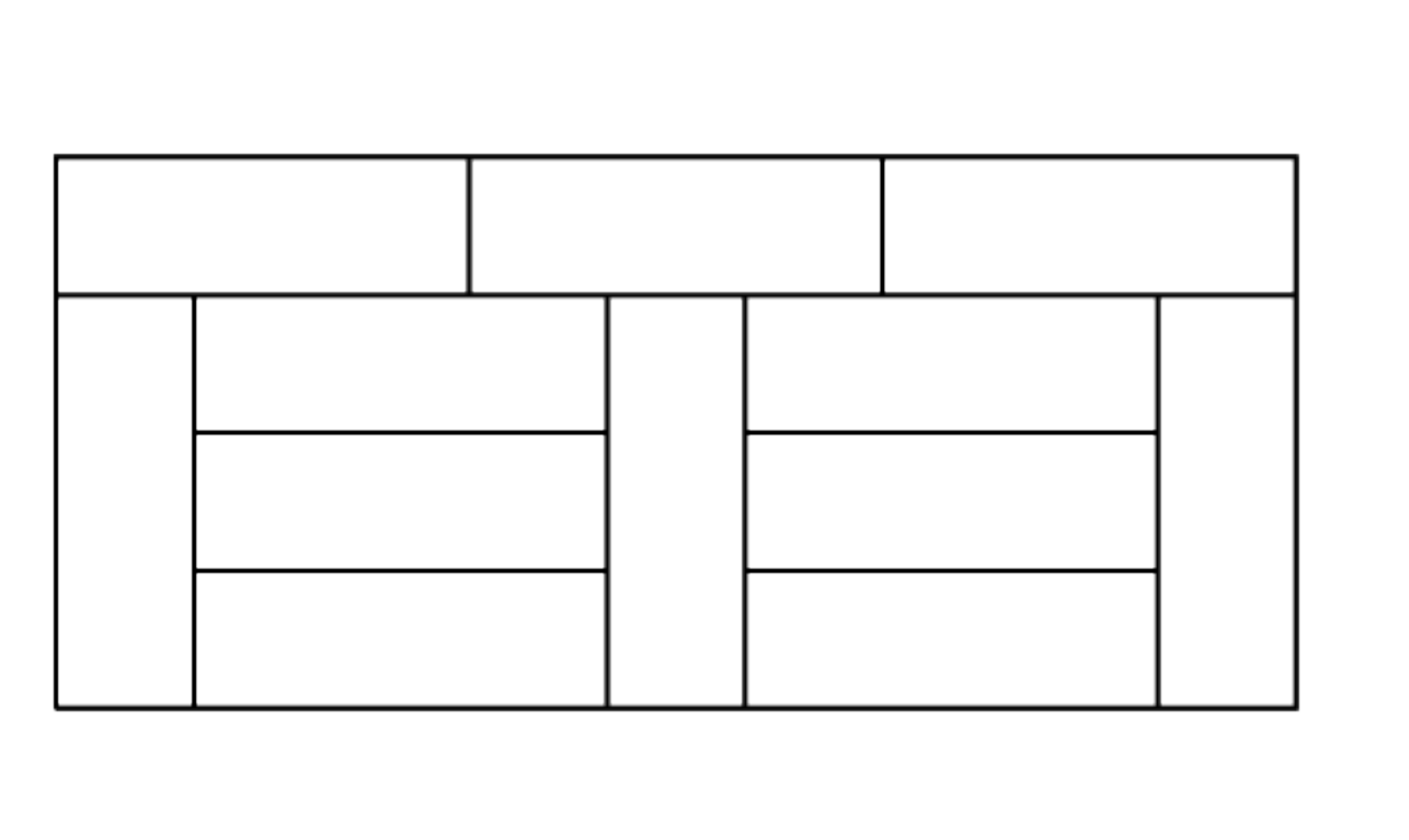}
\\
\hline
\includegraphics[width=0.25\textwidth]{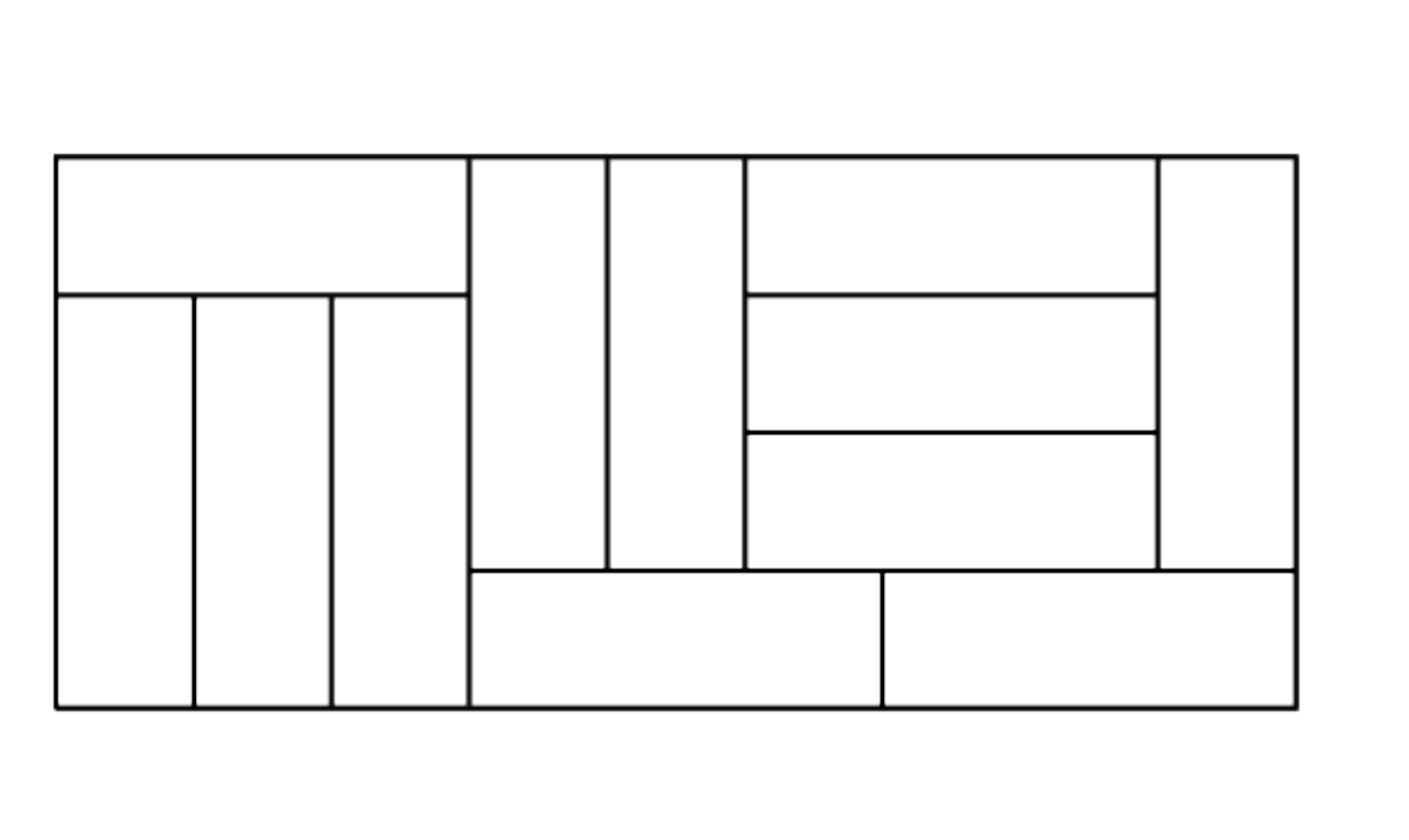}
\includegraphics[width=0.25\textwidth]{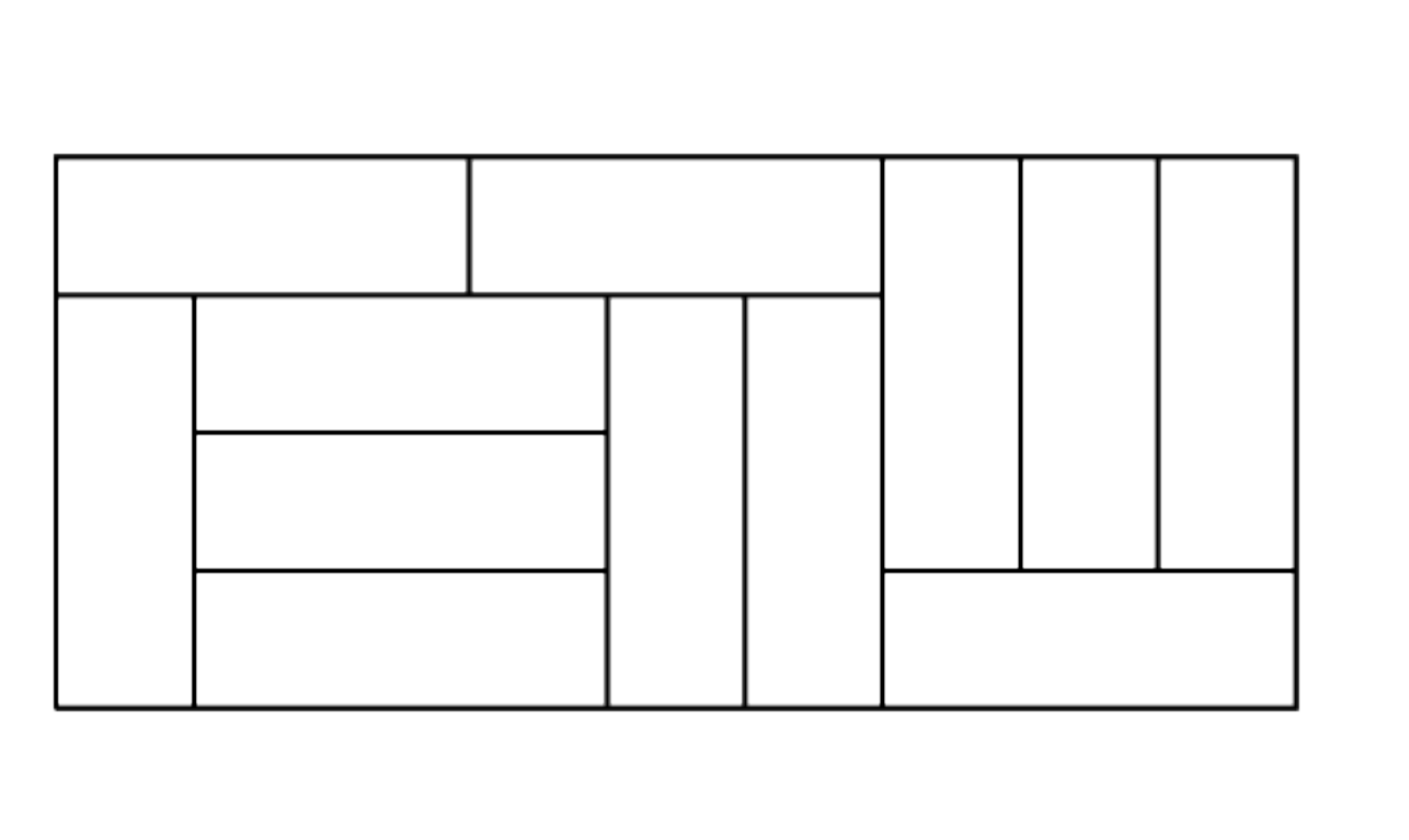}
\includegraphics[width=0.25\textwidth]{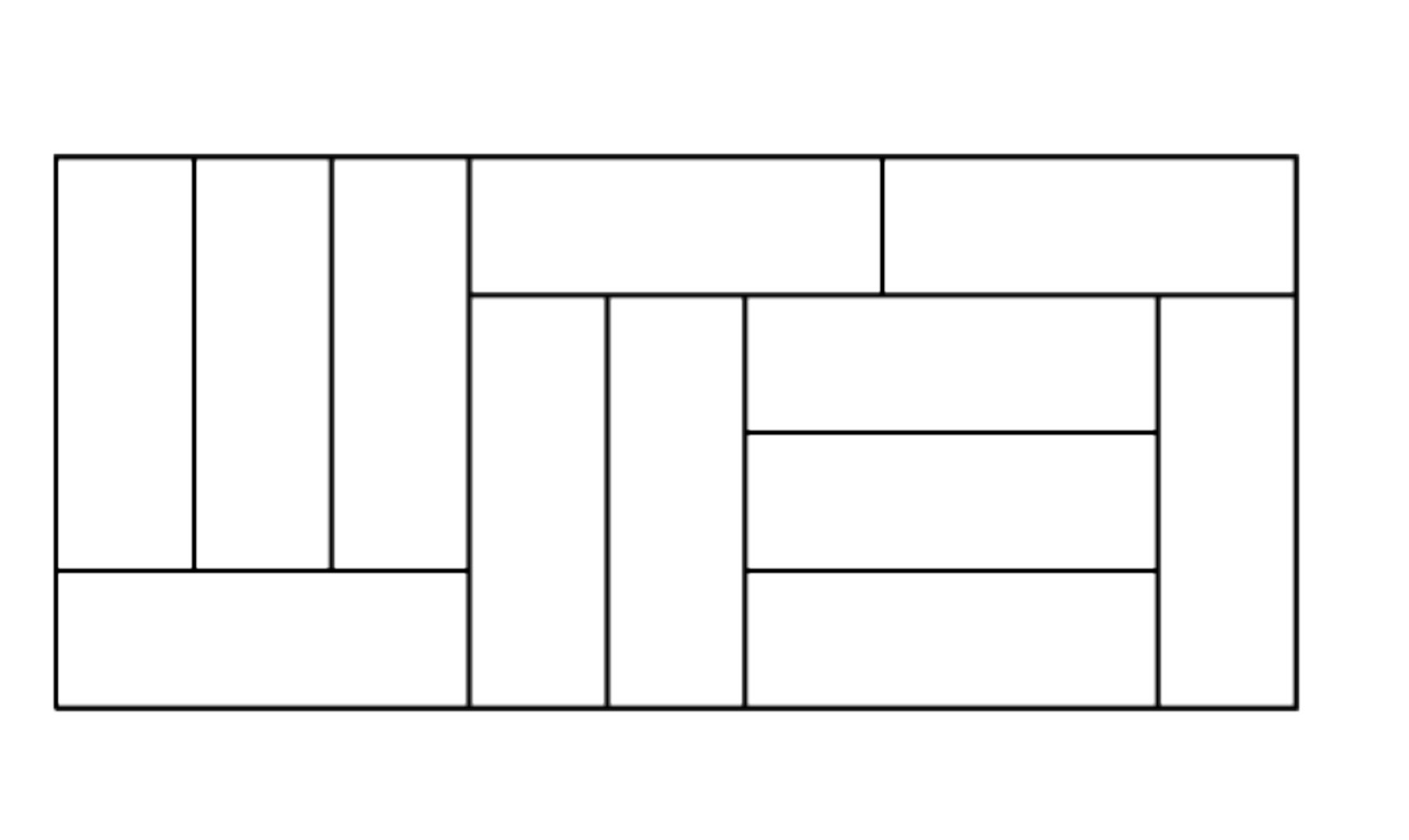}
\includegraphics[width=0.25\textwidth]{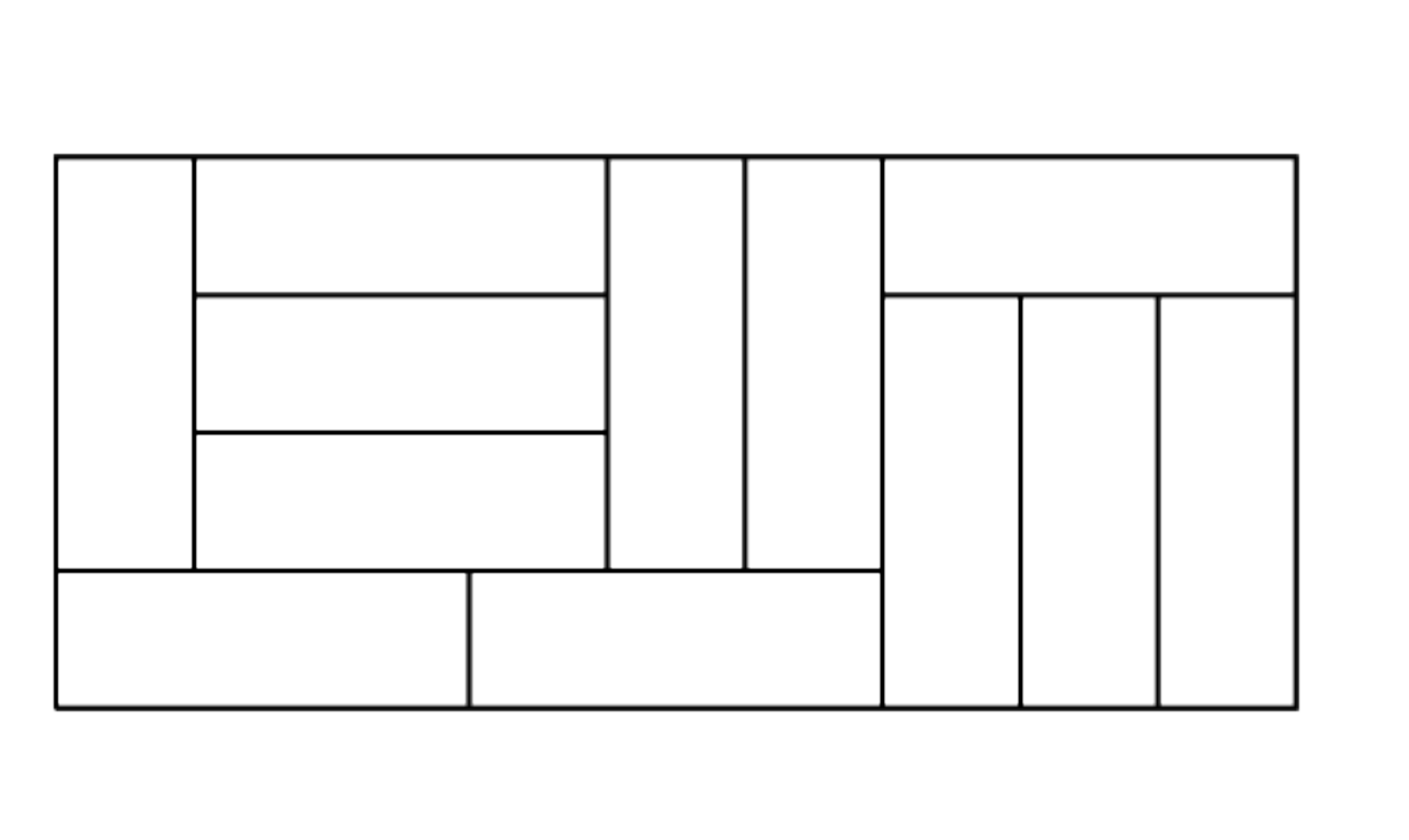}
\\
\hline
\includegraphics[width=0.25\textwidth]{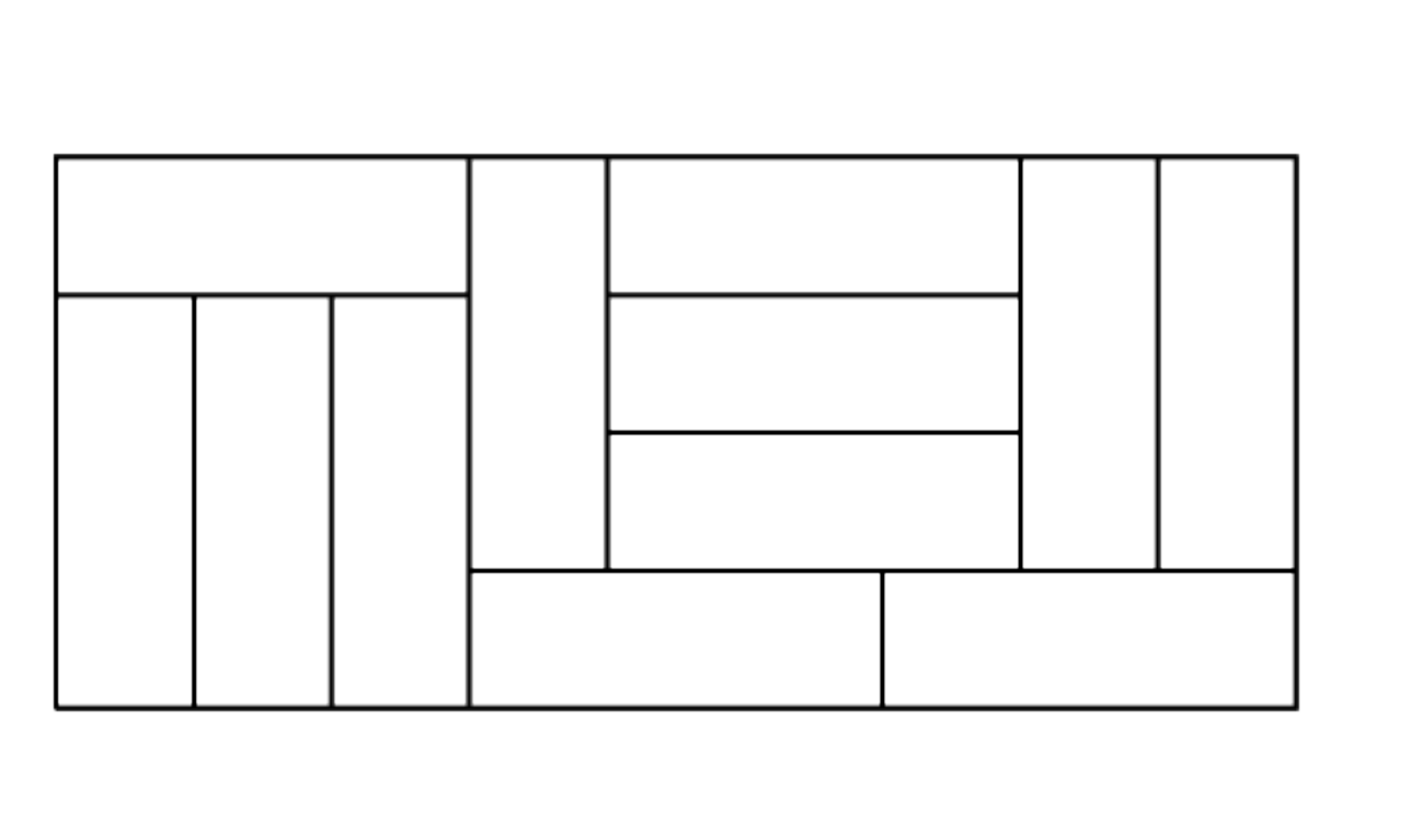}
\includegraphics[width=0.25\textwidth]{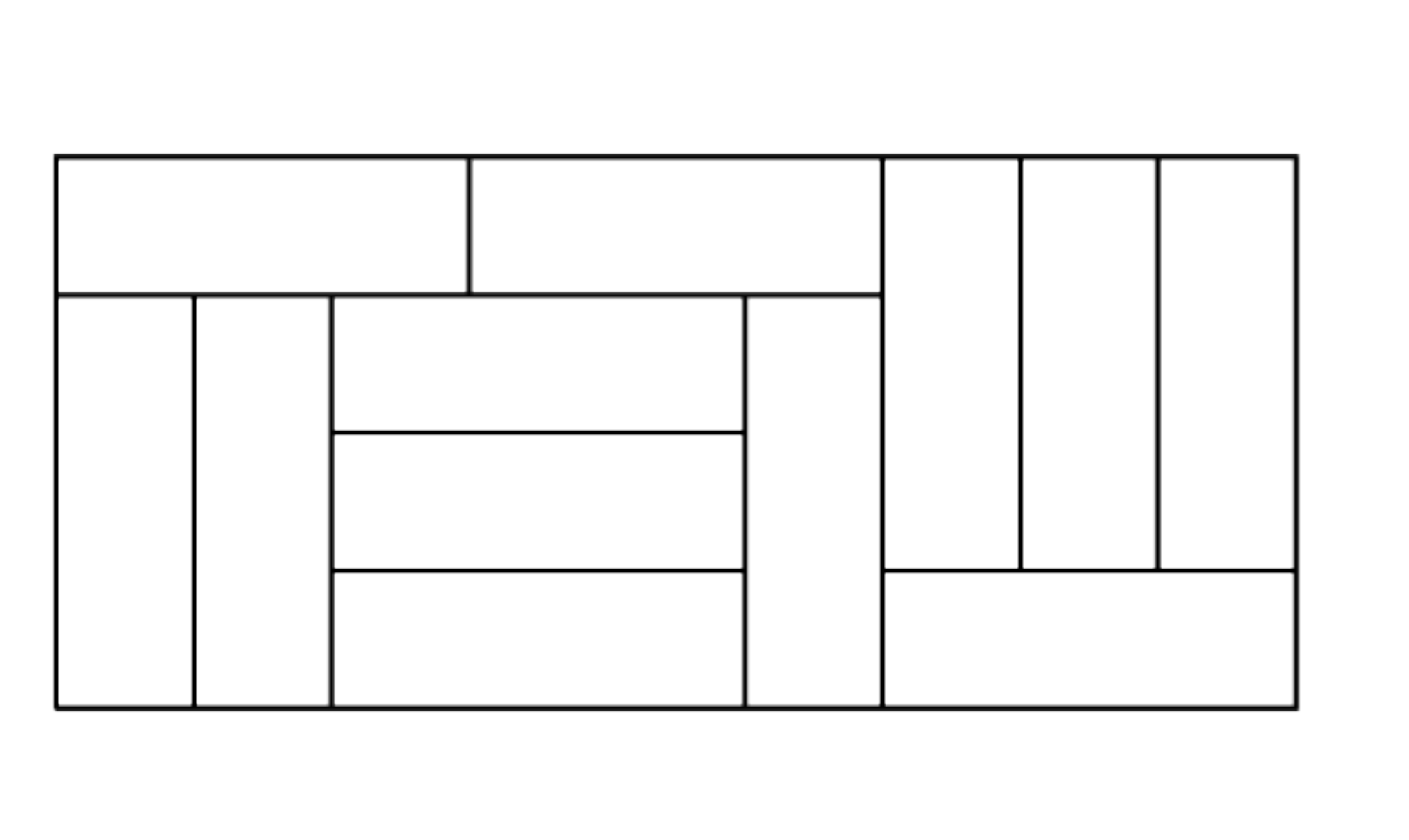}
\includegraphics[width=0.25\textwidth]{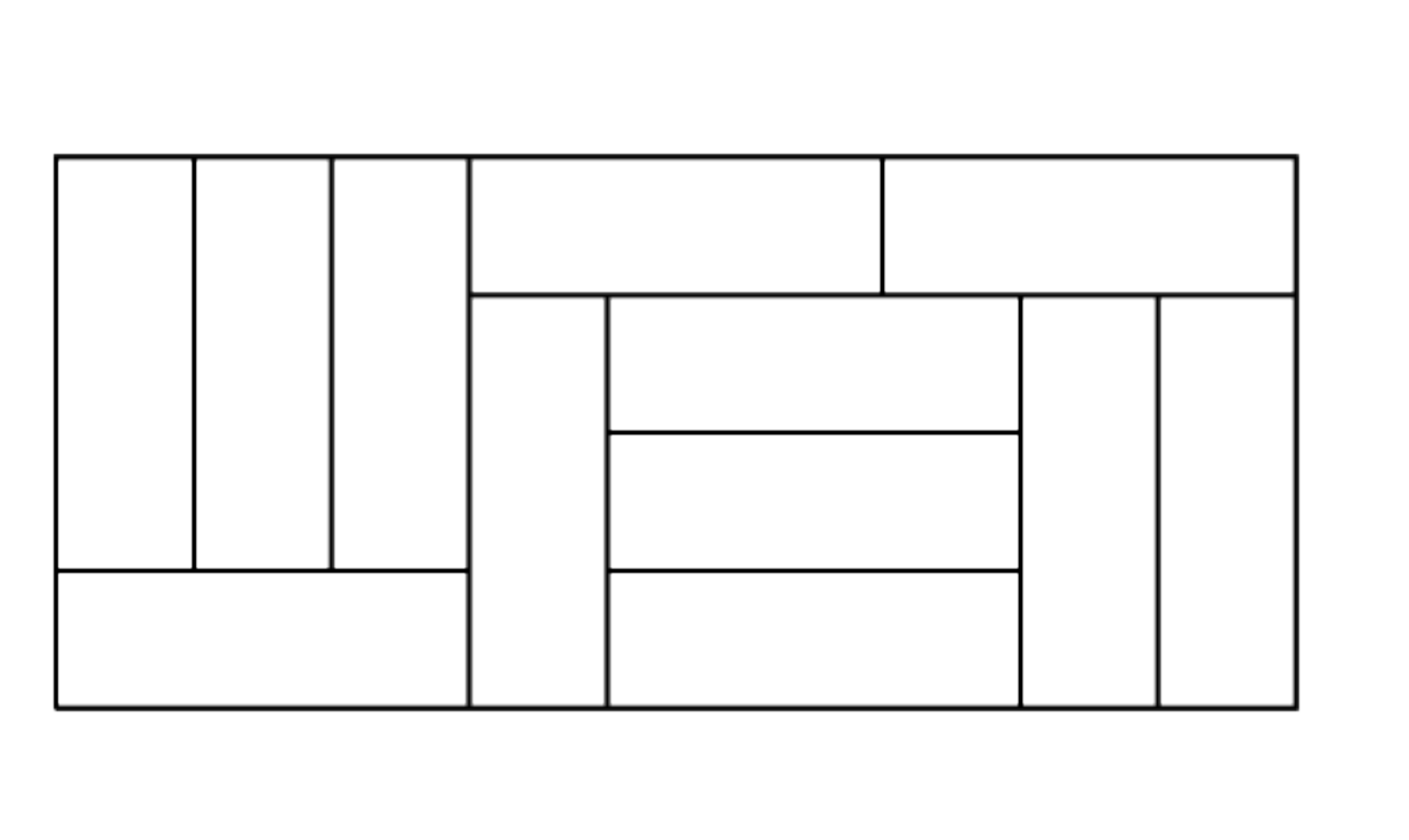}
\includegraphics[width=0.25\textwidth]{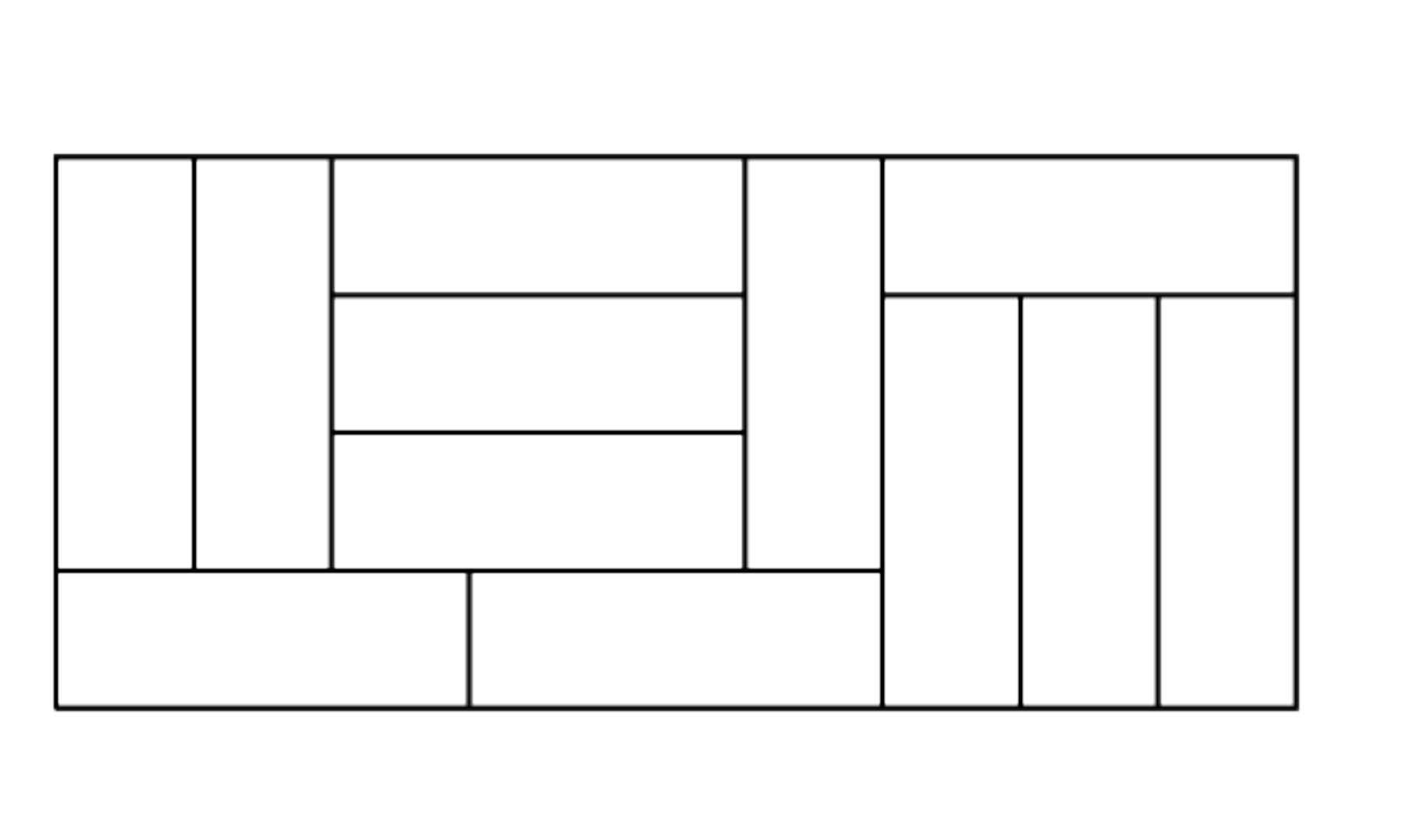}
\end{tabular}
\label{tab.run431}
\end{table}

\small
\begin{table}
\caption{Number $T(n,4)$ and $T_t(n,4)$ of tilings of $4\times n$ boards with 
$1\times 3$ tiles.
}
\begin{tabular}{rr|rrrrrrrrrrrrrrr}
$n$ & & 0 & 1 & 2 & 3 & 4 & 5 & 6 & 7 &8 &9 \\
\hline
3&3&3&0&0&0&0&0&0&0&0&0\\
6&13&6&6&0&1&0&0&0&0&0&0\\
9&57&12&24&16&0&4&0&1&0&0&0\\
12&249&26&66&84&40&16&12&0&4&0&1\\
15&1087&56&176&306&264&134&76&44&12&14&0\\
18&4745&120&452&970&1170&892&504&316&160&78&50\\
21&20713&258&1128&2852&4324&4388&3152&1996&1232&642&368\\
24&90417&554&2762&7986&14414&18070&16298&11784&7808&4810&2694\\
27&394691&1190&6660&21590&44916&66492&72344&61720&45156&30510&18948\\
30&1722917&2556&15868&56862&133338&226324&287668&286806&237554&175262&119330\\
33&7520929&5490&37440&146714&381640&727634&1054508&1206850&1135412&925770&684908\\
\end{tabular}
\label{tab.run4_3_1}
\end{table}
\normalsize

The row sums in Table \ref{tab.run5_3_1} are presumably 
generated by
\begin{conj} (Table \ref{tab.run5_3_1})
\begin{equation}
T(z,3) =
\frac{1-2z^3+z^6}{1-6z^3+3z^6-z^9}.
\end{equation}
\end{conj}

\small
\begin{table}
\caption{Number $T(n,5)$ and $T_t(n,5)$ of tilings of $5\times n$ boards with 
$1\times 3$ tiles.
}
\begin{tabular}{rr|rrrrrrrrrrrrrrr}
$n$ & & 0 & 1 & 2 & 3 & 4 & 5 & 6 & 7 &8 &9 & 10 & 11 \\
\hline
3&4&4&0&0&0&0&0&0&0&0&0&0&0\\
6&22&4&8&9&0&1&0&0&0&0&0&0&0\\
9&121&3&20&34&32&25&0&6&0&1&0&0&0\\
12&664&2&28&90&140&160&124&72&20&21&0&6&0\\
15&3643&2&32&164&388&653&760&678&456&265&112&86&16\\
18&19987&2&36&245&832&1854&2908&3748&3548&2756&1816&1077&552\\
21&109657&2&40&330&1488&4254&8592&14014&17636&18502&15692&11626&7508\\
24&601624&2&44&422&2344&8378&21096&41721&65708&86523&93068&86956&68992\\
27&3300760&2&48&522&3408&14695&44912&105920&200464&317485&417696&475112&461888\\
\end{tabular}
\label{tab.run5_3_1}
\end{table}
\normalsize

\small
\begin{table}
\caption{Number $T(n,6)$ and $T_t(n,6)$ of tilings of $6\times n$ boards with 
$1\times 3$ tiles.
}
\begin{tabular}{rr|rrrrrrrrrrrrrrr}
$n$ & & 0 & 1 & 2 & 3 & 4 & 5 & 6 & 7 &8 &9 \\
\hline
1&1&1&0&0&0&0&0&0&0&0&0\\
2&1&0&1&0&0&0&0&0&0&0&0\\
3&6&5&0&1&0&0&0&0&0&0&0\\
4&13&6&6&0&1&0&0&0&0&0&0\\
5&22&4&8&9&0&1&0&0&0&0&0\\
6&64&0&18&20&24&0&2&0&0&0&0\\
7&155&6&32&41&40&33&0&3&0&0&0\\
8&321&8&46&78&81&60&44&0&4&0&0\\
9&783&6&40&134&184&194&120&91&0&13&0\\
10&1888&0&58&212&416&464&360&220&128&10&18\\
11&4233&10&88&356&676&983&924&651&316&183&20\\
12&9912&16&112&488&1154&1746&2164&1852&1305&588&360\\
13&23494&10&104&634&1844&3574&4576&4744&3808&2332&1032\\
14&54177&0&146&858&3086&6308&9650&10522&9680&6976&4009\\
15&126019&18&216&1280&4548&10800&17676&23085&22848&19265&13036\\
16&295681&32&260&1660&6572&16972&32846&46652&54195&49582&38418\\
17&687690&18&248&2044&9184&28006&57740&93112&116144&121213&101656\\
18&1600185&0&342&2670&13568&43180&100706&172550&240028&270006&257483\\
19&3738332&34&504&3778&18808&65924&162320&313866&471344&589001&609660\\
20&8712992&64&584&4800&25474&95010&264062&545618&911664&1225476&1394156\\
21&20293761&34&568&5754&33940&141048&414236&944758&1678144&2469202&3014824\\
22&47337405&0&778&7402&47404&202330&649428&1559924&3023840&4749534&6313134\\
\end{tabular}
\end{table}
\normalsize

\small
\begin{table}
\caption{Number $T(n,7)$ and $T_t(n,7)$ of tilings of $7\times n$ boards with 
$1\times 3$ tiles.
}
\begin{tabular}{rr|rrrrrrrrrrrrrrr}
$n$ & & 0 & 1 & 2 & 3 & 4 & 5 & 6 & 7 &8 &9 \\
\hline
3&9&7&0&2&0&0&0&0&0&0&0\\
6&155&6&32&41&40&33&0&3&0&0&0\\
9&2861&23&100&374&552&558&572&366&168&131&0\\
12&52817&28&288&1301&3884&7303&8976&9658&8516&5610&3764\\
15&972557&55&612&3900&16356&43495&86916&129728&154264&158222&133548\\
18&17892281&140&1388&11760&54828&196629&520436&1057125&1723036&2301996&2624904\\
\end{tabular}
\end{table}
\normalsize

\small
\begin{table}
\caption{Number $T(n,8)$ and $T_t(n,8)$ of tilings of $8\times n$ boards with 
$1\times 3$ tiles.
}
\begin{tabular}{rr|rrrrrrrrrrrrrrr}
$n$ & & 0 & 1 & 2 & 3 & 4 & 5 & 6 & 7 &8 &9 \\
\hline
3&13&10&0&3&0&0&0&0&0&0&0&\\
6&321&8&46&78&81&60&44&0&4&0&0&\\
9&8133&8&92&471&1052&1580&1552&1298&1008&623&236&\\
12&204975&8&140&1026&4578&12340&22428&31240&34738&31358&25410&\\
15&5158223&10&128&2030&12040&50140&139904&295867&488740&656582&751376&\\
\end{tabular}
\end{table}
\normalsize

\small
\begin{table}
\caption{Number $T(n,9)$ and $T_t(n,9)$ of tilings of $9\times n$ boards with 
$1\times 3$ tiles.
}
\begin{tabular}{rr|rrrrrrrrrrrrrrr}
$n$ & & 0 & 1 & 2 & 3 & 4 & 5 & 6 & 7 &8 &9 \\
\hline
1&1&1&0&0&0&0&0&0&0&0&0&\\
2&1&0&0&1&0&0&0&0&0&0&0&\\
3&19&14&0&4&0&1&0&0&0&0&0&\\
4&57&12&24&16&0&4&0&1&0&0&0&\\
5&121&3&20&34&32&25&0&6&0&1&0&\\
6&783&6&40&134&184&194&120&91&0&13&0&\\
7&2861&23&100&374&552&558&572&366&168&131&0&\\
8&8133&8&92&471&1052&1580&1552&1298&1008&623&236&\\
9&37160&6&80&528&1832&4344&6432&7092&6016&4690&3040&\\
10&143419&22&284&1441&4796&10973&17900&23948&25048&21910&16172&\\
11&468816&46&408&2276&8040&21416&40512&60384&73908&77511&66884&\\
12&1876855&12&220&2062&10084&34134&83864&157514&229952&281691&292240&\\
13&7263468&35&444&3504&18432&69552&190316&402638&669380&917207&1058448&\\
14&25496863&86&1048&7892&38008&137588&387664&879237&1602712&2461860&3207432&\\
15&97187247&97&1100&8848&48832&204383&654700&1688610&3530980&6175615&9156644&\\
\end{tabular}
\end{table}
\normalsize

\clearpage
\subsection{Results (incongruent)}
Counts of tilings with $1\times 3$ tiles where only one representative of the
roto-reflected copies of each tiling is counted are shown
in Tables \ref{tab.run3i_3_1}--\ref{tab.run9i_3_1}.

Row sums $\bar T(n,3)$ of Table \ref{tab.run3i_3_1} appear to
be tabulated in \cite[A102543]{EIS} with row sums
\begin{conj} (Table \ref{tab.run3i_3_1})
\begin{equation}
\bar T(z,3)=
-z^3+\frac{1}{2}[
\frac{1}{1-z-z^3}+\frac{1+z+z^3}{1-z^2-z^6}
].
\end{equation}
\end{conj}
Up to a shift in the indices these numbers appeared already in column 0
of Table \ref{tab.run2i}.

Column $\bar T_0(n,3)$ seems to be \cite[A192928]{EIS}:
\begin{conj} (Table \ref{tab.run3i_3_1})
\begin{equation}
\bar T_0(z,3)=
-z^3+\frac{1}{2}[
 \frac{1+z^3}{1-z-z^4}
+ \frac{1+z+z^3+z^7}{1-z^2-z^8}
].
\end{equation}
\end{conj}

\begin{conj} (Table \ref{tab.run3i_3_1})
\begin{equation}
\bar T_2(z,3)=
\frac{z^2}{2}[
\frac{1-z}{(1-z-z^4)^2}
+\frac{1}{1-z-z^4}
+\frac{z^4}{1-z^2-z^8}
].
\end{equation}
\end{conj}

\small
\begin{table}
\caption{Number $\bar T(n,3)$ and $\bar T_t(n,3)$ of incongruent tilings of $3\times n$ boards with 
$1\times 3$ tiles.
}
\begin{tabular}{rr|rrrrrrrrrrrrrrr}
$n$ & & 0 & 1 & 2 & 3 & 4 & 5 & 6 & 7 &8 &9 & 10 & 11 & 12 & 13 & 14 \\
\hline
1&1&1&0&0&0&0&0&0&0&0&0&0&0&0&0&0\\
2&1&1&0&0&0&0&0&0&0&0&0&0&0&0&0&0\\
3&1&1&0&0&0&0&0&0&0&0&0&0&0&0&0&0\\
4&2&2&0&0&0&0&0&0&0&0&0&0&0&0&0&0\\
5&3&3&0&0&0&0&0&0&0&0&0&0&0&0&0&0\\
6&4&3&0&1&0&0&0&0&0&0&0&0&0&0&0&0\\
7&6&5&0&1&0&0&0&0&0&0&0&0&0&0&0&0\\
8&8&6&0&2&0&0&0&0&0&0&0&0&0&0&0&0\\
9&12&9&0&2&0&1&0&0&0&0&0&0&0&0&0&0\\
10&16&11&0&4&0&1&0&0&0&0&0&0&0&0&0&0\\
11&24&16&0&6&0&2&0&0&0&0&0&0&0&0&0&0\\
12&33&20&0&10&0&2&0&1&0&0&0&0&0&0&0&0\\
13&49&29&0&14&0&5&0&1&0&0&0&0&0&0&0&0\\
14&69&37&0&22&0&8&0&2&0&0&0&0&0&0&0&0\\
15&102&53&0&32&0&14&0&2&0&1&0&0&0&0&0&0\\
16&145&69&0&50&0&20&0&5&0&1&0&0&0&0&0&0\\
17&214&98&0&72&0&33&0&9&0&2&0&0&0&0&0&0\\
18&307&130&0&108&0&50&0&16&0&2&0&1&0&0&0&0\\
19&452&183&0&156&0&82&0&24&0&6&0&1&0&0&0&0\\
20&653&245&0&232&0&122&0&41&0&11&0&2&0&0&0&0\\
21&960&343&0&336&0&191&0&67&0&20&0&2&0&1&0&0\\
22&1393&463&0&493&0&286&0&114&0&30&0&6&0&1&0&0\\
23&2046&646&0&711&0&444&0&178&0&53&0&12&0&2&0&0\\
24&2978&877&0&1036&0&666&0&285&0&89&0&22&0&2&0&1\\
25&4371&1220&0&1492&0&1014&0&447&0&156&0&34&0&7&0&1\\
26&6376&1664&0&2161&0&1511&0&714&0&248&0&62&0&14&0&2\\
27&9354&2310&0&3103&0&2280&0&1114&0&408&0&110&0&26&0&2\\
28&13665&3161&0&4472&0&3390&0&1740&0&656&0&198&0&40&0&7\\
29&20041&4381&0&6408&0&5073&0&2684&0&1079&0&324&0&75&0&15\\
30&29307&6009&0&9201&0&7507&0&4158&0&1723&0&542&0&136&0&28\\
31&42972&8319&0&13155&0&11156&0&6382&0&2757&0&900&0&250&0&44\\
32&62884&11430&0&18828&0&16454&0&9791&0&4349&0&1514&0&414&0&85\\
33&92191&15811&0&26864&0&24320&0&14909&0&6902&0&2482&0&707&0&161\\
34&134974&21751&0&38349&0&35739&0&22694&0&10834&0&4048&0&1198&0&302\\
35&197858&30070&0&54615&0&52568&0&34374&0&17000&0&6540&0&2062&0&510\\
36&289772&41405&0&77788&0&76994&0&51998&0&26442&0&10594&0&3440&0&883\\
37&424746&57216&0&110592&0&112799&0&78326&0&41129&0&17008&0&5713&0&1535\\
38&622198&78836&0&157210&0&164707&0&117802&0&63605&0&27200&0&9394&0&2692
\end{tabular}
\label{tab.run3i_3_1}
\end{table}
\normalsize

\small
\begin{table}
\caption{Number $\bar T(n,4)$ and $\bar T_t(n,4)$ of incongruent tilings of $4\times n$ boards with 
$1\times 3$ tiles.
}
\begin{tabular}{rr|rrrrrrrrrrrrrrr}
$n$ & & 0 & 1 & 2 & 3 & 4 & 5 & 6 & 7 &8 &9 & 10 \\
\hline
3&2&2&0&0&0&0&0&0&0&0&0&0\\
6&5&2&2&0&1&0&0&0&0&0&0&0\\
9&18&4&6&6&0&1&0&1&0&0&0&0\\
12&69&8&18&22&11&4&4&0&1&0&1&0\\
15&287&16&44&83&66&38&19&11&3&5&0&1\\
18&1215&33&116&247&297&226&132&79&42&20&13&3\\
21&5244&69&282&729&1081&1119&788&511&308&165&92&54\\
24&22729&145&697&2010&3617&4534&4097&2953&1968&1205&682&394\\
27&98959&307&1665&5438&11229&16698&18086&15503&11289&7666&4737&2826\\
30&431273&653&3981&14253&33372&56644&71993&71751&59465&43839&29889&18918\\
33&1881481&1393&9360&36778&95410&182136&263627&302018&283853&231688&171227&117565\\
36&8210019&2978&21937&93182&265472&560002&908584&1173497&1244674&1127399&909441&673113\\
\end{tabular}
\label{tab.run4i_3_1}
\end{table}
\normalsize

\small
\begin{table}
\caption{Number $\bar T(n,5)$ and $\bar T_t(n,5)$ of incongruent tilings of $5\times n$ boards with 
$1\times 3$ tiles.
}
\begin{tabular}{rr|rrrrrrrrrrrrrrr}
$n$ & & 0 & 1 & 2 & 3 & 4 & 5 & 6 & 7 &8 &9 & 10 & 11 & 12 \\
\hline
3&3&3&0&0&0&0&0&0&0&0&0&0&0&0\\
6&9&2&2&4&0&1&0&0&0&0&0&0&0&0\\
9&42&2&5&13&8&11&0&2&0&1&0&0&0&0\\
12&192&1&7&28&35&51&31&23&5&8&0&2&0&1\\
15&996&1&8&47&97&187&190&199&114&83&28&25&4&10\\
18&5206&1&9&66&208&485&727&1003&887&749&454&302&138&100\\
21&28091&1&10&89&372&1092&2148&3613&4409&4839&3923&3067&1877&1227\\
24&152212&1&11&111&586&2118&5274&10535&16427&21992&23267&22282&17248&12726\\
27&830974&1&12&138&852&3714&11228&26614&50116&79911&104424&120138&115472&101054\\
30&4543764&1&13&164&1174&5951&21405&59457&131820&245085&381516&515789&596135&615285\\
\end{tabular}
\end{table}
\normalsize

\small
\begin{table}
\caption{Number $\bar T(n,6)$ and $\bar T_t(n,6)$ of incongruent tilings of  $6\times n$ boards with 
$1\times 3$ tiles.
}
\begin{tabular}{rr|rrrrrrrrrrrrrrr}
$n$ & & 0 & 1 & 2 & 3 & 4 & 5 & 6 & 7 &8 &9 & 10 \\
\hline
1&1&1&0&0&0&0&0&0&0&0&0&0\\
2&1&0&1&0&0&0&0&0&0&0&0&0\\
3&4&3&0&1&0&0&0&0&0&0&0&0\\
4&5&2&2&0&1&0&0&0&0&0&0&0\\
5&9&2&2&4&0&1&0&0&0&0&0&0\\
6&11&0&3&3&4&0&1&0&0&0&0&0\\
7&47&3&8&13&10&11&0&2&0&0&0&0\\
8&91&2&13&20&24&16&13&0&3&0&0&0\\
9&219&3&10&38&46&56&30&30&0&5&0&1\\
10&494&0&16&54&107&117&97&57&35&3&7&0\\
11&1106&4&22&99&169&256&231&177&79&54&5&8\\
12&2533&4&28&124&293&441&550&465&341&150&98&14\\
13&5978&4&26&170&461&916&1144&1211&952&607&258&158\\
14&13670&0&37&217&779&1583&2431&2638&2443&1747&1030&414\\
15&31765&7&54&335&1137&2734&4419&5824&5712&4878&3259&1922\\
16&74194&8&67&416&1656&4252&8239&11678&13590&12410&9650&6136\\
17&172508&7&62&529&2296&7071&14435&23380&29036&30431&25414&18417\\
18&400688&0&88&669&3403&10808&25223&43164&60073&67537&64465&50706\\
19&935891&11&126&973&4702&16574&40580&78666&117836&147498&152415&135777\\
20&2179732&16&146&1203&6382&23774&66072&136452&228070&306440&348724&333975\\
21&5076530&11&142&1471&8485&35410&103559&236504&419536&617812&753706&789385\\
22&11837691&0&195&1855&11871&50605&162461&390059&756195&1187512&1578687&1778376\\
\end{tabular}
\end{table}
\normalsize

\small
\begin{table}
\caption{Number $\bar T(n,7)$ and $\bar T_t(n,7)$ of incongruent tilings of $7\times n$ boards with 
$1\times 3$ tiles.
}
\begin{tabular}{rr|rrrrrrrrrrrrrrr}
$n$ & & 0 & 1 & 2 & 3 & 4 & 5 & 6 & 7 &8 &9 & 10 \\
\hline
3&6&5&0&1&0&0&0&0&0&0&0&0\\
6&47&3&8&13&10&11&0&2&0&0&0&0\\
9&769&10&25&105&138&155&143&103&42&42&0&5\\
12&13354&10&72&343&971&1856&2244&2447&2129&1431&941&535\\
15&244096&23&153&1020&4089&11002&21729&32621&38566&39749&33387&25076\\
18&4475868&45&347&2991&13707&49336&130109&264648&430759&575998&656226&647062\\
\end{tabular}
\end{table}
\normalsize

\small
\begin{table}
\caption{Number $\bar T(n,8)$ and $\bar T_t(n,8)$ of incongruent tilings of $8\times n$ boards with 
$1\times 3$ tiles.
}
\begin{tabular}{rr|rrrrrrrrrrrrrrr}
$n$ & & 0 & 1 & 2 & 3 & 4 & 5 & 6 & 7 &8 &9 & 10 & 11 \\
\hline
3&8&6&0&2&0&0&0&0&0&0&0&0&0\\
6&91&2&13&20&24&16&13&0&3&0&0&0&0\\
9&2126&4&23&130&263&417&388&354&252&170&59&58&0\\
12&51493&2&37&258&1153&3091&5637&7819&8734&7849&6398&4582&2831\\
15&1291743&4&32&529&3010&12637&34976&74236&122185&164546&187844&185619&160116\\
\end{tabular}
\end{table}
\normalsize

\small
\begin{table}
\caption{Number $\bar T(n,9)$ and $\bar T_t(n,9)$ of incongruent tilings of $9\times n$ boards with 
$1\times 3$ tiles.
}
\begin{tabular}{rr|rrrrrrrrrrrrrrr}
$n$ & & 0 & 1 & 2 & 3 & 4 & 5 & 6 & 7 &8 &9 & 10 \\
\hline
1&1&1&0&0&0&0&0&0&0&0&0&0\\
2&1&0&0&1&0&0&0&0&0&0&0&0\\
3&12&9&0&2&0&1&0&0&0&0&0&0\\
4&18&4&6&6&0&1&0&1&0&0&0&0\\
5&42&2&5&13&8&11&0&2&0&1&0&0\\
6&219&3&10&38&46&56&30&30&0&5&0&1\\
7&769&10&25&105&138&155&143&103&42&42&0&5\\
8&2126&4&23&130&263&417&388&354&252&170&59&58\\
9&4808&2&10&76&229&579&804&926&752&628&380&258\\
10&36250&7&71&375&1199&2791&4475&6083&6262&5571&4043&2648\\
11&118396&19&102&610&2010&5492&10128&15323&18477&19668&16721&12599\\
12&471234&5&55&542&2521&8640&20966&39641&57488&70824&73060&66147\\
13&1820940&14&111&915&4608&17634&47579&101220&167345&230236&264612&267931\\
14&6383748&28&262&2022&9502&34623&96916&220504&400678&616756&801858&904085\\
\end{tabular}
\label{tab.run9i_3_1}
\end{table}
\normalsize

\clearpage

\section{Tiling with $1\times 4$ tiles}\label{sec.14}
Tilings with $1\times 4$ tiles are
represented by Tables \ref{tab.run3_4_1}--\ref{tab.run9i_4_1}.
\subsection{Results (full count)}

Table \ref{tab.run3_4_1} is simple because the floor width $m=3$ is
too small to allow ``vertical'' placements of the tiles. So there
is only one tiling with a ``ferromagnetic'' alignment of all tiles 
and a maximum number of points on the floor where 4 tiles meet.

\small
\begin{table}
\caption{Number $T(n,3)$ and $T_t(n,3)$ of tilings of $3\times n$ boards with 
$1\times 4$ tiles.
}
\begin{tabular}{rr|rrrrrrrrrrrrrrr}
$n$ & & 0 & 1 & 2 & 3 & 4 & 5 & 6 & 7 &8 &9 & 10 & 11 & 12 & 13 & 14 \\
\hline
4&1&1&0&0&0&0&0&0&0&0&0&0&0&0&0&0\\
8&1&0&0&1&0&0&0&0&0&0&0&0&0&0&0&0\\
12&1&0&0&0&0&1&0&0&0&0&0&0&0&0&0&0\\
16&1&0&0&0&0&0&0&1&0&0&0&0&0&0&0&0\\
20&1&0&0&0&0&0&0&0&0&1&0&0&0&0&0&0\\
24&1&0&0&0&0&0&0&0&0&0&0&1&0&0&0&0\\
28&1&0&0&0&0&0&0&0&0&0&0&0&0&1&0&0\\
32&1&0&0&0&0&0&0&0&0&0&0&0&0&0&0&1\\
\end{tabular}
\label{tab.run3_4_1}
\end{table}
\normalsize

Row sums $T(n,4)$ of Table \ref{tab.run4_4_1}
are \cite[A003269]{EIS} and column $T_0(n,4)$
is \cite[A003520]{EIS}.
For two further columns we propose
\begin{conj} (Table \ref{tab.run4_4_1})
\begin{equation}
T_3(z,4) =
\frac{z^8}{ (1-z+z^2)^2(1-z^2-z^3)^2 }.
\end{equation}
\end{conj}
\begin{conj} (Table \ref{tab.run4_4_1})
\begin{equation}
T_6(z,4) =
z^{12}\frac{1-z}{ (1-z+z^2)^3(1-z^2-z^3)^3 }.
\end{equation}
\end{conj}

\small
\begin{table}
\caption{Number $T(n,4)$ and $T_t(n,4)$ of tilings of $4\times n$ boards with 
$1\times 4$ tiles.
}
\begin{tabular}{rr|rrrrrrrrrrrrrrr}
$n$ & & 0 & 1 & 2 & 3 & 4 & 5 & 6 & 7 &8 &9 & 10 & 11 & 12 & 13 & 14 \\
\hline
1&1&1&0&0&0&0&0&0&0&0&0&0&0&0&0&0\\
2&1&1&0&0&0&0&0&0&0&0&0&0&0&0&0&0\\
3&1&1&0&0&0&0&0&0&0&0&0&0&0&0&0&0\\
4&2&2&0&0&0&0&0&0&0&0&0&0&0&0&0&0\\
5&3&3&0&0&0&0&0&0&0&0&0&0&0&0&0&0\\
6&4&4&0&0&0&0&0&0&0&0&0&0&0&0&0&0\\
7&5&5&0&0&0&0&0&0&0&0&0&0&0&0&0&0\\
8&7&6&0&0&1&0&0&0&0&0&0&0&0&0&0&0\\
9&10&8&0&0&2&0&0&0&0&0&0&0&0&0&0&0\\
10&14&11&0&0&3&0&0&0&0&0&0&0&0&0&0&0\\
11&19&15&0&0&4&0&0&0&0&0&0&0&0&0&0&0\\
12&26&20&0&0&5&0&0&1&0&0&0&0&0&0&0&0\\
13&36&26&0&0&8&0&0&2&0&0&0&0&0&0&0&0\\
14&50&34&0&0&13&0&0&3&0&0&0&0&0&0&0&0\\
15&69&45&0&0&20&0&0&4&0&0&0&0&0&0&0&0\\
16&95&60&0&0&29&0&0&5&0&0&1&0&0&0&0&0\\
17&131&80&0&0&40&0&0&9&0&0&2&0&0&0&0&0\\
18&181&106&0&0&56&0&0&16&0&0&3&0&0&0&0&0\\
19&250&140&0&0&80&0&0&26&0&0&4&0&0&0&0&0\\
20&345&185&0&0&115&0&0&39&0&0&5&0&0&1&0&0\\
21&476&245&0&0&164&0&0&55&0&0&10&0&0&2&0&0\\
22&657&325&0&0&230&0&0&80&0&0&19&0&0&3&0&0\\
23&907&431&0&0&320&0&0&120&0&0&32&0&0&4&0&0\\
24&1252&571&0&0&445&0&0&181&0&0&49&0&0&5&0&0\\
25&1728&756&0&0&620&0&0&269&0&0&70&0&0&11&0&0\\
26&2385&1001&0&0&864&0&0&390&0&0&105&0&0&22&0&0\\
27&3292&1326&0&0&1200&0&0&560&0&0&164&0&0&38&0&0\\
28&4544&1757&0&0&1660&0&0&805&0&0&257&0&0&59&0&0\\
29&6272&2328&0&0&2290&0&0&1161&0&0&394&0&0&85&0&0\\
30&8657&3084&0&0&3155&0&0&1674&0&0&585&0&0&131&0&0\\
31&11949&4085&0&0&4344&0&0&2400&0&0&860&0&0&212&0&0\\
32&16493&5411&0&0&5975&0&0&3420&0&0&1269&0&0&343&0&0\\
33&22765&7168&0&0&8206&0&0&4855&0&0&1882&0&0&539&0&0\\
34&31422&9496&0&0&11252&0&0&6881&0&0&2789&0&0&815&0&0\\
35&43371&12580&0&0&15408&0&0&9744&0&0&4100&0&0&1221&0&0\\
36&59864&16665&0&0&21078&0&0&13775&0&0&5980&0&0&1842&0&0\\
37&82629&22076&0&0&28810&0&0&19426&0&0&8684&0&0&2798&0&0\\
38&114051&29244&0&0&39344&0&0&27327&0&0&12592&0&0&4244&0&0\\
39&157422&38740&0&0&53680&0&0&38364&0&0&18244&0&0&6370&0&0\\
40&217286&51320&0&0&73173&0&0&53778&0&0&26375&0&0&9471&0&0\\
41&299915&67985&0&0&99662&0&0&75285&0&0&38006&0&0&14017&0&0\\
42&413966&90061&0&0&135640&0&0&105245&0&0&54591&0&0&20723&0&0\\
43&571388&119305&0&0&184480&0&0&146908&0&0&78220&0&0&30619&0&0\\
44&788674&158045&0&0&250740&0&0&204765&0&0&111878&0&0&45120&0&0\\
45&1088589&209365&0&0&340578&0&0&285032&0&0&159760&0&0&66222&0&0\\
46&1502555&277350&0&0&462316&0&0&396295&0&0&227726&0&0&96824&0&0\\
47&2073943&367411&0&0&627200&0&0&550380&0&0&323980&0&0&141176&0&0\\
48&2862617&486716&0&0&850420&0&0&763548&0&0&460057&0&0&205453&0&0\\
49&3951206&644761&0&0&1152480&0&0&1058151&0&0&652202&0&0&298455&0&0\\
50&5453761&854126&0&0&1561043&0&0&1464921&0&0&923225&0&0&432643&0&0\\
51&7527704&1131476&0&0&2113420&0&0&2026100&0&0&1305036&0&0&625721&0&0\\
52&10390321&1498887&0&0&2859925&0&0&2799700&0&0&1842184&0&0&902974&0&0\\
\end{tabular}
\label{tab.run4_4_1}
\end{table}
\normalsize

For row sums in Table \ref{tab.run5_4_1} we conjecture
\begin{conj} (Table \ref{tab.run5_4_1})
\begin{equation}
T(z,5)=
\frac{1-3z^4+3z^8-z^{12}}{ 1-6z^4+6z^8-4z^{12}+z^{16} }.
\end{equation}
\end{conj}
\begin{conj} (Table \ref{tab.run5_4_1})
\begin{equation}
T_0(z,5)=
-1+z^4+\frac{2}{1-z^4-3z^8-3z^{12}-z^{16}}.
\end{equation}
\end{conj}

\small
\begin{table}
\caption{Number $T(n,5)$ and $T_t(n,5)$ of tilings of $5\times n$ boards with 
$1\times 4$ tiles.
}
\begin{tabular}{rr|rrrrrrrrrrrrrrr}
$n$ & & 0 & 1 & 2 & 3 & 4 & 5 & 6 & 7 &8 \\
\hline
4&3&3&0&0&0&0&0&0&0&0\\
8&15&8&6&0&0&1&0&0&0&0\\
12&75&20&32&12&6&0&4&0&0&1\\
16&371&52&112&96&48&24&16&18&0&0\\
20&1833&138&368&464&340&192&112&110&56&12\\
24&9057&362&1168&1872&1856&1288&824&640&490&200\\
28&44753&952&3592&6928&8502&7384&5292&3992&3120&1966\\
32&221137&2504&10816&24248&35236&36772&30492&23778&18864&13864\\
36&1092699&6584&32048&81584&136456&165812&158516&133522&109200&86104\\
40&5399327&17314&93754&266592&502880&695648&755996&700938&603264&500044\\
44&26679563&45530&271468&851582&1784664&2763024&3366464&3444380&3169744&2758358\\
48&131831075&119728&779418&2671108&6147774&10513564&14189044&15958800&15823444&14535462\\

\end{tabular}
\label{tab.run5_4_1}
\end{table}
\normalsize

For row sums in Table \ref{tab.run6_4_1} we conjecture
\begin{conj} (Table \ref{tab.run6_4_1})
\begin{equation}
T(z,6)=
\frac{1-3z^4+3z^8-z^{12}}{1-7z^4+6z^8-4z^{12}+z^{16}}.
\end{equation}
\end{conj}

\small
\begin{table}
\caption{Number $T(n,6)$ and $T_t(n,6)$ of tilings of $6\times n$ boards with 
$1\times 4$ tiles.
}
\begin{tabular}{rr|rrrrrrrrrrrrrrr}
$n$ & & 0 & 1 & 2 & 3 & 4 & 5 & 6 & 7 &8 &9 & 10 \\
\hline
2&0&0&0&0&0&0&0&0&0&0&0&0\\
4&4&4&0&0&0&0&0&0&0&0&0&0\\
6&0&0&0&0&0&0&0&0&0&0&0&0\\
8&25&5&10&9&0&0&1&0&0&0&0&0\\
10&0&0&0&0&0&0&0&0&0&0&0&0\\
12&154&5&28&44&43&21&6&0&6&0&0&1\\
14&0&0&0&0&0&0&0&0&0&0&0&0\\
16&943&3&46&155&198&214&164&75&34&20&27&0\\
18&0&0&0&0&0&0&0&0&0&0&0&0\\
20&5773&2&56&342&695&1084&1120&1015&604&323&207&156\\
22&0&0&0&0&0&0&0&0&0&0&0&0\\
24&35344&2&62&565&1860&3795&5342&6427&5674&4360&2609&1763\\
26&0&0&0&0&0&0&0&0&0&0&0&0\\
28&216388&2&68&794&3892&10697&19470&29080&34291&33610&27010&19746\\
30&0&0&0&0&0&0&0&0&0&0&0&0\\
32&1324801&2&74&1027&6788&25048&58416&106132&154688&188377&191619&168341\\
34&0&0&0&0&0&0&0&0&0&0&0&0\\
36&8110882&2&80&1272&10492&50091&149232&328062&572908&839481&1031920&1094524&1018003&849414&656866&480313\\
38&0&0&0&0&0&0&0&0&0&0&0&0&0&0&0&0\\
40&49657576&2&86&1534&15004&88602&331456&882617&1824950&3141422&4555427&5699764&6226786&6056113&5323270&4319508

\end{tabular}
\label{tab.run6_4_1}
\end{table}
\normalsize

For row sums in Table \ref{tab.run7_4_1} we conjecture
\begin{conj} (Table \ref{tab.run7_4_1})
\begin{equation}
T(z,7)=
\frac{1-z^4}{2}
[
\frac{1}{1-2z^2-2z^4+z^8}
+\frac{1}{1+2z^2-2z^4+z^8}
].
\end{equation}
\end{conj}

\small
\begin{table}
\caption{Number $T(n,7)$ and $T_t(n,7)$ of tilings of $7\times n$ boards with 
$1\times 4$ tiles.
}
\begin{tabular}{rr|rrrrrrrrrrrrrrr}
$n$ & & 0 & 1 & 2 & 3 & 4 & 5 & 6 & 7 &8 &9 & 10 & 11 & 12 \\
\hline
4&5&5&0&0&0&0&0&0&0&0&0&0&0&0\\
8&37&2&10&12&12&0&0&1&0&0&0&0&0&0\\
12&269&2&4&38&56&68&54&32&6&0&8&0&0&1\\
16&1949&2&4&40&98&258&386&356&348&228&114&44&24&38\\
20&14121&2&4&48&124&448&1036&1682&2294&2476&2150&1638&940&512\\
24&102313&2&4&56&152&666&1786&4154&8000&11640&14602&16012&14476&11104\\
28&741305&2&4&64&180&916&2696&7608&18120&34864&57616&80422&97260&104246\\
32&5371097&2&4&72&208&1198&3766&12238&33214&76820&157410&273472&412536&553940\\
36&38916077&2&4&80&236&1512&4996&18172&54156&143044&339912&700280&1277312&2062920\\

\end{tabular}
\label{tab.run7_4_1}
\end{table}
\normalsize

\small
\begin{table}
\caption{Number $T(n,8)$ and $T_t(n,8)$ of tilings of $8\times n$ boards with 
$1\times 4$ tiles.
}
\begin{tabular}{rr|rrrrrrrrrrrrrrr}
$n$ & & 0 & 1 & 2 & 3 & 4 & 5 & 6 & 7 &8 &9 \\
\hline
1&1&1&0&0&0&0&0&0&0&0&0\\
2&1&0&1&0&0&0&0&0&0&0&0\\
3&1&0&0&1&0&0&0&0&0&0&0\\
4&7&6&0&0&1&0&0&0&0&0&0\\
5&15&8&6&0&0&1&0&0&0&0&0\\
6&25&5&10&9&0&0&1&0&0&0&0\\
7&37&2&10&12&12&0&0&1&0&0&0\\
8&100&0&10&30&28&30&0&0&2&0&0\\
9&229&13&30&53&40&52&36&2&0&3&0\\
10&454&20&61&87&95&62&74&47&4&0&4\\
11&811&11&80&147&168&154&84&98&58&6&0\\
12&1732&6&65&182&291&344&314&198&187&114&14\\
13&3777&2&56&302&534&741&760&573&296&310&138\\
14&7858&33&132&502&928&1274&1511&1347&938&470&430\\
15&15339&62&240&762&1588&2185&2616&2621&2198&1432&674\\
16&31273&31&342&1060&2354&3721&4467&4930&4637&3962&2541\\
17&65536&12&288&1286&3268&6265&8842&10044&10194&8781&7080\\
18&136600&10&224&1819&5076&10549&16351&20314&20956&19151&15498\\
19&276535&91&490&2715&8198&17134&28388&37836&42196&39726&33754\\
20&562728&182&895&3995&12713&27487&46732&65807&78456&80726&72305\\
21&1159942&95&1288&5404&17816&41810&74992&112822&144618&162166&158642\\
22&2400783&22&1107&6600&23252&62848&123834&198059&269447&316179&330597\\
23&4918159&30&808&8761&33098&95131&202666&344897&492644&607682&657528\\
24&10052140&267&1664&12278&50305&144614&325405&582052&874124&1127120&1277621\\
25&20627526&528&3178&17688&74316&219790&505436&949104&1490540&2029691&2427958\\
26&42480474&287&4568&23867&101648&318272&766146&1516726&2522492&3620544&4577179\\
\end{tabular}
\end{table}
\normalsize

\small
\begin{table}
\caption{Number $T(n,9)$ and $T_t(n,9)$ of tilings of $9\times n$ boards with 
$1\times 4$ tiles.
}
\begin{tabular}{rr|rrrrrrrrrrrrrrr}
$n$ & & 0 & 1 & 2 & 3 & 4 & 5 & 6 & 7 &8 &9 \\
\hline
4&10&8&0&0&2&0&0&0&0&0&0\\
8&229&13&30&53&40&52&36&2&0&3&0\\
12&5764&63&314&534&870&973&870&751&616&350&240\\
16&143765&261&1020&4081&9208&14592&18820&21126&19374&17658&13658\\
20&3556413&668&4296&19034&61038&138550&239868&337187&420514&461175&449884\\
\end{tabular}
\label{tab.run9_4_1}
\end{table}
\normalsize

\clearpage
\subsection{Results (incongruent)}
Counts of tilings with $1\times 4$ tiles where only one representative of the
roto-reflected copies of each tiling is counted are shown
in Tables \ref{tab.run3i_4_1}--\ref{tab.run9i_4_1}.

\small
\begin{table}
\caption{Number of incongruent tilings of $3\times n$ boards with 
$1\times 4$ tiles.
}
\begin{tabular}{rr|rrrrrrrrrrrrrrr}
$n$ & & 0 & 1 & 2 & 3 & 4 & 5 & 6 & 7 &8 &9 & 10 & 11 & 12 & 13 & 14 \\
\hline
4&1&1&0&0&0&0&0&0&0&0&0&0&0&0&0&0\\
8&1&0&0&1&0&0&0&0&0&0&0&0&0&0&0&0\\
12&1&0&0&0&0&1&0&0&0&0&0&0&0&0&0&0\\
16&1&0&0&0&0&0&0&1&0&0&0&0&0&0&0&0\\
20&1&0&0&0&0&0&0&0&0&1&0&0&0&0&0&0\\
24&1&0&0&0&0&0&0&0&0&0&0&1&0&0&0&0\\
28&1&0&0&0&0&0&0&0&0&0&0&0&0&1&0&0\\
32&1&0&0&0&0&0&0&0&0&0&0&0&0&0&0&1\\
36&1&0&0&0&0&0&0&0&0&0&0&0&0&0&0&0\\
\end{tabular}
\label{tab.run3i_4_1}
\end{table}
\normalsize

The row sums in Table \ref{tab.run4i_4_1} are found in
\cite[A192928]{EIS}, repeating those in Table \ref{tab.run3i_3_1}.
\begin{conj} (Table \ref{tab.run4i_4_1})
\begin{equation}
\bar T_0(z,4)=
-z^4+\frac{1}{14}[
-\frac{1-3z-3z^2+2z^3}{1-z^2+z^4}
+\frac{3-2z}{1-z+z^2}
+\frac{8+4z+5z^2+6z^3+3z^4+5z^5}{1-z^4-z^6}
+\frac{4+6z+5z^2}{1-z^2-z^3}
].
\end{equation}
\end{conj}
\begin{conj} (Table \ref{tab.run4i_4_1})
\begin{equation}
\bar T_3(z,4)=
z^8\frac{(z^4+z^3+z^2+z+1)(z-1)^2 }{ (z^4-z^2+1)(1-z^6-z^4)(z^2-z+1)^2(z^3+z^2-1)^2 }.
\end{equation}
\end{conj}

\small
\begin{table}
\caption{Number $\bar T(n,4)$ and $\bar T_t(n,4)$ of incongruent tilings of $4\times n$ boards with 
$1\times 4$ tiles.
}
\begin{tabular}{rr|rrrrrrrrrrrrrrr}
$n$ & & 0 & 1 & 2 & 3 & 4 & 5 & 6 & 7 &8 &9 & 10 & 11 & 12 & 13 & 14 \\
\hline
1&1&1&0&0&0&0&0&0&0&0&0&0&0&0&0&0\\
2&1&1&0&0&0&0&0&0&0&0&0&0&0&0&0&0\\
3&1&1&0&0&0&0&0&0&0&0&0&0&0&0&0&0\\
4&1&1&0&0&0&0&0&0&0&0&0&0&0&0&0&0\\
5&2&2&0&0&0&0&0&0&0&0&0&0&0&0&0&0\\
6&3&3&0&0&0&0&0&0&0&0&0&0&0&0&0&0\\
7&3&3&0&0&0&0&0&0&0&0&0&0&0&0&0&0\\
8&5&4&0&0&1&0&0&0&0&0&0&0&0&0&0&0\\
9&6&5&0&0&1&0&0&0&0&0&0&0&0&0&0&0\\
10&9&7&0&0&2&0&0&0&0&0&0&0&0&0&0&0\\
11&11&9&0&0&2&0&0&0&0&0&0&0&0&0&0&0\\
12&16&12&0&0&3&0&0&1&0&0&0&0&0&0&0&0\\
13&20&15&0&0&4&0&0&1&0&0&0&0&0&0&0&0\\
14&29&20&0&0&7&0&0&2&0&0&0&0&0&0&0&0\\
15&37&25&0&0&10&0&0&2&0&0&0&0&0&0&0&0\\
16&53&34&0&0&15&0&0&3&0&0&1&0&0&0&0&0\\
17&69&43&0&0&20&0&0&5&0&0&1&0&0&0&0&0\\
18&98&58&0&0&29&0&0&9&0&0&2&0&0&0&0&0\\
19&130&74&0&0&40&0&0&14&0&0&2&0&0&0&0&0\\
20&183&99&0&0&59&0&0&21&0&0&3&0&0&1&0&0\\
21&245&128&0&0&82&0&0&29&0&0&5&0&0&1&0&0\\
22&343&171&0&0&117&0&0&43&0&0&10&0&0&2&0&0\\
23&463&223&0&0&160&0&0&62&0&0&16&0&0&2&0&0\\
24&646&297&0&0&225&0&0&95&0&0&25&0&0&3&0&0\\
25&877&388&0&0&310&0&0&137&0&0&35&0&0&6&0&0\\
26&1220&516&0&0&435&0&0&201&0&0&54&0&0&12&0&0\\
27&1664&676&0&0&600&0&0&284&0&0&82&0&0&20&0&0\\
28&2310&899&0&0&834&0&0&411&0&0&131&0&0&31&0&0\\
29&3161&1181&0&0&1145&0&0&587&0&0&197&0&0&44&0&0\\
30&4381&1569&0&0&1583&0&0&849&0&0&296&0&0&69&0&0\\
31&6009&2065&0&0&2172&0&0&1210&0&0&430&0&0&108&0&0\\
32&8319&2741&0&0&2995&0&0&1728&0&0&639&0&0&177&0&0\\
33&11430&3614&0&0&4103&0&0&2442&0&0&941&0&0&272&0&0\\
34&15811&4795&0&0&5636&0&0&3467&0&0&1400&0&0&415&0&0\\
35&21751&6330&0&0&7704&0&0&4892&0&0&2050&0&0&615&0&0\\
36&30070&8395&0&0&10552&0&0&6925&0&0&2998&0&0&932&0&0\\
37&41405&11091&0&0&14405&0&0&9741&0&0&4342&0&0&1407&0&0\\
38&57216&14705&0&0&19689&0&0&13716&0&0&6308&0&0&2138&0&0\\
39&78836&19440&0&0&26840&0&0&19222&0&0&9122&0&0&3198&0&0\\
40&108906&25770&0&0&36609&0&0&26962&0&0&13205&0&0&4761&0&0\\
41&150130&34085&0&0&49831&0&0&37700&0&0&19003&0&0&7028&0&0\\
42&207346&45176&0&0&67850&0&0&52725&0&0&27320&0&0&10401&0&0\\
43&285932&59775&0&0&92240&0&0&73536&0&0&39110&0&0&15337&0&0\\
44&394838&79215&0&0&125410&0&0&102527&0&0&55972&0&0&22618&0&0\\
45&544623&104845&0&0&170289&0&0&142631&0&0&79880&0&0&33151&0&0\\
46&751969&138930&0&0&231211&0&0&198350&0&0&113908&0&0&48496&0&0\\
47&1037425&183921&0&0&313600&0&0&275350&0&0&161990&0&0&70648&0&0\\
48&1432263&243696&0&0&425280&0&0&382056&0&0&230091&0&0&102847&0&0\\
49&1976229&322666&0&0&576240&0&0&529298&0&0&326101&0&0&149318&0&0\\
\end{tabular}
\label{tab.run4i_4_1}
\end{table}
\normalsize

\small
\begin{table}
\caption{Number $\bar T(n,5)$ and $\bar T_t(n,5)$ of incongruent tilings of $5\times n$ boards with 
$1\times 4$ tiles.
}
\begin{tabular}{rr|rrrrrrrrrrrrrrr}
$n$ & & 0 & 1 & 2 & 3 & 4 & 5 & 6 & 7 &8 &9 \\
\hline
4&2&2&0&0&0&0&0&0&0&0&0\\
8&6&3&2&0&0&1&0&0&0&0&0\\
12&23&6&8&5&2&0&1&0&0&1&0\\
16&103&16&30&26&13&6&4&6&0&0&1\\
20&478&37&92&124&86&52&29&28&14&3&6\\
24&2314&98&297&478&470&326&208&167&124&53&34\\
28&11285&245&898&1756&2128&1872&1327&1008&782&497&271\\
32&55529&646&2717&6099&8832&9225&7641&5973&4726&3489&2208\\
36&273652&1664&8012&20469&34121&41562&39643&33464&27314&21569&15483\\
40&1351040&4381&23473&66774&125799&174068&189092&175382&150883&125132&97403\\
44&6672248&11430&67867&213114&446184&691163&841663&861532&792499&689894&570063\\

\end{tabular}
\end{table}
\normalsize

\small
\begin{table}
\caption{Number $\bar T(n,6)$ and $\bar T_t(n,6)$ of incongruent tilings of $6\times n$ boards with 
$1\times 4$ tiles.
}
\begin{tabular}{rr|rrrrrrrrrrrrrrr}
$n$ & & 0 & 1 & 2 & 3 & 4 & 5 & 6 & 7 &8 &9 & 10 & 11 & 12 \\
\hline
2&0&0&0&0&0&0&0&0&0&0&0&0&0&0\\
4&3&3&0&0&0&0&0&0&0&0&0&0&0&0\\
6&0&0&0&0&0&0&0&0&0&0&0&0&0&0\\
8&11&3&3&4&0&0&1&0&0&0&0&0&0&0\\
10&0&0&0&0&0&0&0&0&0&0&0&0&0&0\\
12&53&3&7&16&12&10&2&0&2&0&0&1&0&0\\
14&0&0&0&0&0&0&0&0&0&0&0&0&0&0\\
16&278&2&12&49&52&66&45&24&10&5&10&0&0&2\\
18&0&0&0&0&0&0&0&0&0&0&0&0&0&0\\
20&1578&1&14&97&175&304&292&285&162&95&60&41&25&7\\
22&0&0&0&0&0&0&0&0&0&0&0&0&0&0\\
24&9262&1&16&154&469&1006&1355&1696&1465&1156&690&464&314&205\\
26&0&0&0&0&0&0&0&0&0&0&0&0&0&0\\
28&55530&1&17&209&973&2750&4885&7478&8692&8666&6905&5117&3426&2356\\
30&0&0&0&0&0&0&0&0&0&0&0&0&0&0\\
32&336015&1&19&268&1704&6352&14625&26894&38865&47763&48399&42806&33568&24361\\
34&0&0&0&0&0&0&0&0&0&0&0&0&0&0\\
36&2044505&1&20&329&2623&12610&37312&82511&143431&211334&259060&275895&256334&214503\\

\end{tabular}
\end{table}
\normalsize

\small
\begin{table}
\caption{Number $\bar T(n,7)$ and $\bar T_t(n,7)$ of incongruent tilings of $7\times n$ boards with 
$1\times 4$ tiles.
}
\begin{tabular}{rr|rrrrrrrrrrrrrrr}
$n$ & & 0 & 1 & 2 & 3 & 4 & 5 & 6 & 7 &8 &9 & 10 & 11 & 12 & 13 \\
\hline
4&3&3&0&0&0&0&0&0&0&0&0&0&0&0&0\\
8&13&1&3&4&4&0&0&1&0&0&0&0&0&0&0\\
12&77&1&1&11&14&19&14&12&2&0&2&0&0&1&0\\
16&513&1&1&11&26&69&100&92&91&59&31&11&6&12&0\\
20&3599&1&1&13&31&118&259&432&575&634&539&423&237&138&79\\
24&25763&1&1&15&39&173&449&1053&2015&2928&3671&4021&3638&2789&1902\\
28&185823&1&1&17&45&237&674&1922&4531&8760&14406&20186&24325&26147&24062\\
\end{tabular}
\end{table}
\normalsize

\small
\begin{table}
\caption{Number $\bar T(n,8)$ and $\bar T_t(n,8)$ of incongruent tilings of $8\times n$ boards with 
$1\times 4$ tiles.
}
\begin{tabular}{rr|rrrrrrrrrrrrrrr}
$n$ & & 0 & 1 & 2 & 3 & 4 & 5 & 6 & 7 &8 &9 & 10 \\
\hline
1&1&1&0&0&0&0&0&0&0&0&0&0\\
2&1&0&1&0&0&0&0&0&0&0&0&0\\
3&1&0&0&1&0&0&0&0&0&0&0&0\\
4&5&4&0&0&1&0&0&0&0&0&0&0\\
5&6&3&2&0&0&1&0&0&0&0&0&0\\
6&11&3&3&4&0&0&1&0&0&0&0&0\\
7&13&1&3&4&4&0&0&1&0&0&0&0\\
8&19&0&2&6&4&6&0&0&1&0&0&0\\
9&70&5&8&17&11&15&11&1&0&2&0&0\\
10&138&6&18&25&28&20&21&16&1&0&\\
11&230&5&21&41&45&45&23&28&17&2&0&3\\
12&496&2&19&49&79&96&86&61&52&40&4&0\\
13&1014&1&14&84&137&197&197&159&80&85&39&8\\
14&2106&12&34&138&243&335&394&355&255&134&115&58\\
15&3993&16&61&199&405&568&668&680&567&385&178&152\\
16&8152&12&87&280&600&962&1146&1273&1191&1036&670&388\\
17&16803&4&73&330&824&1604&2235&2569&2589&2258&1810&1135\\
18&34946&3&57&477&1282&2693&4141&5169&5327&4889&3960&3034\\
19&70170&28&123&696&2063&4335&7150&9560&10646&10076&8548&6487\\
20&142629&47&230&1015&3210&6945&11769&16609&19765&20411&18272&14915\\
21&292627&29&324&1374&4472&10536&18822&28383&36316&40825&39922&34669\\
22&605021&7&283&1667&5846&15827&31088&49804&67646&79514&83092&76865\\
23&1236177&9&202&2223&8289&23894&50784&86534&123502&152467&164960&161521\\
24&2524938&77&418&3114&12625&36308&81557&145961&219059&282657&320298&323951\\
25&5173937&133&796&4450&18609&55108&126525&237764&373162&508477&608074&650492\\
\end{tabular}
\end{table}
\normalsize

\small
\begin{table}
\caption{Number $\bar T(n,9)$ and $\bar T_t(n,9)$ of incongruent tilings of $9\times n$ boards with 
$1\times 4$ tiles.
}
\begin{tabular}{rr|rrrrrrrrrrrrrrr}
$n$ & & 0 & 1 & 2 & 3 & 4 & 5 & 6 & 7 &8 &9 & 10 \\
\hline
4&6&5&0&0&1&0&0&0&0&0&0&0\\
8&70&5&8&17&11&15&11&1&0&2&0&0\\
12&1505&19&80&144&219&254&223&196&157&94&62&45\\
16&36239&73&258&1048&2311&3697&4711&5325&4855&4454&3428&2431\\
20&890546&179&1074&4803&15272&34758&60008&84484&105184&115485&112524&100040\\
\end{tabular}
\label{tab.run9i_4_1}
\end{table}
\normalsize

\clearpage

\section{Tiling with $2\times 3$ tiles}\label{sec.23}
Tilings with $2\times 3$ tiles are
represented by Tables \ref{tab.run3_3_2}--\ref{tab.run9i_3_2}.
\subsection{Results (full count)}

The pattern in Table \ref{tab.run3_3_2} is obvious:
all hexominoes have the same orientation along the strip.

\small
\begin{table}
\caption{Number $T(n,3)$ and $T_t(n,3)$ of tilings of $3\times n$ boards with 
$2\times 3$ tiles.
}
\begin{tabular}{rr|rrrrrrrrrrrrrrr}
$n$ & & 0 & 1 & 2 & 3 & 4 & 5 & 6 & 7 &8 &9 & 10 & 11 & 12 & 13 & 14 \\
\hline
2&1&1&0&0&0&0&0&0&0&0&0&0&0&0&0&0\\
4&1&1&0&0&0&0&0&0&0&0&0&0&0&0&0&0\\
6&1&1&0&0&0&0&0&0&0&0&0&0&0&0&0&0\\
8&1&1&0&0&0&0&0&0&0&0&0&0&0&0&0&0\\
10&1&1&0&0&0&0&0&0&0&0&0&0&0&0&0&0\\
12&1&1&0&0&0&0&0&0&0&0&0&0&0&0&0&0\\
14&1&1&0&0&0&0&0&0&0&0&0&0&0&0&0&0\\
16&1&1&0&0&0&0&0&0&0&0&0&0&0&0&0&0\\
\end{tabular}
\label{tab.run3_3_2}
\end{table}
\normalsize

The pattern in Table \ref{tab.run4_3_2} is also obvious:
all hexominoes have the same orientation along the strip,
because there is only one composition of $m=4$ into parts
of the side lengths $t_n$ and $t_m$.

\small
\begin{table}
\caption{Number $T(n,4)$ and $T_t(n,4)$ of tilings of $4\times n$ boards with 
$2\times 3$ tiles.
}
\begin{tabular}{rr|rrrrrrrrrrrrrrr}
$n$ & & 0 & 1 & 2 & 3 & 4 & 5 & 6 & 7 &8 &9 & 10 & 11 & 12 & 13 & 14 \\
\hline
3&1&1&0&0&0&0&0&0&0&0&0&0&0&0&0&0\\
6&1&0&1&0&0&0&0&0&0&0&0&0&0&0&0&0\\
9&1&0&0&1&0&0&0&0&0&0&0&0&0&0&0&0\\
12&1&0&0&0&1&0&0&0&0&0&0&0&0&0&0&0\\
15&1&0&0&0&0&1&0&0&0&0&0&0&0&0&0&0\\
18&1&0&0&0&0&0&1&0&0&0&0&0&0&0&0&0\\
21&1&0&0&0&0&0&0&1&0&0&0&0&0&0&0&0\\
24&1&0&0&0&0&0&0&0&1&0&0&0&0&0&0&0\\
27&1&0&0&0&0&0&0&0&0&1&0&0&0&0&0&0\\
30&1&0&0&0&0&0&0&0&0&0&1&0&0&0&0&0\\
33&1&0&0&0&0&0&0&0&0&0&0&1&0&0&0&0\\
36&1&0&0&0&0&0&0&0&0&0&0&0&1&0&0&0\\
\end{tabular}
\label{tab.run4_3_2}
\end{table}
\normalsize

In Table \ref{tab.run5_3_2} the row sums are powers of 2
because the width $m=5$ allows only tiling with super-tiles
of shape $6\times 5$ containing 5 tiles of shape $2\times 3$
with two orientations. Each of these blocks may be inserted into
the floor with two different orientations. The Tatami tilings are
$T_0(z,5)=2/(1-z^6)$ because the two orientations must be filled
in alteratingly to avoid 4-crossings. Distributing the
places along the ``long'' edge of the floor where the
Tatami violations may occur leads to $T_t(6n,5)=2\binom{n-1}{t}$
\cite[A028326]{EIS}.

\small
\begin{table}
\caption{Number $T(n,5)$ and $T_t(n,5)$ of tilings of $5\times n$ boards with 
$2\times 3$ tiles.
}
\begin{tabular}{rr|rrrrrrrrrrrrrrr}
$n$ & & 0 & 1 & 2 & 3 & 4 & 5 & 6 & 7 &8 &9 & 10 & 11 & 12 & 13 & 14 \\
\hline
6&2&2&0&0&0&0&0&0&0&0&0&0&0&0&0&0\\
12&4&2&2&0&0&0&0&0&0&0&0&0&0&0&0&0\\
18&8&2&4&2&0&0&0&0&0&0&0&0&0&0&0&0\\
24&16&2&6&6&2&0&0&0&0&0&0&0&0&0&0&0\\
30&32&2&8&12&8&2&0&0&0&0&0&0&0&0&0&0\\
36&64&2&10&20&20&10&2&0&0&0&0&0&0&0&0&0\\
42&128&2&12&30&40&30&12&2&0&0&0&0&0&0&0&0\\
48&256&2&14&42&70&70&42&14&2&0&0&0&0&0&0&0\\
54&512&2&16&56&112&140&112&56&16&2&0&0&0&0&0&0
\end{tabular}
\label{tab.run5_3_2}
\end{table}
\normalsize

In Table \ref{tab.run6_3_2}, the floor width $m=6$
allows only consecutive placements of super-tiles of $6\times 2$
(with 2 tiles)
or super-tiles $6\times 3$ (with 3 tiles) along the long floor axis.
Therefore the row sums $T(n,6)$ are basically the
Padovan sequence \cite[A000931]{EIS}.
\begin{thm} (Table \ref{tab.run6_3_2})
\begin{equation}
T(z,6)=
\frac{1}{ 1-z^2-z^3 }.
\end{equation}
\end{thm}
Column $t=0$ is periodic with period length 5 because the only
way to avoid points where 4 tiles meet is to alternate the
two forms of super-tiles:
\begin{thm} (Table \ref{tab.run6_3_2})
\begin{equation}
T_0(z,6)=
\frac{(1+z)(1+z^2)(1-z+z^2)}{1-z^5}.
\end{equation}
\end{thm}
The other columns are basically given by considering the
number of compositions of $n$ into parts of 2 or 3, where
pairs of adjacent 2 increase $t$ by 1 and pairs of adjacent 3 increase
$t$ by 2.
\begin{conj} (Table \ref{tab.run6_3_2})
\begin{equation}
T_1(z,6)
=
 z^4\frac{(1+z)^2(z^2-z+1)^2}{(1-z^5)^2}.
\end{equation}
\end{conj}
\begin{conj} (Table \ref{tab.run6_3_2})
\begin{equation}
T_2(z,6)
=
z^6\frac{2+2z^2+2z^3+z^4-z^5+z^6-2z^7-z^9}{ (1-z^5)^3}.
\end{equation}
\end{conj}
The appearance of one more factor $1-z^5$ in the denominator
of the previous three equations each time $t$ increases by 1 is
a repeated convolution with $1/(1-z^5)=1+z^5+z^{10}+z^{15}+\cdots$
and understood by repeated attachment of blocks two super-tiles with
total dimension $6\times 5$ to the shorter floors of length $n-5$,
$n-10$ and so on with adaptation of the up-down aligment of the new
blocks to define one more point where 4 tiles meet at the interface.

\small
\begin{table}
\caption{Number $T(n,6)$ and $T_t(n,6)$ of tilings of $6\times n$ boards with 
$2\times 3$ tiles.
}
\begin{tabular}{rr|rrrrrrrrrrrrrrr}
$n$ & & 0 & 1 & 2 & 3 & 4 & 5 & 6 & 7 &8 &9 & 10 & 11 & 12 & 13 \\
\hline
1&0&0&0&0&0&0&0&0&0&0&0&0&0&0&0\\
2&1&1&0&0&0&0&0&0&0&0&0&0&0&0&0\\
3&1&1&0&0&0&0&0&0&0&0&0&0&0&0&0\\
4&1&0&1&0&0&0&0&0&0&0&0&0&0&0&0\\
5&2&2&0&0&0&0&0&0&0&0&0&0&0&0&0\\
6&2&0&0&2&0&0&0&0&0&0&0&0&0&0&0\\
7&3&1&2&0&0&0&0&0&0&0&0&0&0&0&0\\
8&4&1&0&2&1&0&0&0&0&0&0&0&0&0&0\\
9&5&0&2&2&0&1&0&0&0&0&0&0&0&0&0\\
10&7&2&1&1&2&1&0&0&0&0&0&0&0&0&0\\
11&9&0&0&5&2&2&0&0&0&0&0&0&0&0&0\\
12&12&1&4&1&2&2&1&1&0&0&0&0&0&0&0\\
13&16&1&0&4&6&3&2&0&0&0&0&0&0&0&0\\
14&21&0&3&6&1&6&2&3&0&0&0&0&0&0&0\\
15&28&2&2&2&8&7&4&2&0&1&0&0&0&0&0\\
16&37&0&0&9&8&7&7&3&3&0&0&0&0&0&0\\
17&49&1&6&3&6&12&8&9&2&2&0&0&0&0&0\\
18&65&1&0&6&16&13&13&8&4&3&0&1&0&0&0\\
19&86&0&4&12&4&18&16&17&10&3&2&0&0&0&0\\
20&114&2&3&3&18&24&21&19&9&10&3&2&0&0&0\\
21&151&0&0&14&20&17&32&24&26&11&4&2&0&1&0\\
22&200&1&8&6&12&36&33&42&25&20&11&4&2&0&0\\
23&265&1&0&8&32&36&42&48&36&35&12&11&2&2&0\\
24&351&0&5&20&10&40&60&63&66&36&31&12&5&2&0\\
25&465&2&4&4&32&59&66&79&66&67&44&25&12&3&2\\
26&616&0&0&20&40&35&90&100&114&93&52&42&13&13&2\\
27&816&1&10&10&20&80&96&136&128&116&102&59&38&13&4\\
28&1081&1&0&10&55&80&101&165&166&186&123&94&53&28&14\\
29&1432&0&6&30&20&75&160&184&252&204&198&141&80&51&14\\
30&1897&2&5&5&50&120&162&228&268&303&279&198&141&71&43\\
31&2513&0&0&27&70&65&200&300&364&420&322&312&184&135&64\\
32&3329&1&12&15&30&150&226&353&436&472&521&414&318&193&96\\
33&4410&1&0&12&86&155&206&430&548&677&646&553&458&287&196\\
34&5842&0&7&42&35&126&350&455&730&780&850&830&612&483&250\\
35&7739&2&6&6&72&216&342&534&804&1042&1164&1048&912&664&450\\
36&10252&0&0&35&112&112&385&735&960&1380&1350&1475&1240&978&693\\
37&13581&1&14&21&42&252&462&791&1170&1505&1940&1922&1744&1414&958\\
38&17991&1&0&14&126&273&378&945&1470&1965&2400&2397&2420&1929&1533\\
39&23833&0&8&56&56&196&672&994&1780&2335&2860&3420&3104&2852&2074\\
40&31572&2&7&7&98&357&651&1092&1995&2961&3765&4149&4207&3842&3075\\
41&41824&0&0&44&168&182&672&1568&2220&3720&4420&5340&5680&5129&4490\\
42&55405&1&16&28&56&392&854&1596&2682&4032&5890&6852&7248&7101&5982\\
43&73396&1&0&16&176&448&644&1848&3416&4899&7170&8308&9560&9358&8472\\
44&97229&0&9&72&84&288&1176&1974&3852&5913&8118&11250&12041&12580&11435\\
45&128801&2&8&8&128&554&1148&2030&4340&7322&10254&13356&15444&16489&15532\\
46&170625&0&0&54&240&282&1092&3024&4656&8748&12181&16137&20380&21058&21290\\
\end{tabular}
\label{tab.run6_3_2}
\end{table}
\normalsize

The number of tilings in Table \ref{tab.run7_3_2} are $T(n,3)=3^{n/6}$, powers of 3,
because there are only 3 compositions of 7 in parts of the
two side lengths $t_n$ and $t_m$ of the tile. The tiling consists
of a linear package of super-tiles of shape $7\times 6$ that
contain 7 tiles of shape $2\times 3$ each.

\small
\begin{table}
\caption{Number $T(n,7)$ and $T_t(n,7)$  of tilings of $7\times n$ boards with 
$2\times 3$ tiles.
}
\begin{tabular}{rr|rrrrrrrrrrrrrrr}
$n$ & & 0 & 1 & 2 & 3 & 4 & 5 & 6 & 7 &8 &9 & 10 & 11 & 12 & 13 & 14 \\
\hline
6&3&1&2&0&0&0&0&0&0&0&0&0&0&0&0&0\\
12&9&0&0&7&0&2&0&0&0&0&0&0&0&0&0&0\\
18&27&0&0&0&8&9&8&0&2&0&0&0&0&0&0&0\\
24&81&0&0&0&0&8&24&25&12&10&0&2&0&0&0&0\\
30&243&0&0&0&0&0&8&40&66&61&38&16&12&0&2&0\\
36&729&0&0&0&0&0&0&8&56&138&184&153&100&54&20&14\\
42&2187&0&0&0&0&0&0&0&8&72&242&436&496&409&262&148
\end{tabular}
\label{tab.run7_3_2}
\end{table}
\normalsize

For the row sums in Table \ref{tab.run8_3_2} we find:
\begin{conj} (Table \ref{tab.run8_3_2})
\begin{equation}
T(z,8) =
\frac{1-2z^3-z^6}{1-3z^3+z^9-z^{12}}.
\end{equation}
\end{conj}
An appropriate fit to column $T_0$ is
\begin{conj} (Table \ref{tab.run8_3_2})
\begin{equation}
T_0(z,8) =
\frac{1+2z^3+3z^6+4z^9+5z^{12}+4z^{15}+3z^{18}} {1+z^3+z^6-z^{12}-z^{15}-z^{18}}.
\end{equation}
\end{conj}

\small
\begin{table}
\caption{Number $T(n,8)$ and $T_t(n,8)$ of tilings of $8\times n$ boards with 
$2\times 3$ tiles.
}
\begin{tabular}{rr|rrrrrrrrrrrrrrr}
$n$ & & 0 & 1 & 2 & 3 & 4 & 5 & 6 & 7 &8 &9 & 10 & 11 & 12 \\
\hline
3&1&1&0&0&0&0&0&0&0&0&0&0&0&0\\
6&4&1&0&2&1&0&0&0&0&0&0&0&0&0\\
9&11&2&4&0&4&0&0&1&0&0&0&0&0&0\\
12&33&3&0&9&10&4&0&6&0&0&1&0&0&0\\
15&96&1&4&14&22&12&24&6&8&0&4&0&0&1\\
18&281&2&4&31&21&57&38&51&36&18&6&12&0&4\\
21&821&1&16&14&68&106&118&145&124&66&84&36&20&6\\
24&2400&3&4&40&119&171&262&415&314&335&261&186&122&75\\
27&7015&2&8&58&170&263&632&772&926&1074&864&782&520&365\\
30&20505&1&12&95&162&526&1098&1558&2407&2551&2744&2671&2058&1576\\
33&59936&1&24&62&284&907&1646&3184&4960&6134&7694&7719&6772&6401\\
36&175193&4&8&99&447&1174&2698&6016&9313&13852&18880&19957&22022&19724\\
39&512089&1&16&132&584&1492&4608&9862&17156&29873&40004&51252&58818&59298\\
42&1496836&1&16&200&545&2322&7156&15128&31920&55890&84242&119677&144697&162460\\
45&4375251&2&32&136&796&3478&9568&24175&55092&100161&170268&255318&334008&420009\\
\end{tabular}
\label{tab.run8_3_2}
\end{table}
\normalsize

\small
\begin{table}
\caption{Number $T(n,9)$ and $T_t(n,9)$ of tilings of $9\times n$ boards with 
$2\times 3$ tiles.
}
\begin{tabular}{rr|rrrrrrrrrrrrrrr}
$n$ & & 0 & 1 & 2 & 3 & 4 & 5 & 6 & 7 &8 &9 & 10 & 11 & 12 \\
\hline
2&1&1&0&0&0&0&0&0&0&0&0&0&0&0\\
4&1&0&0&1&0&0&0&0&0&0&0&0&0&0\\
6&5&0&2&2&0&1&0&0&0&0&0&0&0&0\\
8&11&2&4&0&4&0&0&1&0&0&0&0&0&0\\
10&19&0&8&0&4&2&4&0&0&1&0&0&0&0\\
12&45&0&0&8&8&8&10&4&6&0&0&1&0&0\\
14&105&0&0&12&12&24&30&8&4&10&4&0&0&1\\
16&219&2&0&20&26&38&22&44&36&8&6&12&4&0\\
18&475&0&8&16&24&68&64&52&70&70&62&12&8&16\\
20&1061&0&0&16&36&72&120&178&172&144&114&98&68&8\\
22&2313&0&0&10&46&108&238&330&300&362&318&184&162&128\\
24&5027&2&0&12&92&158&288&512&684&662&582&644&524&312\\
26&11035&0&10&16&68&240&474&636&1002&1452&1542&1446&1188&1108\\
28&24173&0&0&26&68&206&584&1154&1828&2422&3032&3442&3122&2406\\
30&52793&0&0&12&66&260&752&1724&2990&4434&5782&6240&6448&6554\\
32&115499&2&0&12&106&376&924&2388&4590&7000&9774&12316&13694&13558\\
34&252827&0&12&20&80&492&1350&2902&5996&11076&16704&22040&27320&30388\\
\end{tabular}
\label{tab.run9_3_2}
\end{table}
\normalsize

\clearpage
\subsection{Results (incongruent)}
Counts of tilings with $2\times 3$ tiles where only one representative of the
roto-reflected copies of each tiling is counted are shown
in Tables \ref{tab.run3i_3_2}--\ref{tab.run9i_3_2}.

\small
\begin{table}
\caption{Number of $\bar T(n,3)$ and $\bar T_t(n,3)$ incongruent tilings of $3\times n$ boards with 
$2\times 3$ tiles.
This is the same as Table \ref{tab.run3_3_2} because there is
only one tiling for each $n$.
}
\begin{tabular}{rr|rrrrrrrrrrrrrrr}
$n$ & & 0 & 1 & 2 & 3 & 4 & 5 & 6 & 7 &8 &9 & 10 & 11 & 12 & 13 & 14 \\
\hline
2&1&1&0&0&0&0&0&0&0&0&0&0&0&0&0&0\\
4&1&1&0&0&0&0&0&0&0&0&0&0&0&0&0&0\\
6&1&1&0&0&0&0&0&0&0&0&0&0&0&0&0&0\\
8&1&1&0&0&0&0&0&0&0&0&0&0&0&0&0&0\\
10&1&1&0&0&0&0&0&0&0&0&0&0&0&0&0&0\\
12&1&1&0&0&0&0&0&0&0&0&0&0&0&0&0&0\\
\end{tabular}
\label{tab.run3i_3_2}
\end{table}
\normalsize

\small
\begin{table}
\caption{Number $\bar T(n,4)$ and $\bar T_t(n,4)$ of incongruent tilings of $4\times n$ boards with 
$2\times 3$ tiles.
This is the same as Table \ref{tab.run4_3_2} because there is
only one tiling for each $n$.
}
\begin{tabular}{rr|rrrrrrrrrrrrrrr}
$n$ & & 0 & 1 & 2 & 3 & 4 & 5 & 6 & 7 &8 &9 & 10 & 11 & 12 & 13 & 14 \\
\hline
3&1&1&0&0&0&0&0&0&0&0&0&0&0&0&0&0\\
6&1&0&1&0&0&0&0&0&0&0&0&0&0&0&0&0\\
9&1&0&0&1&0&0&0&0&0&0&0&0&0&0&0&0\\
12&1&0&0&0&1&0&0&0&0&0&0&0&0&0&0&0\\
15&1&0&0&0&0&1&0&0&0&0&0&0&0&0&0&0\\
18&1&0&0&0&0&0&1&0&0&0&0&0&0&0&0&0\\
\end{tabular}
\end{table}
\normalsize

The row sums of Table \ref{tab.run5i_3_2} are fitted by the 
following generating function:
\begin{conj} (Table \ref{tab.run5i_3_2})
\begin{equation}
\bar T(z,5)=
\frac{1}{4}[
1+\frac{1}{1-2z^6}
]
+\frac{1+z^6}{2(1-2z^{12})}.
\end{equation}
\end{conj}
Followup columns appear to have similarly simple forms:
\begin{conj} (Table \ref{tab.run5i_3_2})
\begin{equation}
\bar T_0(z,5)= \frac{1} {1-z^6}.
\end{equation}
\end{conj}
\begin{conj} (Table \ref{tab.run5i_3_2})
\begin{equation}
\bar T_1(z,5)= z^{12}\frac{1} {(1+z^6)(1-z^6)^2}.
\end{equation}
\end{conj}
\begin{conj} (Table \ref{tab.run5i_3_2})
\begin{equation}
\bar T_2(z,5)= z^{18}\frac{1} {(1+z^6)(1-z^6)^3}.
\end{equation}
\end{conj}
\begin{conj} (Table \ref{tab.run5i_3_2})
\begin{equation}
\bar T_3(z,5)= z^{24}\frac{1+z^{12}} {(1+z^6)^2(1-z^6)^4}.
\end{equation}
\end{conj}
\begin{conj} (Table \ref{tab.run5i_3_2})
\begin{equation}
\bar T_4(z,5)= z^{30}\frac{1+z^{12}} {(1+z^6)^2(1-z^6)^5}.
\end{equation}
\end{conj}

\small
\begin{table}
\caption{Number $\bar T(n,5)$ and $\bar T_t(n,5)$ of incongruent tilings of $5\times n$ boards with 
$2\times 3$ tiles.
}
\begin{tabular}{rr|rrrrrrrrrrrrrrr}
$n$ & & 0 & 1 & 2 & 3 & 4 & 5 & 6 & 7 &8 &9 & 10 & 11 & 12 & 13 & 14 \\
\hline
6&1&1&0&0&0&0&0&0&0&0&0&0&0&0&0&0\\
12&2&1&1&0&0&0&0&0&0&0&0&0&0&0&0&0\\
18&3&1&1&1&0&0&0&0&0&0&0&0&0&0&0&0\\
24&6&1&2&2&1&0&0&0&0&0&0&0&0&0&0&0\\
30&10&1&2&4&2&1&0&0&0&0&0&0&0&0&0&0\\
36&20&1&3&6&6&3&1&0&0&0&0&0&0&0&0&0\\
42&36&1&3&9&10&9&3&1&0&0&0&0&0&0&0&0\\
48&72&1&4&12&19&19&12&4&1&0&0&0&0&0&0&0\\
54&136&1&4&16&28&38&28&16&4&1&0&0&0&0&0&0\\
60&272&1&5&20&44&66&66&44&20&5&1&0&0&0&0&0\\
66&528&1&5&25&60&110&126&110&60&25&5&1&0&0&0&0\\
72&1056&1&6&30&85&170&236&236&170&85&30&6&1&0&0&0\\
78&2080&1&6&36&110&255&396&472&396&255&110&36&6&1&0&0\\
84&4160&1&7&42&146&365&651&868&868&651&365&146&42&7&1&0\\
90&8256&1&7&49&182&511&1001&1519&1716&1519&1001&511&182&49&7&1\\
96&16512&1&8&56&231&693&1512&2520&3235&3235&2520&1512&693&231&56&8\\
102&32896&1&8&64&280&924&2184&4032&5720&6470&5720&4032&2184&924&280&64\\
108&65792&1&9&72&344&1204&3108&6216&9752&12190&12190&9752&6216&3108&1204&344
\end{tabular}
\label{tab.run5i_3_2}
\end{table}
\normalsize

In Table \ref{tab.run6i_3_2} we conjecture
\begin{conj} (Table \ref{tab.run6i_3_2}, row sum)
\begin{equation}
\bar T(z,6)
=
1-z^6+\frac{1}{2}[
\frac{1}{1-z^2-z^3}
+\frac{1+z^2+z^3}{1-z^4-z^6}
].
\end{equation}
\end{conj}
The column $t=0$ is periodic with period length 5 for the same
reason as in Table \ref{tab.run6_3_2}, but entries of 2 are reduced to 1
because the two tilings counted in \ref{tab.run6_3_2} are congruent
to each other.
\begin{conj} (Table \ref{tab.run6i_3_2})
\begin{equation}
\bar T_1(z,6)
=
z^4 +z^7\frac{1+z^2+z^3+z^5+z^7-z^{12}}{ (1+z^5)(1-z^5)^2}.
\end{equation}
\end{conj}
\begin{conj} (Table \ref{tab.run6i_3_2})
\begin{equation}
\bar T_2(z,6)
=
z^6 +z^8\frac{1+z+z^2+3z^3+z^4+z^6-z^7-z^{10}-4z^{13}+2z^{18}}{ (1+z^5)(1-z^5)^3}
.
\end{equation}
\end{conj}

\small
\begin{table}
\caption{Number $\bar T(n,6)$ and $\bar T_t(n,6)$ of incongruent tilings of $6\times n$ boards with 
$2\times 3$ tiles.
}
\begin{tabular}{rr|rrrrrrrrrrrrrrr}
$n$ & & 0 & 1 & 2 & 3 & 4 & 5 & 6 & 7 &8 &9 & 10 & 11 & 12 & 13 & 14 \\
\hline
1&0&0&0&0&0&0&0&0&0&0&0&0&0&0&0&0\\
2&1&1&0&0&0&0&0&0&0&0&0&0&0&0&0&0\\
3&1&1&0&0&0&0&0&0&0&0&0&0&0&0&0&0\\
4&1&0&1&0&0&0&0&0&0&0&0&0&0&0&0&0\\
5&1&1&0&0&0&0&0&0&0&0&0&0&0&0&0&0\\
6&1&0&0&1&0&0&0&0&0&0&0&0&0&0&0&0\\
7&2&1&1&0&0&0&0&0&0&0&0&0&0&0&0&0\\
8&3&1&0&1&1&0&0&0&0&0&0&0&0&0&0&0\\
9&3&0&1&1&0&1&0&0&0&0&0&0&0&0&0&0\\
10&5&1&1&1&1&1&0&0&0&0&0&0&0&0&0&0\\
11&5&0&0&3&1&1&0&0&0&0&0&0&0&0&0&0\\
12&8&1&2&1&1&1&1&1&0&0&0&0&0&0&0&0\\
13&9&1&0&2&3&2&1&0&0&0&0&0&0&0&0&0\\
14&13&0&2&3&1&4&1&2&0&0&0&0&0&0&0&0\\
15&15&1&1&1&4&4&2&1&0&1&0&0&0&0&0&0\\
16&22&0&0&6&4&4&4&2&2&0&0&0&0&0&0&0\\
17&26&1&3&2&3&6&4&5&1&1&0&0&0&0&0&0\\
18&37&1&0&3&9&7&7&5&2&2&0&1&0&0&0&0\\
19&45&0&2&6&2&10&8&9&5&2&1&0&0&0&0&0\\
20&63&1&2&2&9&14&11&10&5&6&2&1&0&0&0&0\\
21&78&0&0&8&10&9&16&12&13&6&2&1&0&1&0&0\\
22&108&1&4&4&6&18&18&23&13&11&6&3&1&0&0&0\\
23&136&1&0&4&16&19&21&25&18&18&6&6&1&1&0&0\\
24&186&0&3&10&6&22&30&34&34&19&16&7&3&1&0&1\\
25&237&1&2&2&16&31&33&40&33&35&22&13&6&2&1&0\\
26&322&0&0&12&20&19&46&51&59&49&27&22&7&8&1&1\\
27&414&1&5&6&10&40&48&70&64&59&51&30&19&7&2&1\\
28&559&1&0&5&29&41&52&86&83&96&63&50&27&15&8&2\\
29&724&0&3&15&10&39&80&94&126&104&99&72&40&26&7&7\\
30&973&1&3&3&25&64&82&116&136&155&142&102&72&37&22&9\\
31&1267&0&0&15&35&34&100&150&182&213&161&157&92&70&32&16\\
32&1697&1&6&9&15&75&116&181&220&241&262&212&161&100&49&30\\
33&2219&1&0&6&43&79&103&218&274&341&323&280&229&145&98&43\\
34&2964&0&4&21&19&66&175&234&368&394&428&421&309&245&127&86\\
35&3888&1&3&3&36&111&171&269&402&524&582&528&456&335&225&134\\
36&5183&0&0&20&56&59&194&369&485&699&679&744&623&498&349&196\\
37&6815&1&7&12&21&126&231&400&585&757&970&964&872&712&479&338\\
38&9071&1&0&7&65&138&192&480&735&992&1206&1209&1214&973&772&498\\
39&11949&0&4&28&28&100&336&502&890&1173&1430&1716&1552&1431&1037&749\\
40&15886&1&4&4&49&185&327&551&1002&1488&1890&2089&2111&1931&1542&1162\\
41&20955&0&0&24&84&94&336&784&1110&1869&2210&2676&2840&2573&2245&1580\\
42&27835&1&8&16&28&196&432&806&1346&2030&2948&3442&3634&3569&2998&2439\\
43&36755&1&0&8&88&226&322&930&1708&2457&3585&4166&4780&4688&4236&3454\\
44&48790&0&5&36&44&148&588&1000&1932&2967&4068&5643&6031&6311&5730&4768\\
45&64476&1&4&4&64&282&574&1020&2170&3666&5127&6693&7722&8257&7766&6846\\
46&85545&0&0&30&120&146&548&1514&2338&4395&6101&8091&10199&10560&10660&9158\\
47&113115&1&9&20&36&288&732&1496&2748&4758&7734&10179&12345&13974&13830&12872\\
48&150021&1&0&9&121&350&524&1672&3570&5493&9220&12225&15541&17995&18352&17591\\
49&198460&0&5&45&60&205&960&1828&3816&6641&10146&15822&19323&22667&24235&23110\\
50&263136&1&5&5&81&419&956&1768&4292&8158&12348&18501&24015&28850&31413&31450
\end{tabular}
\label{tab.run6i_3_2}
\end{table}
\normalsize

The fit to row sums of Table \ref{tab.run7i_3_2} is:
\begin{conj} (Table \ref{tab.run7i_3_2})
\begin{equation}
\bar T(z,7)
=
\frac{1-2z^6-4z^{12}+6z^{18}}{ (1-z^6)(1-3z^6)(1-3z^{12}) }.
\end{equation}
\end{conj}

\small
\begin{table}
\caption{Number $\bar T(n,7)$ and $\bar T_t(n,7)$ of incongruent tilings of $7\times n$ boards with 
$2\times 3$ tiles.
}
\begin{tabular}{rr|rrrrrrrrrrrrrrr}
$n$ & & 0 & 1 & 2 & 3 & 4 & 5 & 6 & 7 &8 &9 & 10 & 11 & 12 & 13 & 14 \\
\hline
6&2&1&1&0&0&0&0&0&0&0&0&0&0&0&0&0\\
12&4&0&0&3&0&1&0&0&0&0&0&0&0&0&0&0\\
18&10&0&0&0&3&4&2&0&1&0&0&0&0&0&0&0\\
24&25&0&0&0&0&3&6&9&3&3&0&1&0&0&0&0\\
30&70&0&0&0&0&0&3&12&18&17&11&5&3&0&1&0\\
36&196&0&0&0&0&0&0&3&14&39&46&43&25&16&5&4\\
42&574&0&0&0&0&0&0&0&3&20&64&114&128&107&68&38\\
48&1681&0&0&0&0&0&0&0&0&3&22&101&220&338&339&282\\
54&5002&0&0&0&0&0&0&0&0&0&3&28&142&404&741&961\\
60&14884&0&0&0&0&0&0&0&0&0&0&3&30&195&650&1485\\
66&44530&0&0&0&0&0&0&0&0&0&0&0&3&36&252&1014\\
72&133225&0&0&0&0&0&0&0&0&0&0&0&0&3&38&321
\end{tabular}
\label{tab.run7i_3_2}
\end{table}
\normalsize

\small
\begin{table}
\caption{Number $\bar T(n,8)$ and $\bar T_t(n,8)$ of incongruent tilings of $8\times n$ boards with 
$2\times 3$ tiles.
}
\begin{tabular}{rr|rrrrrrrrrrrrrrr}
$n$ & & 0 & 1 & 2 & 3 & 4 & 5 & 6 & 7 &8 &9 & 10 & 11 & 12 & 13 \\
\hline
3&1&1&0&0&0&0&0&0&0&0&0&0&0&0&0\\
6&3&1&0&1&1&0&0&0&0&0&0&0&0&0&0\\
9&4&1&1&0&1&0&0&1&0&0&0&0&0&0&0\\
12&13&2&0&3&3&2&0&2&0&0&1&0&0&0&0\\
15&28&1&1&4&6&3&6&3&2&0&1&0&0&1&0\\
18&84&1&1&10&6&17&10&16&9&5&3&4&0&1&0\\
21&216&1&4&4&18&29&30&37&31&18&22&9&5&3&3\\
24&639&2&1&12&31&49&66&110&81&92&66&48&32&22&11\\
27&1784&1&2&16&43&67&160&199&232&272&218&199&131&92&64\\
30&5238&1&3&29&42&138&277&404&603&656&693&682&518&407&303\\
33&15068&1&6&17&72&231&413&805&1243&1545&1925&1933&1697&1618&1142\\
36&44118&2&2&28&113&306&679&1534&2334&3499&4731&5032&5519&4973&4350\\
39&128257&1&4&35&148&379&1155&2481&4292&7486&10010&12840&14708&14853&14288\\
42&375126&1&4&56&138&597&1793&3828&7992&14044&21079&30027&36198&40734&42919\\
45&1094470&1&8&36&200&880&2400&6063&13779&25073&42582&63892&83522&105067&114960\\
48&3199848&2&4&53&279&1112&3295&9728&21969&44151&81627&127447&186220&246533&292094\\
\end{tabular}
\end{table}
\normalsize

\small
\begin{table}
\caption{Number $\bar T(n,9)$ and $\bar T_t(n,9)$ of incongruent tilings of $9\times n$ boards with 
$2\times 3$ tiles.
}
\begin{tabular}{rr|rrrrrrrrrrrrrrr}
$n$ & & 0 & 1 & 2 & 3 & 4 & 5 & 6 & 7 &8 &9 & 10 & 11 & 12 & 13 & 14 \\
\hline
2&1&1&0&0&0&0&0&0&0&0&0&0&0&0&0&0\\
4&1&0&0&1&0&0&0&0&0&0&0&0&0&0&0&0\\
6&3&0&1&1&0&1&0&0&0&0&0&0&0&0&0&0\\
8&4&1&1&0&1&0&0&1&0&0&0&0&0&0&0&0\\
10&7&0&3&0&1&1&1&0&0&1&0&0&0&0&0&0\\
12&15&0&0&3&2&3&3&1&2&0&0&1&0&0&0&0\\
14&31&0&0&4&4&6&8&3&1&3&1&0&0&1&0&0\\
16&62&1&0&5&7&12&6&12&9&3&2&3&1&0&0&1\\
18&131&0&2&6&6&17&17&15&19&19&17&4&2&5&1&0\\
20&279&0&0&5&9&21&31&46&43&38&29&27&17&3&4&4\\
22&603&0&0&3&13&29&61&87&78&92&81&48&42&34&20&4\\
24&1289&1&0&3&23&42&72&133&172&172&146&163&133&82&57&50\\
26&2810&0&3&4&19&61&121&163&255&368&389&372&299&280&204&96\\
28&6115&0&0&8&17&55&146&293&460&615&760&874&784&613&481&413\\
30&13315&0&0&3&18&68&192&436&753&1117&1452&1572&1619&1652&1324&979\\
32&29026&1&0&3&27&98&231&610&1151&1765&2448&3095&3430&3411&3232&2981\\
34&63463&0&3&6&20&127&340&732&1507&2786&4191&5531&6851&7618&7535&6694
\end{tabular}
\label{tab.run9i_3_2}
\end{table}
\normalsize

\clearpage
\section{Classified by Slide Lines}\label{sec.sli}

Table \ref{tab.run3s}--\ref{tab.run6s} show the $T(n,m)$ and the rows $\hat T_s(n,m)$
for columns labeled by $s$.

\subsection{Results (full count)}

Row sums of Table \ref{tab.run3s} are those of
Table \ref{tab.run3}. The column with no slide lines
is
\begin{conj} (Table \ref{tab.run3s})
\begin{equation}
\hat T_0(z,3) = 1+2z^4\frac{1}{ (1-z^2)(1-4z^2+z^4)(1-3z^2+z^4) }.
\end{equation}
\end{conj}
The column with $s=1$ slide lines is apparently
\begin{conj} (Table \ref{tab.run3s})
\begin{equation}
\hat T_1(z,3) = 2z^2\frac{1}{ (1-z^2)(1-3z^2+z^4) },
\end{equation}
\end{conj}
twice the sequence \cite[A027941]{EIS}.

\begin{thm} (Table \ref{tab.run3s})
$\hat T_2(n,3)=1$ for even $n$ and $\hat T_2(z,3)=z^2/(1-z^2)$, because there is one tiling with 2 slide lines; this consists
of a stack of 3 copies of the $1\times n$ floor tilings.
\end{thm}
The previous equations for $\hat T_s(z,3)$ are compatible
with the sum rule (\ref{eq.sumrul2}) and with Equation (\ref{eq.run3}).

\small
\begin{table}
\caption{Number $T(n,3)$ and $\hat T_s(n,3)$ of tilings of $3\times n$ boards with dominos.
}
\begin{tabular}{rr|rrrrrrrrrrrrrrr}
$n$ & & 0 & 1 & 2 \\
\hline
2&3&0&2&1\\
4&11&2&8&1\\
6&41&16&24&1\\
8&153&86&66&1\\
10&571&394&176&1\\
12&2131&1666&464&1\\
14&7953&6734&1218&1\\
16&29681&26488&3192&1\\
18&110771&102410&8360&1\\
20&413403&391512&21890&1\\
22&1542841&1485528&57312&1\\
24&5757961&5607912&150048&1\\
26&21489003&21096168&392834&1\\
\end{tabular}
\label{tab.run3s}
\end{table}
\normalsize

\begin{conj} (Table \ref{tab.run4s})
\begin{equation}
\hat T_1(z,4) = 
\frac{2(1-z^2)}{1-4z^2+z^4}
+\frac{3-2z}{5(1-3z+z^2)}
+\frac{17-2z}{5(1-z^2)}
-\frac{3}{1-z-z^2}
-\frac{3}{1+z-z^2}
.
\end{equation}
\end{conj}
\begin{conj} (Table \ref{tab.run4s})
\begin{equation}
\hat T_2(z,4) = 3z^2\frac{1}{(1-z^2)(1-3z^2+z^4)}.
\end{equation}
\end{conj}
Row sums of Table \ref{tab.run4s} are those in
Table \ref{tab.run4}. The column with $s=0$ slide lines
is deduced from the columns and the sum rule:
\begin{conj} (Table \ref{tab.run4s})
\begin{multline}
\hat T_0(z,4)=
1-\frac{2-2z^2}{1-4z^2+z^4}
-\frac{3-2z}{5(1-3z+z^2)}
-\frac{7-2z}{5(1-z^2)}
+\frac{3}{2(1-z-z^2)}
\\
+\frac{3}{2(1+z-z^2)}
+\frac{1-z^2}{1-z-5z^2-z^3+z^4}
.
\end{multline}
\end{conj}

\small
\begin{table}
\caption{Number $T(n,4)$ and $\hat T_s(n,4)$ of tilings of $4\times n$ boards with dominos.
}
\begin{tabular}{rr|rrrrrrrrrrrrrrr}
$n$ & & 0 & 1 & 2 & 3 \\
\hline
1&1&0&1&0&0\\
2&5&0&1&3&1\\
3&11&2&9&0&0\\
4&36&3&20&12&1\\
5&95&31&64&0&0\\
6&281&68&176&36&1\\
7&781&340&441&0&0\\
8&2245&884&1261&99&1\\
9&6336&3311&3025&0&0\\
10&18061&9264&8532&264&1\\
11&51205&30469&20736&0&0\\
12&145601&87748&57156&696&1\\
13&413351&271222&142129&0&0\\
14&1174500&788323&384349&1827&1\\
15&3335651&2361482&974169&0&0\\
16&9475901&6870920&2600192&4788&1\\
17&26915305&20238249&6677056&0&0\\
18&76455961&58766200&17677220&12540&1\\
19&217172736&171407511&45765225&0&0\\
20&616891945&496283060&120576049&32835&1\\
\end{tabular}
\label{tab.run4s}
\end{table}
\normalsize

\small
\begin{table}
\caption{Number $T(n,5)$ and $\hat T_s(n,5)$ of tilings of $5\times n$ boards with dominos.
}
\begin{tabular}{rr|rrrrrrrrrrrrrrr}
$n$ & & 0 & 1 & 2 & 3 & 4 & \\
\hline
2&8&0&0&3&4&1\\
4&95&2&22&54&16&1\\
6&1183&134&520&480&48&1\\
8&14824&3722&7444&3525&132&1\\
10&185921&73282&87872&24414&352&1\\
12&2332097&1216178&948520&166470&928&1\\
\end{tabular}
\label{tab.run5s}
\end{table}
\normalsize

\small
\begin{table}
\caption{Number $T(n,6)$ and $\hat T_s(n,6)$ of tilings of $6\times n$ boards with dominos.
}
\begin{tabular}{rr|rrrrrrrrrrrrrrr}
$n$ & & 0 & 1 & 2 & 3 & 4 & 5 \\
\hline
1&1&0&0&1&0&0&0\\
2&13&0&0&1&6&5&1\\
3&41&2&12&27&0&0&0\\
4&281&3&32&121&104&20&1\\
5&1183&175&496&512&0&0&0\\
6&6728&499&2156&3084&928&60&1\\
7&31529&7988&14280&9261&0&0&0\\
8&167089&31244&73184&55617&6878&165&1\\
9&817991&287406&364210&166375&0&0&0\\
10&4213133&1315092&1932264&917296&48040&440&1
\end{tabular}
\label{tab.run6s}
\end{table}
\normalsize

\small
\begin{table}
\caption{Slide-line free $\hat T_0(n,m)$ of tilings of $m\times n$ boards with dominos.
}
\begin{tabular}{r|rrrrrrrrrrrrrrr}
$m$ & 1 & 2 & 3 & 4 & 5 & 6 & 7 &8 &9 & 10 & 11 & 12 & 13 & 14 \\
\hline
1 & 0 & 1 &0 & 1 &0 & 1 &0 & 1 &0 & 1 &0 & 1 &0 & 1 \\
2 & 1 & 1 &3 & 4 &8 & 12 &21 & 33 &55 & 88 &144 & 232 &377 & 609 \\
3 & 0 & 0 &0 & 2 &0 & 16 &0 & 86 &0 & 394 &0 & 1666 &0 & 6734 \\
4 & 0 & 0 &2 & 3 &31 & 68 &340 & 884 &3311 & 9264 & 30469 & 87748 & 271222 & 788323 \\
5 & 0 & 0 &0 & 2 &0 & 134 &0 & 3722 &0 & 73282 & 0 & 1216178 & 0 & 18339220 \\
\end{tabular}
\label{tab.runms}
\end{table}
\normalsize

The main use of the classification of tilings by their number of slide lines
is that the full number of tilings is computable from the tilings without slide
lines by summing over all compositions of $m$ into positive
parts and building the product of all possibilities of the subfloor stacks,
or summing over all partitions of $m$ into positive parts and inserting
the multinomial coefficients as multiplicities:
\begin{multline}
T(n,m)
=
\sum_{m= m_1+m_2+m_3+\cdots}
\hat T_0(n,m_1)
\hat T_0(n,m_2)
\cdots
\\
=
\sum_{\pi(m)=\{m_1^{\alpha_1};m_2^{\alpha_2};m_3^{\alpha_3};\cdots\}}
\frac{(\alpha_1+\alpha_2+\cdots)!}{\alpha_1!\alpha_2!\cdots}
\hat T_0(n,m_1)^{\alpha_1}
\hat T_0(n,m_2)^{\alpha_2}
\cdots
.
\end{multline}
This transformation of the sequence $\hat T_0(n,m)$ into $T(n,m)$
is known as the Invert Transform \cite{BernsteinLAA226} and Cameron's
A-transform \cite{CameronDM75}, see Deutsch's comment in \cite[A005178]{EIS}.

That same strategy delivers $\hat T_s(n,m)$ for $s>0$ by restriction
of the number of parts in these two equations whenever
$\hat T_0(n,m')$ are known for all $1<=m'<=m-s$.
In this respect, Table \ref{tab.runms}, extracted from Tables \ref{tab.run3s}--\ref{tab.run6s},
contains all other columns of these tables.

\begin{exa}
For fixed $n=12$ and the partitions $4=\{4^1;1^13^1;2^2;1^22^1;1^4\}$
the sum in column $n$ of Table \ref{tab.runms} is
\begin{multline}
T(12,4)
=
\hat T_0(12,4)
+
2 \hat T_0(12,1)T_0(12,3)
+
\hat T_0(12,2)^2
+
3 \hat T_0(12,1)^2 T_0(12,2)
+
\hat T_0(12,1)^4
\\
= 87748
+2\times 1 \times 1666
+232^2
+3\times 1\times 232
+1^4
=145601
\end{multline}
as seen in row $n=12$ in Table \ref{tab.run4s}.
\end{exa}

\begin{rem}
A further reduction of the information may be defined if the tabulation
is not only refined along the index $s$ of the slide lines but along
a pair of indices that test an orthogonal grid of crossed slide lines that
cut through the floor.
\end{rem}

\subsection{Results (icongruental)}

A fit to $\hat{\bar T}_1(n,3)$ in Table \ref{tab.run3si} is 
\begin{conj} (Table \ref{tab.run3si})
\begin{equation}
\hat {\bar T}_1(z,3) = z^2\frac{1-2z^2}{(1-z^2)(1-3z^2+z^4)(1-z^2-z^4)}.
\end{equation}
\end{conj}
The generating function for $\hat{\bar T}_0$ follows then from (\ref{eq.run3i})
with $\hat {\bar T}_2(z,3)=z^2/(1-z^2)$ and with the sum rule:
\begin{conj} (Table \ref{tab.run3si})
\begin{multline}
\hat {\bar T}_0(z,3) =
1+\frac{1-z^2}{4(1-4z^2+z^4)}
-\frac{1+z^2}{2(1-z^2-z^4)}
+\frac{1}{4(1-z^2)}
-\frac{1}{4(1-z-z^2)}
\\
-\frac{1}{4(1+z-z^2)}
+\frac{1+2z^2-z^6}{2(1-4z^4+z^8)}
.
\end{multline}
\end{conj}

\small
\begin{table}
\caption{Number $\bar T(n,3)$ and $\hat {\bar T}_s(n,3)$ of incongruent
tilings of $3\times n$ boards with dominos.
}
\begin{tabular}{rr|rrrrrrrrrrrrrrr}
$n$ & & 0 & 1 & 2 \\
\hline
2&2&0&1&1&\\
4&5&1&3&1&\\
6&14&5&8&1&\\
8&46&25&20&1&\\
10&156&105&50&1&\\
12&561&434&126&1&\\
14&2037&1715&321&1&\\
16&7525&6699&825&1&\\
18&27874&25739&2134&1&\\
20&103741&98196&5544&1&\\
22&386386&371941&14444&1\\
24&1440946&1403245&37700&1\\
26&5374772&5276258&98513&1\\
\end{tabular}
\label{tab.run3si}
\end{table}
\normalsize

\begin{figure}
\caption{The $\bar T(4,4)=9$ tilings recorded in Tables \ref{tab.run4i} and \ref{tab.run4si}.}
\includegraphics[width=0.32\textwidth]{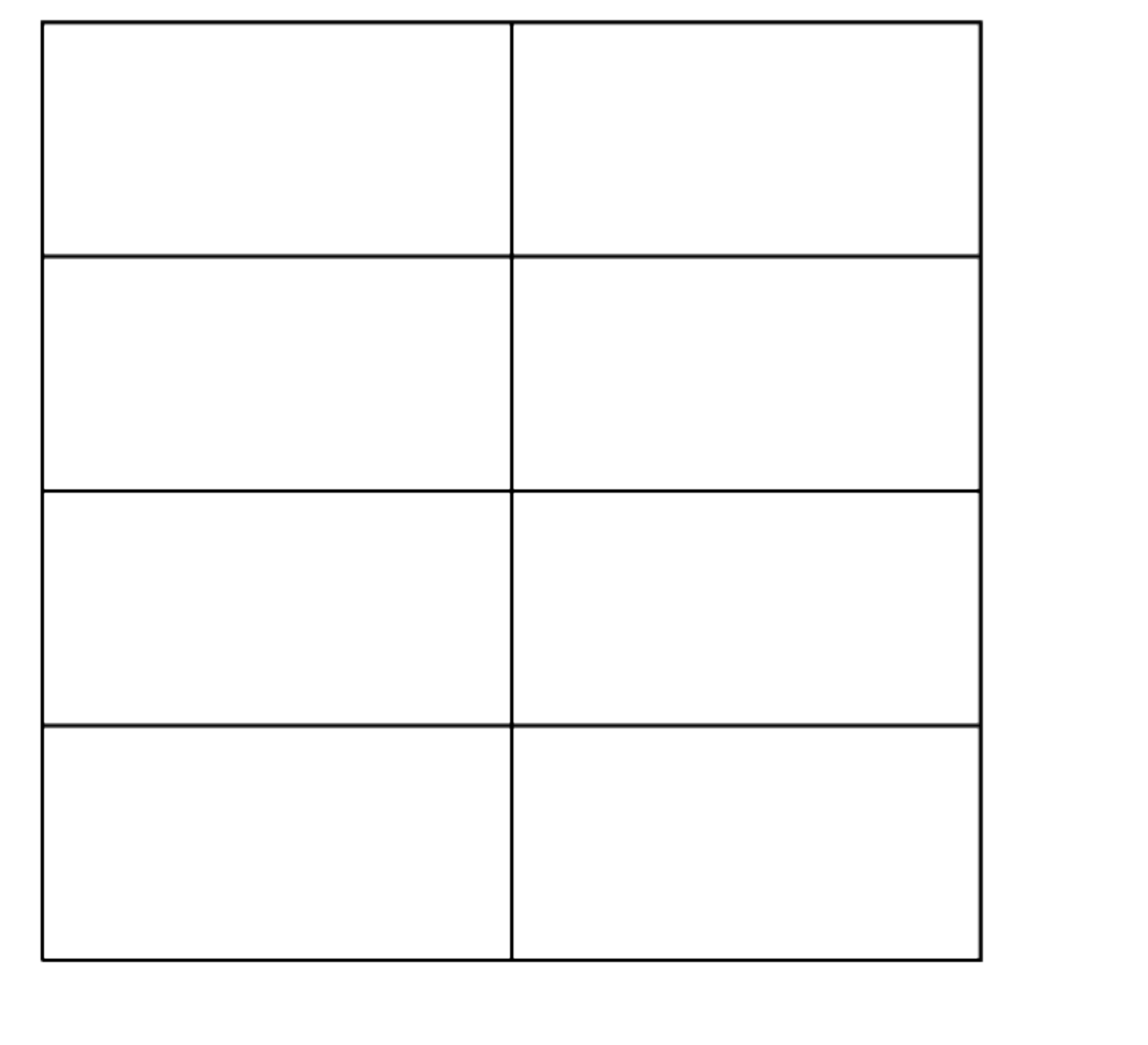}
\includegraphics[width=0.32\textwidth]{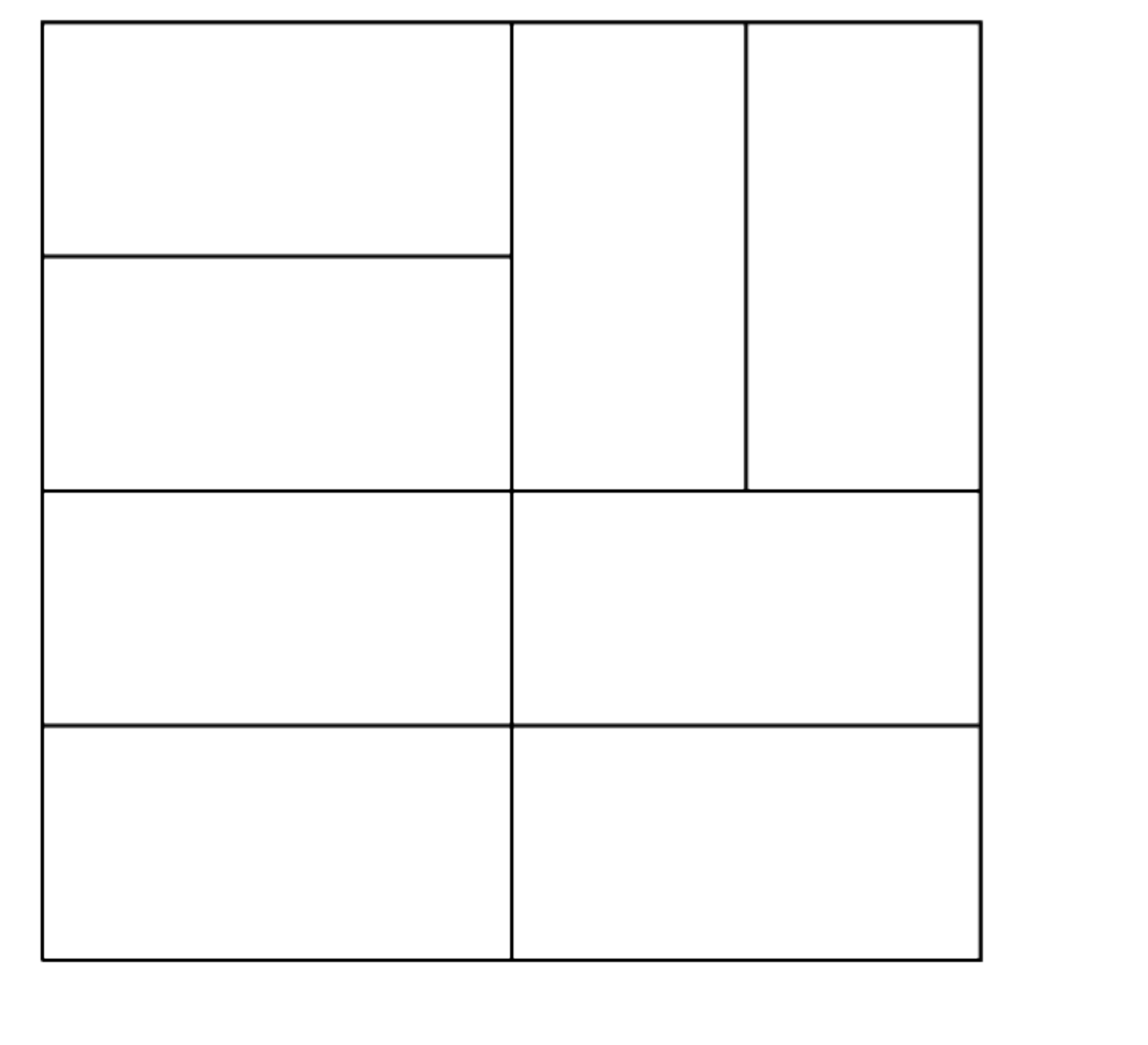}
\includegraphics[width=0.32\textwidth]{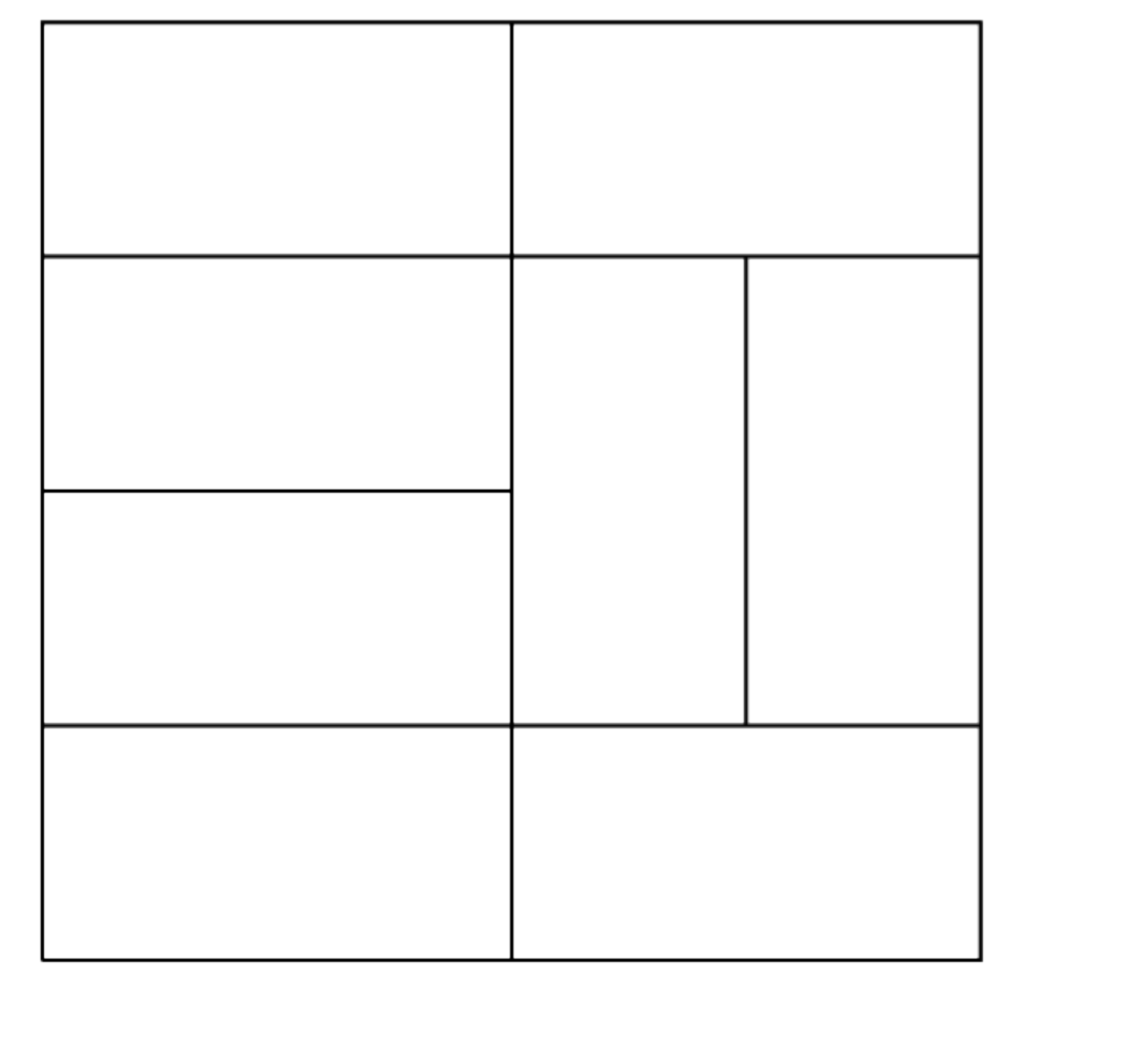}
\includegraphics[width=0.32\textwidth]{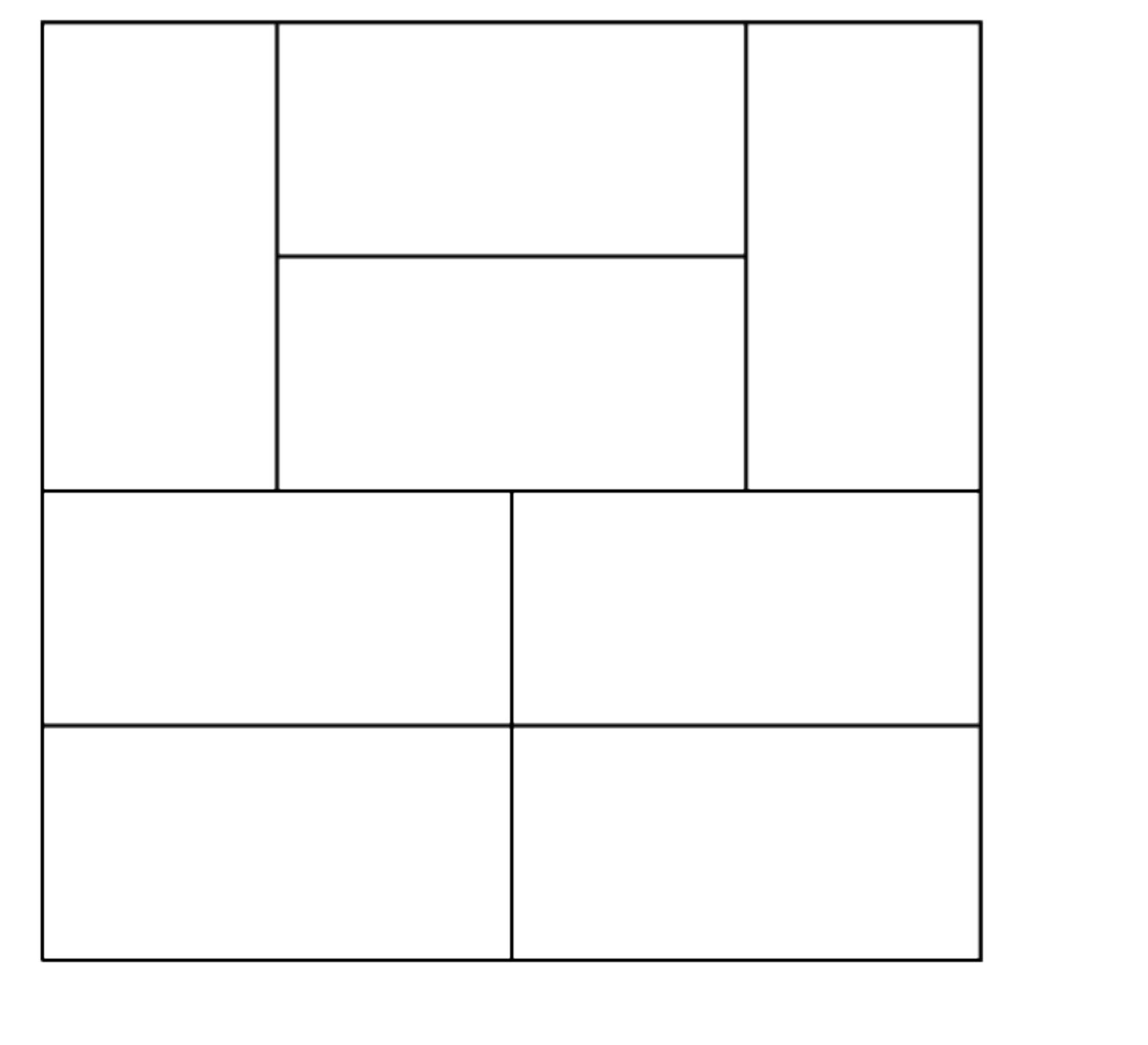}
\includegraphics[width=0.32\textwidth]{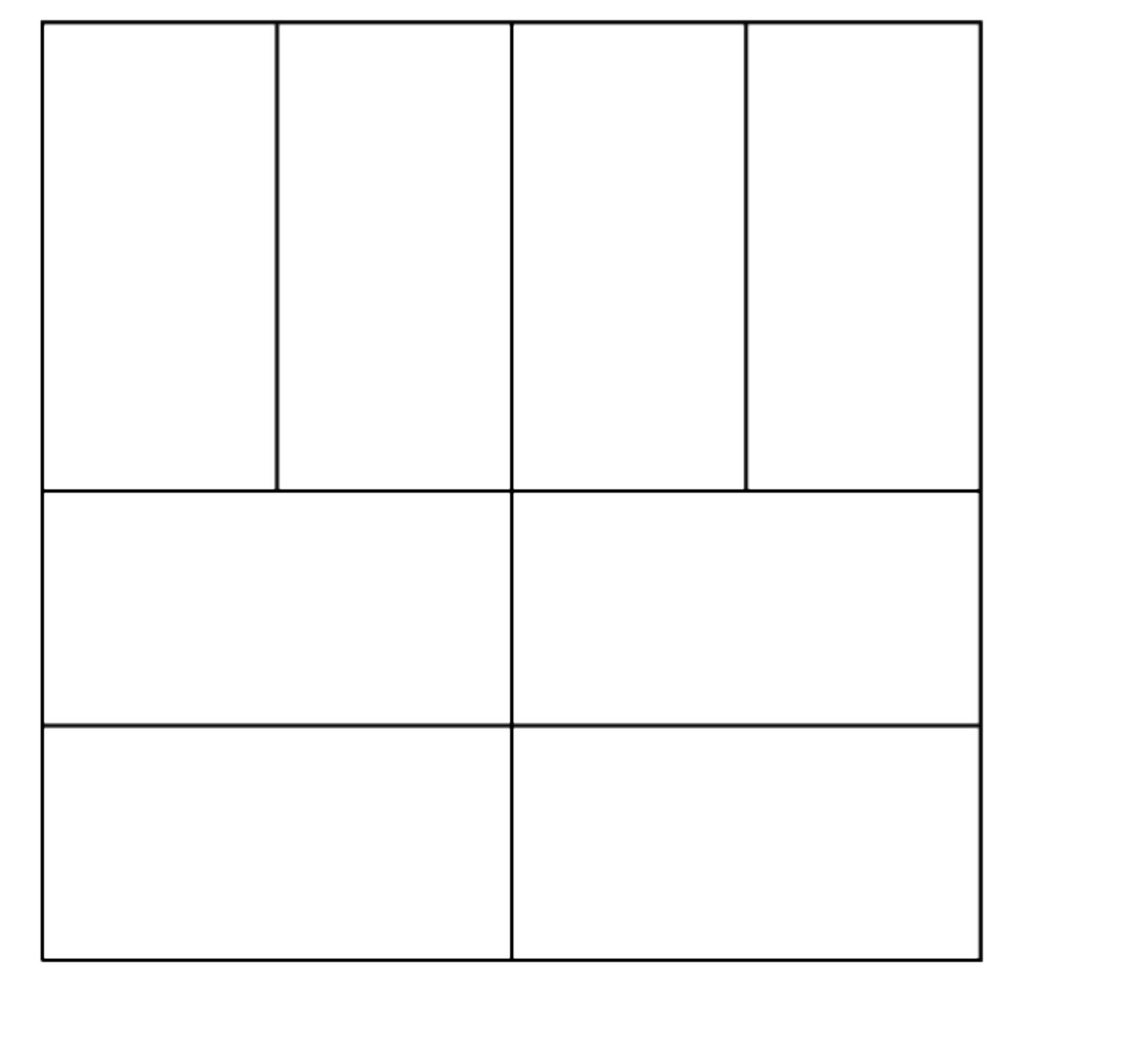}
\includegraphics[width=0.32\textwidth]{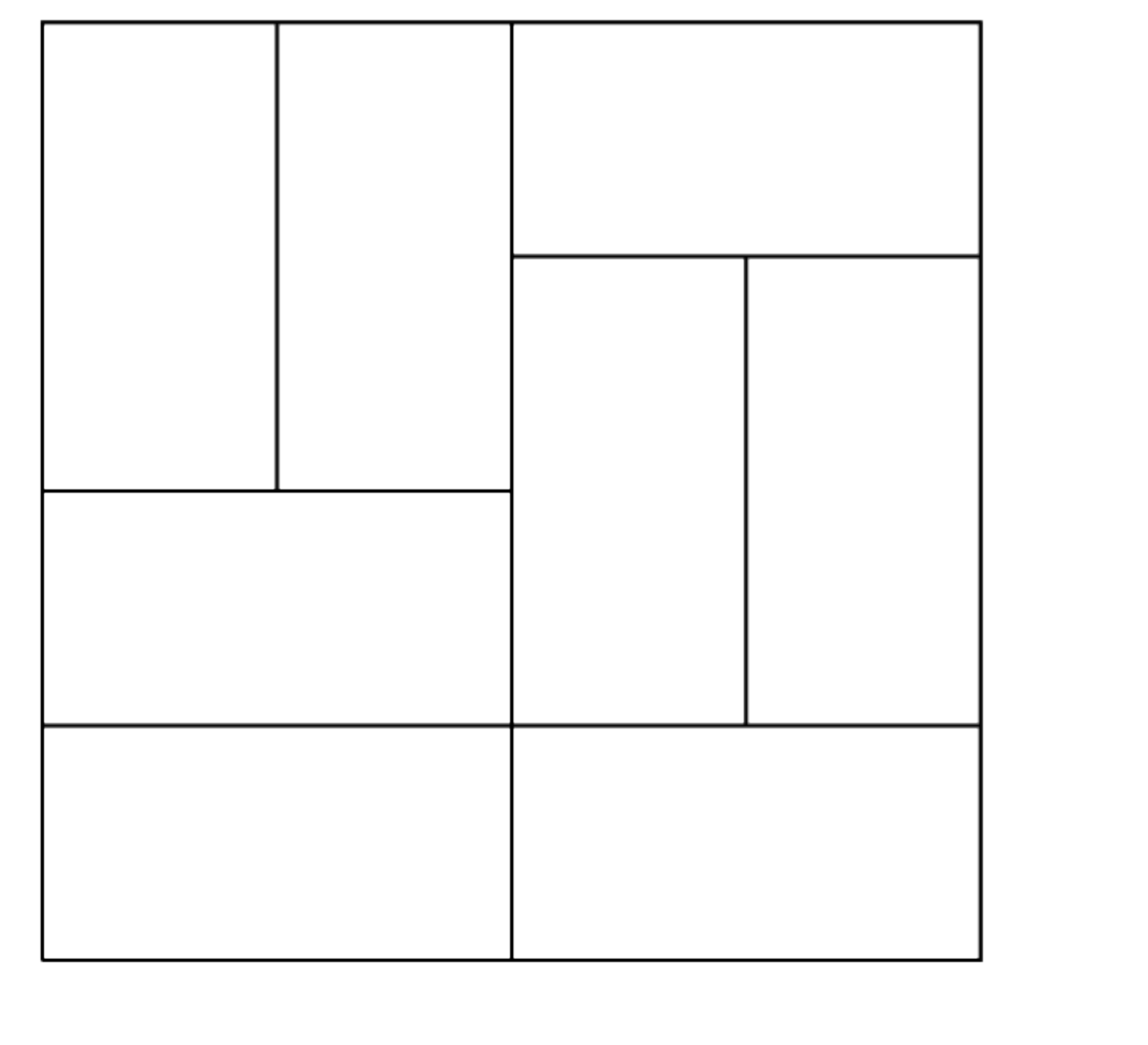}
\includegraphics[width=0.32\textwidth]{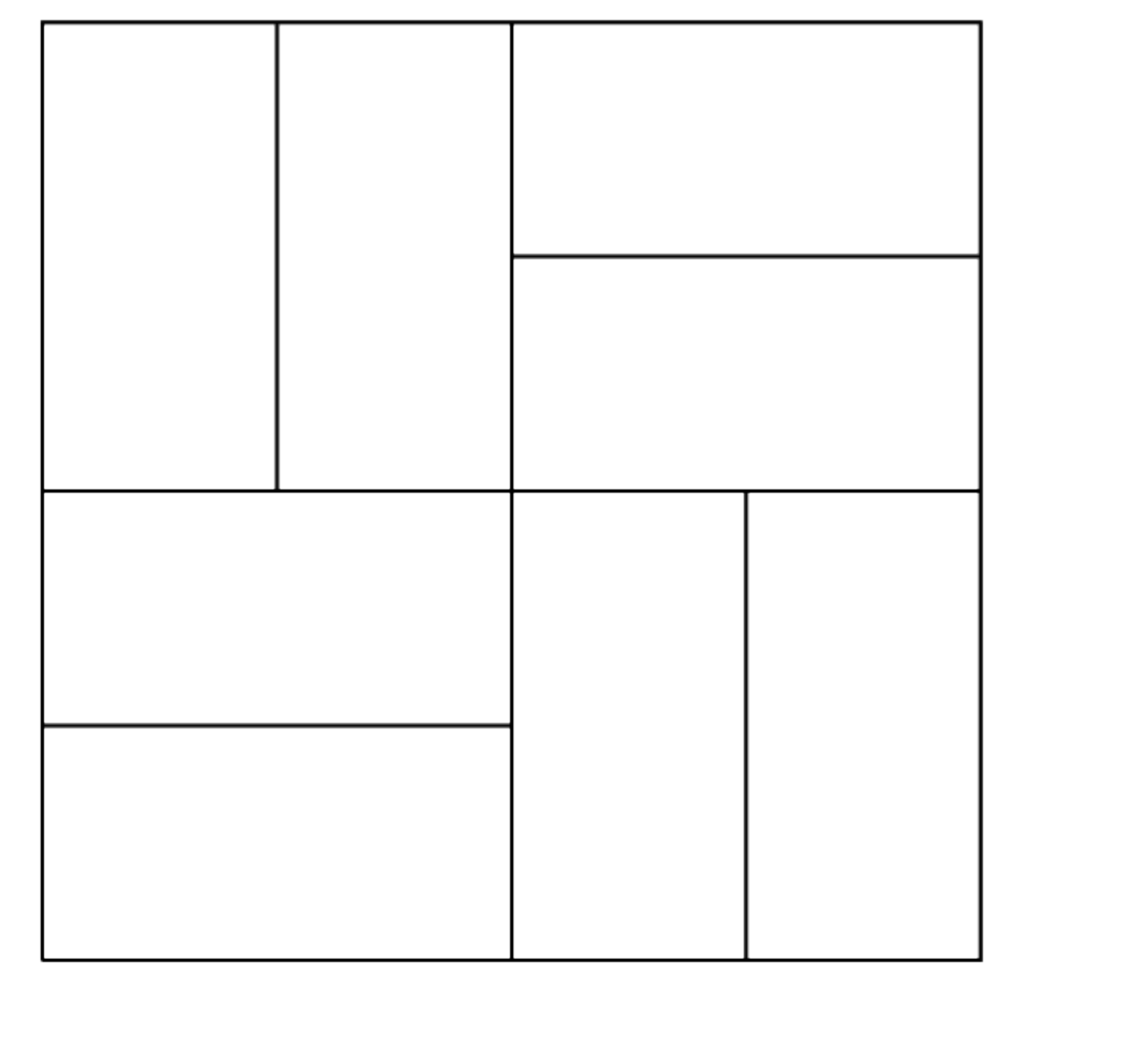}
\includegraphics[width=0.32\textwidth]{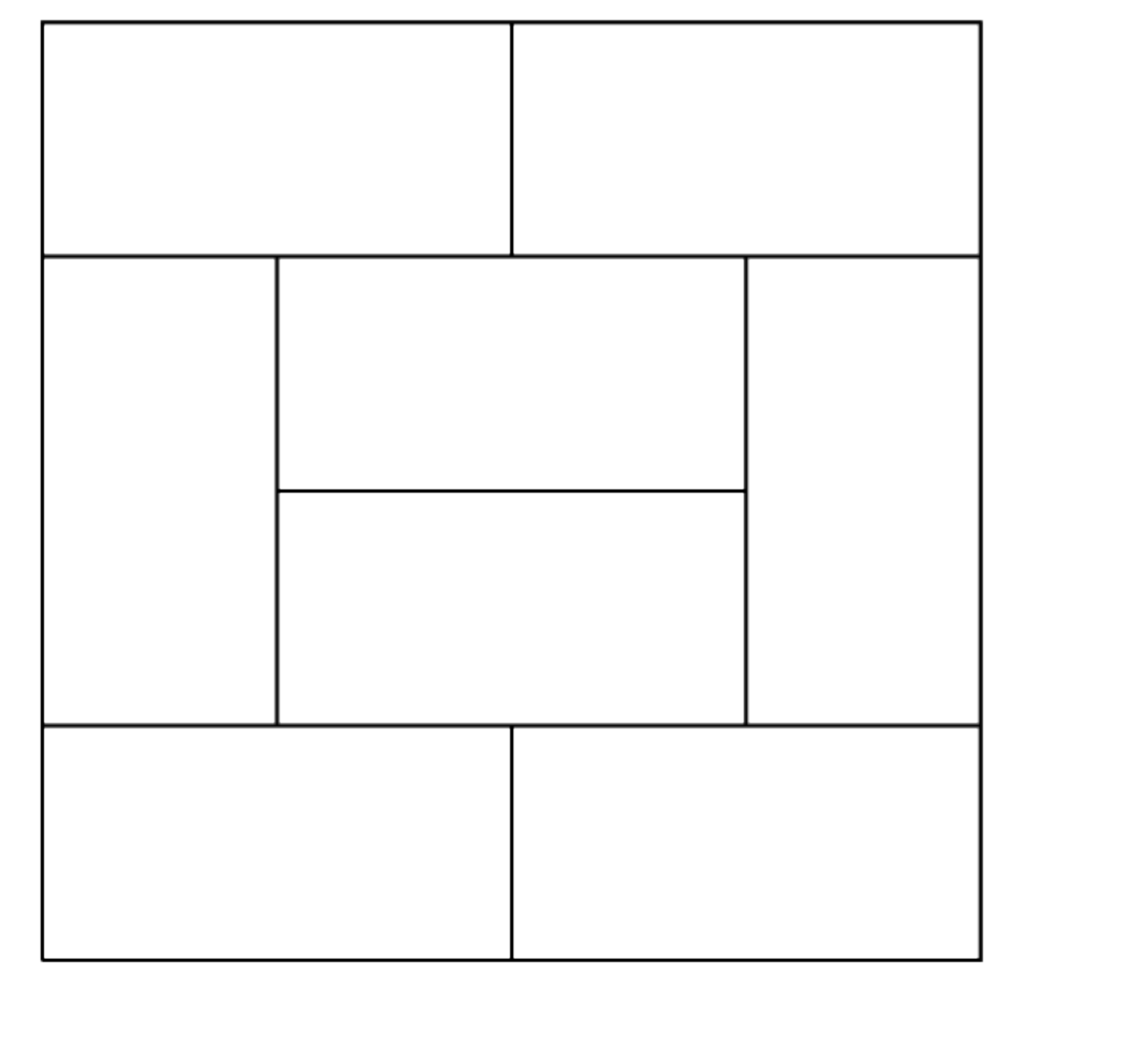}
\includegraphics[width=0.32\textwidth]{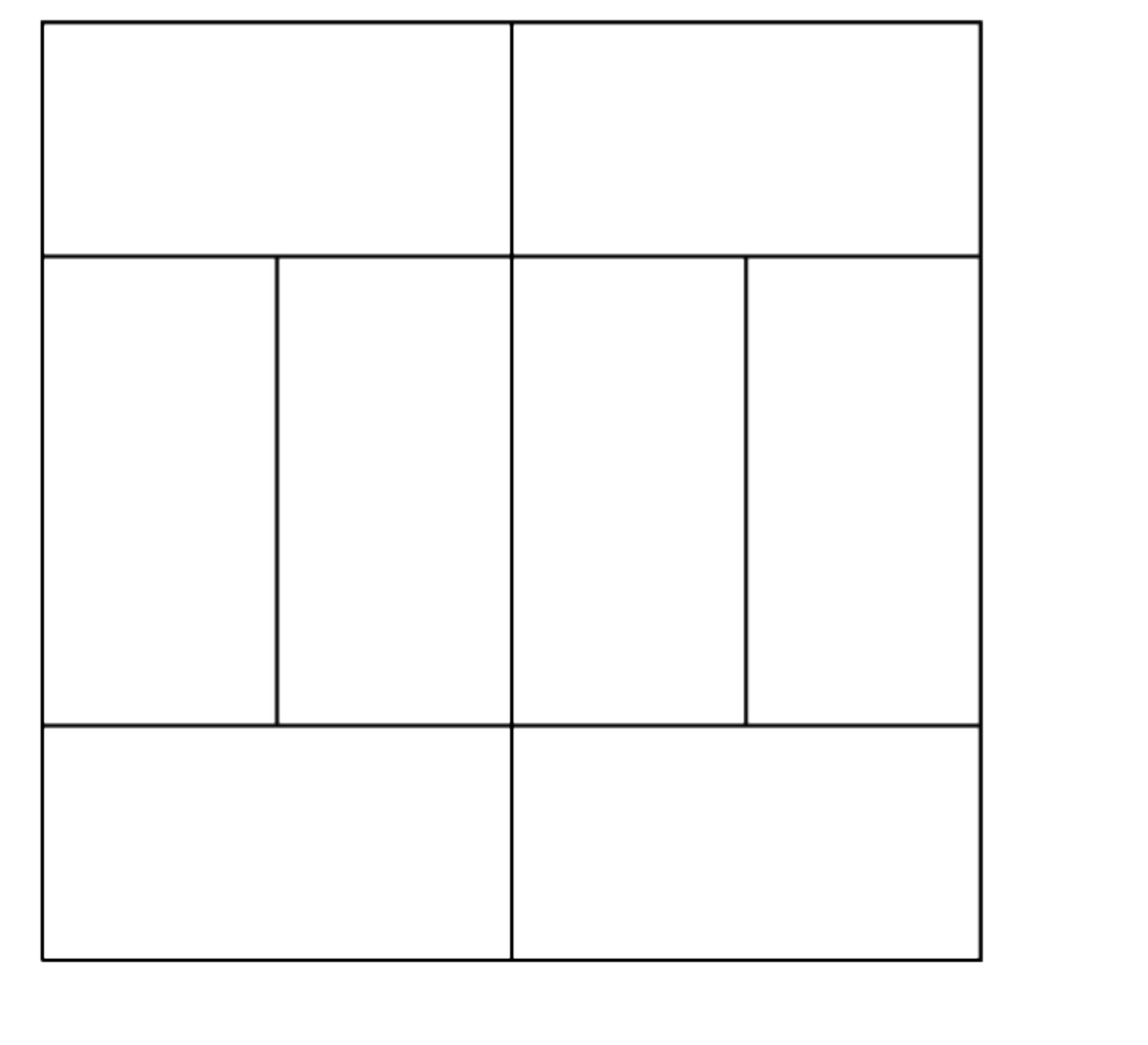}
\label{fig.run4si}
\end{figure}

\begin{rem}
Definition \ref{def.slid} is poor if $n=m$ and the matching condition is
fulfilled, because we have defined slide lines to
run parallel to the ``long'' edge. For square floors,
a tiling rotation by 90 degrees
may not conserve the slide line count as defined.
So the $\hat{\bar T}_s(n,n)$ depend on which of the tilings of a congruential
octet is chosen as the representative, and in that respect row $n=4$ in
Table \ref{tab.run4si} is undefined. Consider the upper left configuration
in Figure \ref{fig.run4si} for example: as shown it has 3 slide lines,
but if rotated by 90 degrees only 1.
\end{rem}

\small
\begin{table}
\caption{Number $\bar T(n,4)$ and $\hat {\bar T}_s(n,4)$ of incongruent
tilings of $4\times n$ boards with dominos.
}
\begin{tabular}{rr|rrrrrrrrrrrrrrr}
$n$ & & 0 & 1 & 2 & 3\\
\hline
1&1&0&1&0&0\\
2&4&0&1&2&1\\
3&5&1&4&0&0\\
4&9&0&2&6&1\\
5&33&12&21&0&0\\
6&98&26&55&16&1\\
7&230&107&123&0&0\\
8&658&271&346&40&1\\
9&1725&935&790&0&0\\
10&4876&2554&2221&100&1\\
11&13378&8106&5272&0&0\\
12&37794&23021&14520&252&1\\
13&105761&69998&35763&0&0\\
14&299221&201871&96707&642&1\\
15&844219&600073&244146&0&0\\
16&2392040&1738736&651653&1650&1\\
17&6773154&5102309&1670845&0&0\\
\end{tabular}
\label{tab.run4si}
\end{table}
\normalsize

The tables are also defined for tiles of other shapes;
one set for the full counts and one set for incongruent counts.
We eventually list two examples of them.

The heuristics for $1\times 3$ tiles splits the 
row sums (\ref{eq.run3_3_1}) as follows:
$\hat T_1(n,3)=0$,
$\hat T_2(z,3)=z^3/(1-z^3)$ and
\begin{conj} (Tilings of $3\times n$ floors
with $1\times 3$ tiles
without slide lines
\cite[A099560]{EIS})
\begin{equation}
\hat T_0(z,3) = 1+\frac{z}{(1-z^3)(1-z-z^3)}.
\end{equation}
\end{conj}

The heuristics for $1\times 3$ tiles splits the 
row sums (\ref{eq.run4_3_1}) as follows:
$\hat T_2(n,4)=0$,
$\hat T_3(z,4)=z^3/(1-z^3)$ and
\begin{conj} (Tilings of $4\times n$ floors
with $1\times 3$ tiles
with one slide line)
\begin{equation}
\hat T_1(z,4) = 2z^3\frac{1}{(1-z^3)(1-4z^3+3z^6-z^9)}.
\end{equation}
\end{conj}
\begin{conj} (Tilings of $4\times n$ floors
with $1\times 3$ tiles
without slide lines)
\begin{equation}
\hat T_0(z,4) = 1+2z^6\frac{1}{(1-z^3)(1-5z^3+3z^6-z^9)(1-4z^3+3z^6-z^9)}.
\end{equation}
\end{conj}

\appendix
\section{C++ Program}
The tables have been computed by a C++ program which is
listed in the ancillary directory. Details of the implementation
are documented in the source code. The code is compiled with
\begin{verbatim}
make
\end{verbatim}
and a suite of tests is then run with
\begin{verbatim}
make test
\end{verbatim}
A reminder of the options of the main program is  shown
by calling the main program without any arguments:
\begin{verbatim}
tatamiMain
\end{verbatim}

The program generates all available tilings with an exhaustive recursive 
placement of new tiles into the current floor, starting with the empty floor.
To support elimination of congruential tilings in case the \texttt{-i} option
has been supplied, a ordering of tilings is defined by 
mapping each tiling onto a binary sequence, then admitting/counting only
those tilings that are represented by the smallest binary number amongst
the 4 or 8 roto-reflected congruential copies.

\bibliographystyle{amsplain}
\bibliography{all}

\end{document}